\documentclass[a4paper,12pt]{article}
 \usepackage[a4paper,top=1cm,bottom=1.25cm,left=2.5cm,right=2cm,marginparwidth=0.75cm]{geometry}
\usepackage{amsmath}
\usepackage{amssymb}
\usepackage{mathtools} %loads amsmath as well
\usepackage{enumerate}% http://ctan.org/pkg/enumerate\jmath{{{{k}}}}
\usepackage{graphicx}
\usepackage{float}
\usepackage{cite}
\usepackage{color}
\usepackage{url}
\usepackage{subcaption}
\usepackage{caption}
\usepackage{hyperref}
\usepackage{amssymb}
\usepackage{amsmath}
\usepackage{amsthm}
\usepackage{float}
\usepackage{siunitx}
\usepackage{booktabs}
\usepackage[english]{babel}
\usepackage[nottoc]{tocbibind}
\usepackage{lipsum}
\usepackage{booktabs}
\usepackage{siunitx}
\usepackage{makecell}

\usepackage{changepage}
\usepackage{orcidlink}
\usepackage{hyperref}
\definecolor{DarkBlue}{rgb}{0.00,0.00,0.50}
\hypersetup{colorlinks=true,
            linkcolor=blue,
            filecolor=magenta,
            urlcolor=cyan,
            citecolor=DarkBlue,
            pdftitle={MSc UoR Computer Science Report},
            }
\usepackage{url}
\usepackage{caption}
\usepackage{array}
\usepackage{rotating}
\usepackage{algorithm,algorithmic}
\usepackage{amsthm}
\newtheorem{definition}{Definition}
\let\olddefinition\definition
\renewcommand{\definition}{\olddefinition\normalfont}

\newtheorem{prb}{Problem}
\let\oldprb\prb
\renewcommand{\prb}{\oldprb\normalfont}
\newtheorem{lemma}{Lemma}
\let\oldlemma\lemma
\renewcommand{\lemma}{\oldlemma\normalfont}
\newtheorem{pol}{Polynomial}
\let\oldpol\pol
\renewcommand{\pol}{\oldpol\normalfont}

\usepackage{authblk}
\usepackage{blindtext}
\usepackage{datetime}
\newdateformat{monthdayyeardate}{%
  \monthname[\THEMONTH]~\THEDAY, \THEYEAR}
\begin{document}
%\date{ }
 %\today
 \monthdayyeardate
\title{Comparative Analysis of Polynomials with Their Computational Costs}

\author{ Qasim Khan}
\author{   Anthony Suen}
\affil{{ Department of Mathematics and Information Technology,  The Education University of Hong Kong, 10 Lo Ping Road, Tai Po, N.T, Hong} Kong}
\affil{\textit {qasimkhan@s.eduhk.hk  \ \
acksuen@eduhk.hk}}
\maketitle

\begin{abstract}

In this article, we explore the effectiveness of two polynomial methods in solving non-linear time and space fractional partial differential equations. We first outline the general methodology and then apply it to five distinct experiments. The proposed method, noted for its simplicity, demonstrates a high degree of accuracy. Comparative analysis with existing techniques reveals that our approach yields more precise solutions. The results, presented through graphs and tables, indicate that He's and Daftardar-Jafari polynomials significantly enhance accuracy. Additionally, we provide an in-depth discussion on the computational costs associated with these polynomials. Due to its straightforward implementation, proposed  method can  be extended for application to a broader range of problems.
\end{abstract}

\section{Introduction}
Fractional Calculus (FC) is a mathematical field that extends the traditional concepts of integrals and derivatives to non-integer orders. This area of study has gained significant attention due to its wide-ranging applications in modeling various real-world phenomena. For instance, FC has been effectively applied to model epidemic spreading, providing insights into how diseases propagate through populations (see \cite{costRw1}). Additionally, it has been utilized in the analysis of earthquakes, offering a more nuanced understanding of seismic activities. In the biomedical and biological fields, fractional calculus has facilitated the development of models that capture complex biological processes and interactions. The growing interest in FC has led researchers to formulate sophisticated mathematical models that incorporate fractional integrals and derivatives. These models often address complex systems where traditional integer-order derivatives fall short. To analyze and solve these fractional models, researchers have employed innovative techniques from various disciplines, including robotic technology, genetic algorithms, and applications in physics, economics, and finance (see \cite{cost7,cost8,cost9}). As the field has evolved, a notable advancement is the development of fractional partial differential equations (FPDEs). These equations enhance modeling precision in several domains, including mechanics, plasma physics, finance, biomathematics, and fluid mechanics. The versatility of FPDEs has made them invaluable in applied sciences, allowing for more accurate representations of dynamic systems (see \cite{cost6}).

Fractional differential equations (FDEs) have emerged as powerful tools for analyzing and modeling the behavior of various physical phenomena. Their capability to capture complex dynamics makes them particularly valuable in contexts where traditional integer-order differential equations may fall short. This increased accuracy is well-documented in the literature, where FDEs are frequently employed to provide more precise models of diverse phenomena (see Section 2, Reference \cite{costa1}). One notable application of FDEs is in the study of nonlinear oscillations. Researchers have found that incorporating fractional derivatives leads to a better understanding of oscillatory behavior in systems that exhibit memory or hereditary properties (see \cite{cost5}, \cite{cost6}). Similarly, the non-linear Korteweg-De Vries (KdV) equation, which models wave propagation in shallow water, has been effectively adapted to fractional forms, enhancing its predictive capabilities (see \cite{costt1}). FDEs have also been applied to model advection-diffusion processes in two-dimensional semiconductor systems, where the inclusion of fractional derivatives allows for a more nuanced representation of charge carrier dynamics (see \cite{costt2}). Additionally, these equations have been utilized to analyze the dynamics of national soccer leagues, where they help in understanding complex interactions and patterns within the league's performance metrics. Moreover, space-fractional diffusion processes have been extensively studied using FDEs, providing insights into phenomena such as anomalous diffusion, where particles spread in a non-standard manner across a medium (see \cite{costt3,cost1,cost2,cost3,cost4}). The ability of FDEs to model these intricate behaviors highlights their significance in both theoretical and applied research across various scientific fields.

Recently, mathematicians have shown a growing interest in developing and applying a variety of efficient numerical and analytical techniques to solve fractional partial differential equations (FPDEs) and their systems. This surge in interest reflects the complexity and significance of FPDEs in modeling real-world phenomena. Several well-established techniques have been created or implemented for this purpose. Among them, the Laplace Adomian Decomposition Method (LADM) \cite{cost54} and the Adomian Decomposition Method (ADM) are notable for their effectiveness in breaking down complex equations into simpler components. The Finite Difference Method \cite{costr3} is widely used for its straightforward approach to discretizing differential equations, making it suitable for numerical solutions. Other innovative methods include Feng’s First Integral Method \cite{costr4}, which offers a unique perspective in solving FPDEs, and the Homotopy Analysis Method (HAM) \cite{costr24}, which provides a systematic way to construct solutions. The Homotopy Perturbation Transform Method (HPTM) \cite{costr25} and Meshless Method (MM) \cite{costr27} are also significant for their flexibility in handling various types of equations. Moreover, the Aboodh Transform Iterative Method \cite{costRw2} and the Modified Homotopy Perturbation Method (MHPM) \cite{costr26} have been developed to enhance the accuracy of solutions. Techniques such as the Multiple Exponential Function Algorithms \cite{cost23}, Shifted Chebyshev-Gauss-Lobatto Collocation \cite{costZ1}, and the Operational Matrix Method \cite{costZ2} further exemplify the breadth of methods available for tackling FPDEs. The Variational Iteration Method (VIM) \cite{costvim} also plays a crucial role in providing approximate solutions to these complex equations. Overall, while obtaining analytical and approximate solutions for non-linear FPDEs and their systems can be a challenging endeavor, it is a fascinating area of research that attracts mathematicians. Researchers are particularly keen to discover accurate and effective techniques for solving FPDEs, as the analysis of fractional solutions significantly contributes to understanding the actual dynamics of various physical phenomena. This growing interest in FPDEs has led to increased popularity and numerous studies within the mathematical community \cite{cost11,cost12,cost13,cost14,cost15,cost16,cost17,cost18}.

Controlling the non-linear components of differential equations remains a significant challenge for researchers. Non-linearities often introduce complex behaviors that can complicate both analysis and solution processes. These challenges stem from the intricate interactions within non-linear systems, which can lead to phenomena such as chaos, bifurcations, and multiple equilibria. As a result, developing effective methods to manage and control these non-linear parts is crucial for advancing our understanding and application of differential equations in various fields. Researchers continue to explore innovative approaches and techniques to tackle these complexities, striving for more robust solutions in both theoretical and practical contexts. In \cite{costg1}, one said that most of the non-linear systems are impossible to solve analytically.
In solving non-linear differential equations, the decomposition and iterative handling of non-linear terms are crucial for achieving accurate and efficient solutions. Various polynomial methods have been developed to address this challenge, including Daftardar-Jafari polynomials, Adomian polynomials, and He’s polynomials. Each of these methods offers a unique approach to decomposing non-linear terms, facilitating the iterative process required for solving complex differential equations. Here, we provide an overview of these polynomial techniques and their applications. Daftardar-Jafari polynomials are employed in iterative methods to break down non-linear terms systematically. This approach enhances the efficiency of solving non-linear functional equations. The polynomials are defined recursively, allowing for the iterative approximation of solutions (see ref \cite{cost105,costDJ2}).
Adomian polynomials are a cornerstone of the Adomian Decomposition Method (ADM), which is extensively used for solving non-linear differential equations. These polynomials decompose the non-linear operator into a series, simplifying the solution process (see ref. \cite{costAD2,costAD3}).
He’s polynomials are integral to variational iteration methods and other iterative techniques. They facilitate the management of non-linear terms, improving the convergence and accuracy of the solutions for non-linear differential equations (see ref. \cite{costhes,cost1,cost3,cost5}).
The polynomial methods of Daftardar-Jafari, Adomian, and He's provide robust frameworks for handling non-linear terms in differential equations. Their recursive and systematic approaches enable iterative solutions that are both accurate and computationally efficient. These methods are foundational tools in the mathematical toolbox for solving a wide range of non-linear fractional order integro-differential equations and non-linear fractional order partial differential equations, with applications spanning various scientific and engineering disciplines (see, recent results in \cite{costqasim1,costqasim2,costqasim3,costqasim4,costqasim5}).

 This research article tackles non-linear fractional-order PDEs using the Iterative Laplace Transform Method (ILTM) along with He's and D-J polynomials. This approach is notable for its higher accuracy and robustness when applied to fractional differential equations. The study carefully examines the effects of these polynomials and presents the results using graphs and tables. ILTM is an advanced technique that iteratively applies the Laplace transform to solve differential equations. The iterative nature of this method allows for improved convergence and accuracy, making it suitable for non-linear and fractional-order equations. By transforming the problem into the Laplace domain, the differential equation is converted into an algebraic equation, which is often simpler to solve. He's polynomials and D-J polynomials are utilized to approximate the solutions of the transformed equations. These polynomials are instrumental in breaking down complex non-linear terms into manageable forms, facilitating the iterative process of ILTM. The choice of polynomials can significantly influence the convergence and accuracy of the solution, which is why their effects are closely examined in the study. The solutions obtained through ILTM are analyzed using graphs and tables. Graphical analysis helps visually assess the solution's behavior over time and space, while tabular data provides quantitative insights into the accuracy and convergence rates. This dual approach ensures a comprehensive evaluation of the method's performance. The success of ILTM in solving the fractional-order Advection equation paves the way for its application to other fractional problems. For instance, fractional diffusion equations, fractional wave equations, and fractional Schrödinger equations could benefit from this approach. The methodology can be adapted by applying similar polynomial approximations and iterative procedures under various integral transforms. This research shows a strong and accurate method for solving non-linear fractional-order PDEs using ILTM and polynomial approximations. The method's potential to be extended to other fractional problems highlights its importance and opens up many opportunities for future research and applications.

The rest of the paper is organised as follows. Section \ref{sec2}
 presents an overview of fundamental definitions and the methodology utilized in the study. In Section \ref{sec3}, numerical problems are addressed to illustrate the method's applications. Section \ref{result}
 discusses the obtained results, providing detailed analysis and discussion. Finally, Section \ref{conclusion}
 concludes the paper by summarizing the main findings and proposing directions for future research.
\section{Preliminaries Concept and Methodology }\label{sec2}
This section is divided into two subsections. In section \ref{subsec1}, we discuss the essential preliminaries relevant to the current topic, which are crucial for successfully completing the research task at hand. These foundational concepts provide the necessary background and context, enabling a deeper understanding of the methodologies and techniques employed in our study. Section \ref{subsec2}  builds upon the foundational concepts outlined in Section \ref{subsec1}. Here, we delve into the specific methodologies used in our research, explaining how they are applied to fractional PDEs.
\subsection{Basic Definitions}\label{subsec1}

\begin{definition}
A bunch of literature discusses various definitions of fractional derivative operators, including the Riemann-Liouville, Caputo, Caputo-Fabrizio, and Atangana-Baleanu fractional derivatives. This paper focuses on one of the most widely used derivative operators in fractional calculus, specifically the Caputo derivative operator. Introduced by Michele Caputo in 1967 (as referenced in \cite{cost43}), the Caputo operator for fractional derivatives of order $\alpha$ is defined as follows:
\begin{align*}
 D_{{{t}}}^{\alpha}{{{u}}}({{{t}}})= \frac{1}{\Gamma(n-\alpha)} \int_0^{{{t}}}({{{t}}}-{{{t}}}_0)^{n-\alpha-1} f^{(n)}({{{t}}}_0)d{{{t}}}_0,
\end{align*}
where $n-1<\alpha\leq n+1, \  \ n\in \mathbb{N}$ and ${{{t}}}>0$. Throughout this paper, the parameter $\alpha$ represents the fractional order and $D_{{t}}^{\alpha}$ denotes the Caputo derivative operator.
\end{definition}
\begin{definition}
The Laplace transform is a powerful integral transform used to convert a function of time, $g({{{t}}})$ , into a function of a complex variable $s$. It is particularly useful for solving differential equations and analyzing linear time-invariant systems. The Laplace transform ${\mathcal{L}}[g({{{t}}})]$ of a function
$g({{{t}}})$ is defined as \cite{cost6}:
\begin{equation*}
G(s)={\mathcal{L}}[g({{{t}}})]=\int_0^\infty e^{-s\chi}g({{{t}}})d{{{t}}}.
\end{equation*}
If ${\mathcal{L}}[f({{{t}}})]$ and ${\mathcal{L}}[g({{{t}}})]$ exist, then   the following properties hold:
\begin{equation*}
{\mathcal{L}}[f({{{t}}})+g({{{t}}})]={\mathcal{L}}[f({{{t}}})]+{\mathcal{L}}[g({{{t}}})], \ \ \ {\mathcal{L}}[cf({{{t}}})]=c{\mathcal{L}}[f({{{t}}})],
\end{equation*}
where $c$ is a constant.
\end{definition}
\begin{definition}
The Laplace transform of Caputo operator is given as \cite{cost43}
\begin{equation*}
\begin{split}
    {\mathcal{L}} (D_{{{t}}}^{\alpha}{{{u}}}({{{t}}}))&=s^\alpha {\mathcal{L}}[{{{u}}}({{{t}}})]-\sum_{k=0}^{n-1}s^{n-k-1}U^k(0^+), \ \ \ n-1<\alpha\leq n,\\
 {\mathcal{L}} (D_{{{t}}}^{\alpha}{{{u}}}({{{t}}}))&=s^\alpha U(s)-s^{\alpha-1}U(0).
\end{split}
\end{equation*}

\end{definition}
\begin{definition}\label{D-J}
\textnormal{   Daftardar-Gejji and Jafari polynomials, often abbreviated as D-J polynomials (see ref. \cite{cost105}), are utilized to represent the non-linear term in a given problem and are defined as follows: }\\
\begin{equation}\label{s4}
N \left(\sum_{{j}=0}^\infty{{{u}}_{j}({{{t}}})}\right)=N({{{u}}_0({{{t}}})}+\sum_{{j}=0}^\infty \bigg[{N} \bigg(\sum_{i=0}^{{j}} {{u}}_i({{{t}}})\bigg)-{N} \bigg(\sum_{i=0}^{{j}-1 }{{u}}_i({{{t}}})\bigg)\bigg].
\end{equation}
\end{definition}

\begin{definition}\label{hes}
	\textnormal{
The  He’s polynomials  denoted by $\mathbf{H} $ to express non-linear terms in a given problem is defined  as\cite{costhes}}
\begin{equation}\label{40}
N{u( {{{t}}})}=\sum_{{{j}=0}}^{\infty}\mathbf{H}_{{j}},
\end{equation}
where $\mathbf{H}_{{{j}}}$ is given by
\begin{equation*}
\mathbf{H}_{{{j}}}({{{{u}}}}_0,{{{{u}}}}_1,\cdots{{{{u}}}}_{{{j}}})=\frac{1}{{{j}}!}\left[\frac{d^{{j}}}{d{{{c_2}}}^{{j}}}\left\{N\sum_{{j}=0}^{\infty}({{{c_2}}}^{{j}}u_{j})\right\}\right]_{{{{c_2}}}=0},
\qquad {j}=0,1,\cdots.
\end{equation*}
\end{definition}
\begin{definition}
  The Mittag-Leffler function is a special function that generalizes the exponential function and plays an important role in various areas of mathematical analysis, including fractional calculus, complex analysis, and the theory of differential equations. It was introduced by the Swedish mathematician Gösta Mittag-Leffler in \cite{costmit39}, defined as
  \begin{equation*}
			{E}_{{\alpha}}({{{t}}})=\sum_{n=0}^{\infty}\frac{{{{t}}}^{n}}{\Gamma
				({{\alpha}}{n}+1)},\ \ \	\textnormal{ where} \ \ {{{t}}} \ \in \ {\mathbb{C}}{.}
		\end{equation*}
\end{definition}

We include the following lemma which contains some formulae for Laplace transform. The proof can be found in \cite{costl1,costl2,costl3,costl4,costl6}.

\begin{lemma}\label{l1}

Let $f({{{t}}})$ and $g({{{t}}})$ be two piecewise continuous functions defined on the interval $[0, \infty)$, and these functions are of exponential order.
Let $F(s) = {\mathcal{L}}[f({{{t}}})]$ and $G(s) = {\mathcal{L}}[g({{{t}}})]$, where ${\mathcal{L}}$ denotes the Laplace transform.
Additionally, let ${{{c_2}}}$, ${{{c_1}}}$, $\rho$, and $\zeta$ be constant parameters.
Under these conditions, the following properties (1-8) are satisfied:
\begin{enumerate}
\item	${\mathcal{L}}[{{{c_1}}} f({{{{t}}}})+{{{c_2}}} g({{{{t}}}})]={{{c_1}}} F(s)+{{{c_2}}} G(s)$.
\item	${\mathcal{L}}^{-1}[{{{c_1}}} F({s})+{{{c_2}}} G({s})]={{{c_1}}} f({{{{t}}}})+{{{c_2}}} g({{{{t}}}}).$
\item	${\mathcal{L}}[e^{\rho{{{{t}}}}}f({{{{t}}}})]=F(s-\rho).$
\item ${\mathcal{L}}[ f(\zeta{{{{t}}}})]=\frac{1}{\zeta}F(\frac{s}{\zeta}), \ \ \zeta>0.$
\item $\lim_{s\rightarrow\infty}sF(s)=f(0).$
\item     ${\mathcal{L}}[J_{{{{t}}}}^{\alpha}f({{{{t}}}})]=\frac{F({s})}{s^{\alpha}}, \ \ \alpha>0.$
\item ${\mathcal{L}}[D_{{{{t}}}}^{\alpha}f({{{{t}}}})]=s^{\alpha}F(s)-\sum_{k=0}^{m-1}s^{\alpha-k-1}f^{(k)}(0), \ \ m-1<\alpha\leq m.$
\item ${\mathcal{L}}[D_{{{{t}}}}^{n\alpha}f({{{{t}}}})]=s^{n\alpha}F(s)-\sum_{k=0}^{n-1}s^{(n-k)\alpha-1}\left(D_{{{{t}}}}^{k\alpha}f\right)(0), \ \ 0<\alpha\leq1.$
\end{enumerate}
\end{lemma}
\subsection{Research Methodology}\label{subsec2}

To understand the fundamental methodology of the Iterative Laplace Transform Method (ILTM) with two different polynomials, we will consider the following two distinct cases: In \ref{pol1}, we will apply the ILTM using He's Polynomial, examining its specific characteristics and how it influences the iterative process. in the same way, we use D-J polynomials in \ref{pol2}.

Let us consider a general non-homogenous FPDEs of the form,
\begin{equation}\label{a1}
D_{{{t}}}^{\alpha+m}{{{u}}}({{{{x}}}},{{{t}}})=f({{{{x}}}},{{{t}}})+{{\mathcal{R}}} {{{u}}}({{{{x}}}},{{{t}}})+\mathcal{N} {{{u}}}({{{{x}}}},{{{t}}}), \ \ \ \ m-1<\alpha\leq m, \ \ m\in \mathbb{N},
\end{equation}
with initial condition
\begin{equation*}
\frac{\partial^k {{{u}}}({{{{x}}}},0)}{\partial {{{t}}}^k}={{{t}}}_k({{{{x}}}}), \ \ \ \ k=0,1,\cdots,n-1.
\end{equation*}
In Eq. (\ref{a1}), the term ${{\mathcal{R}}}$ denotes a linear operator, and $\mathcal{N}$ represents a nonlinear operator. The expression $f({{{{x}}}},{{{t}}})$ is a source term, where ${{{x}}}$ is a variable and ${{{t}}}$ is a parameter. Applying Laplace transform to Eq. (\ref{a1}), we obtain
\begin{equation}\label{a2}
{\mathcal{L}}\bigg({{{u}}}({{{{x}}}},{{{t}}})\bigg)=\vartheta({{{{x}}}},s)+\bigg(\frac{1}{s^{\alpha+m}}\bigg){\mathcal{L}}\bigg({{\mathcal{R}}} {{{u}}}({{{{x}}}},{{{t}}})+\mathcal{N} {{{u}}}({{{{x}}}},{{{t}}})\bigg),
\end{equation}
where
\begin{equation*}
\vartheta({{{{x}}}},s)=\frac{1}{s^{m+1}}\bigg(s^m{{{t}}}_0({{{{x}}}})+{s^{m-1}}{{{t}}}_1({{{{x}}}})+...
+{{{t}}}_m({{{{x}}}})\bigg)+\frac{1}{s^{\alpha+m}}f({{{{x}}}},s).
\end{equation*}
Taking the inverse Laplace Transform on Eq. (\ref{a2}), we get
\begin{equation}\label{a3}
{{{u}}}({{{{x}}}},{{{t}}})=\vartheta({{{{x}}}},{{{t}}})+{\mathcal{L}}^{-1}\bigg(\frac{1}{s^{\alpha+m}}\bigg){\mathcal{L}}\bigg({{\mathcal{R}}} {{{u}}}({{{{x}}}},{{{t}}})+\mathcal{N} {{{u}}}({{{{x}}}},{{{t}}})\bigg).
\end{equation}
The NIM procedure  is  implemented as
\begin{equation}\label{a4}
{{{u}}}({{{{x}}}},{{{t}}})=\sum_{{{{{k}}}}=0}^\infty {{{u}}}_{{{{k}}}}({{{{x}}}},{{{t}}}).
\end{equation}
Since ${{\mathcal{R}}}$ is linear
\begin{equation}\label{a5}
{{\mathcal{R}}} \bigg(\sum_{{{{{k}}}}=0}^\infty {{{u}}}_{{{{k}}}}({{{{x}}}},{{{t}}})\bigg)=\sum_{{{{{k}}}}=0}^\infty {{\mathcal{R}}} \bigg({{{u}}}_{{{{k}}}}({{{{x}}}},{{{t}}})\bigg).
\end{equation}
The nonlinear term $\mathcal{N}$ has been controlled, using two well-known  polynomials,

\begin{pol}\label{pol1} Recall, He's polynomials defined in definition \ref{hes} and use it for the non-linear term,
in the view of  Eqs. (\ref{a2}-\ref{a4}) and Eq. (\ref{a5}).
\begin{equation*}
\sum_{{{{{k}}}}=0}^\infty {{{u}}}_{{{{k}}}}({{{{x}}}},{{{t}}})=\vartheta({{{{x}}}},{{{t}}})+{\mathcal{L}}^{-1}\bigg[\bigg(\frac{1}{s^{\alpha}}\bigg){\mathcal{L}}\bigg(\sum_{{{{{k}}}}=0}^\infty {{\mathcal{R}}} {{{u}}}_{{{{k}}}}({{{{x}}}},{{{t}}})\bigg)\bigg]+{\mathcal{L}}^{-1}\bigg[\bigg(\frac{1}{s^{\alpha}}\mathcal{L}(\mathbf{H}_{{{k}}}({{{{u}}}}_0,{{{{u}}}}_1,\cdots{{{{u}}}}_{{{k}}}))\bigg) \bigg].
\end{equation*}
The ILTM algorithm is
\begin{equation*}
{{{u}}}_0({{{{x}}}},{{{t}}})=\vartheta({{{{x}}}},{{{t}}}),\qquad {{{u}}}_1 ({{{{x}}}},{{{t}}})={\mathcal{L}}^{-1}\bigg[\bigg(\frac{1}{s^{\alpha}}\bigg){\mathcal{L}}\bigg({{\mathcal{R}}}({{{u}}}_0({{{{x}}}},{{{t}}}))+\mathbf{H}_0({{{u}}}_0({{{{x}}}},{{{t}}}))\bigg)\bigg],
\end{equation*}
\ \ \ \ \ \ \ \ \ \ \ \ \ \ \ \ \ \ \ \ \ \ \ \ \ \ \ \ \ \  \ \ \ \ \ \ \ \ \ \ \ \ \ \ \ \ \ \vdots
\begin{equation*}
{{{u}}}_{{{{{k}}}}+1}({{{{x}}}},{{{t}}})= {\mathcal{L}}^{-1}\bigg[\bigg(\frac{1}{s^{\alpha}}\bigg)
{\mathcal{L}}\bigg({{\mathcal{R}}}({{{u}}}_{{{{{k}}}}}({{{{x}}}},{{{t}}}))
+\mathbf{H}_{{{k}}}({{{{}}}}{{{{u}}}}_1,\cdots{{{{u}}}}_{{{k}}})) \bigg)\bigg].
\end{equation*}

\begin{pol}\label{pol2} Recall the  D-J polynomials defined in definition \ref{D-J} and we use it for the non-linear term
 \begin{equation*}
\mathcal{N} \bigg(\sum_{{{{{k}}}}=0}^\infty {{{u}}}_{{{{k}}}}({{{{x}}}},{{{t}}})\bigg)=\mathcal{N} ({{{u}}}_0({{{{x}}}},{{{t}}}))+\sum_{{{{{k}}}}=0}^\infty \bigg[\mathcal{N} \bigg(\sum_{i=0}^{{{{k}}}} {{{u}}}_i({{{{x}}}},{{{t}}})\bigg)-\mathcal{N} \bigg(\sum_{i=0}^{{{{{k}}}}-1}{{{u}}}_i({{{{x}}}},{{{t}}})\bigg)\bigg].
\end{equation*}
In view of  Eqs. (\ref{a2}-\ref{a4}) and Eq. (\ref{a5}), we have
\begin{multline*}
\sum_{{{{{k}}}}=0}^\infty {{{u}}}_{{{{k}}}}({{{{x}}}},{{{t}}})=\vartheta({{{{x}}}},{{{t}}})+{\mathcal{L}}^{-1}\bigg[\bigg(\frac{1}{s^{\alpha}}\bigg){\mathcal{L}}\bigg(\sum_{{{{{k}}}}=0}^\infty {{\mathcal{R}}} {{{u}}}_{{{{k}}}}({{{{x}}}},{{{t}}})\bigg)\bigg]+{\mathcal{L}}^{-1}\bigg[\bigg(\frac{1}{s^{\alpha}}\bigg)\\
{\mathcal{L}}\bigg(\mathcal{N}({{{u}}}_0({{{{x}}}},{{{t}}}))\bigg)
+\sum_{{{{{k}}}}=0}^\infty \bigg\{\mathcal{N}(\sum_{i=0}^{{{{k}}}} {{{u}}}_i({{{{x}}}},{{{t}}})\bigg)-\mathcal{N} \bigg(\sum_{i=0}^{{{{{k}}}}-1} {{{u}}}_i({{{{x}}}},{{{t}}})\bigg)\bigg\}\bigg].
\end{multline*}
The ILTM algorithm is given by
\begin{equation*}
{{{u}}}_0({{{{x}}}},{{{t}}})=\vartheta({{{{x}}}},{{{t}}}),\qquad {{{u}}}_1 ({{{{x}}}},{{{t}}})={\mathcal{L}}^{-1}\bigg[\bigg(\frac{1}{s^{\alpha}}\bigg){\mathcal{L}}\bigg({{\mathcal{R}}}({{{u}}}_0({{{{x}}}},{{{t}}}))+\mathcal{N}({{{u}}}_0({{{{x}}}},{{{t}}}))\bigg)\bigg],
\end{equation*}
\ \ \ \ \ \ \ \ \ \ \ \ \ \ \ \ \ \ \ \ \ \ \ \ \ \ \ \ \ \  \ \ \ \ \ \ \ \ \ \ \ \ \ \ \ \ \ \vdots
\begin{equation*}
{{{u}}}_{{{{{k}}}}+1}({{{{x}}}},{{{t}}})= {\mathcal{L}}^{-1}\bigg[\bigg(\frac{1}{s^{\alpha}}\bigg)
{\mathcal{L}}\bigg({{\mathcal{R}}}({{{u}}}_{{{{{k}}}}}({{{{x}}}},{{{t}}}))
+\sum_{j=1}^\infty \bigg\{\mathcal{N}\bigg(\sum_{i=0}^n {{{u}}}_i({{{{x}}}},{{{t}}})\bigg)-\mathcal{N} \bigg(\sum_{i=0}^{n-1} {{{u}}}_i({{{{x}}}},{{{t}}})\bigg)\bigg\}\bigg)\bigg].
\end{equation*}
\end{pol}
The ILTM series form approximate solution is then obtained by
\begin{equation*}
{{{u}}}({{{{x}}}},{{{t}}})={{{u}}}_0+{{{u}}}_1+{{{u}}}_2+...+{{{u}}}_{{{{{k}}}}}.
\end{equation*}
\end{pol}

		The algorithms for implementing He's   polynomials and Daftardar-Jafari polynomials are outlined below.
We compare Algorithm \ref{al1}  and Algorithm \ref{al2}  to assess their efficiency in handling polynomials, which will be given in Section~\ref{sec3}. Our findings reveal that Algorithm  is more effective in dealing with non-linear cases. The results and discussions, along with corresponding graphical and tabular representations, will be presented in Section \ref{result}.
		\begin{algorithm}
			\caption{ILTM steps with He's polynomials}
			\begin{algorithmic}[1]
				\STATE Apply the Laplace transform to equation \ref{a2}
				\STATE Use the properties of the Laplace transform to simplify the algebraic equation
				\STATE Use the initial guess and source term for 1st iteration
				\STATE Control the non-linear term by He's polynomials
				\STATE Use inverse Laplace transform to convert the expression from s domain to time domain
				%\STATE Express the decomposed solution as an infinite series
				%\STATE Use the recursive relations to compute the first few terms of the series solution
				\STATE Sum of the  first few terms of the series  to obtain an approximate solution
				\STATE Compare the obtained approximate solution with other existing techniques to test the accuracy
			\end{algorithmic}\label{al1}	
		\end{algorithm}
		\begin{algorithm}
			\caption{ILTM steps with Daftardar-Jafari Polynomials }
			\begin{algorithmic}[1]
				\STATE Apply the Laplace transform to equation \ref{a2}
				\STATE Use the properties of the Laplace transform to simplify the algebraic equation
				\STATE Use the initial guess and source term for 1st iteration
				\STATE Control the non-linear term by D-J polynomials
				\STATE Use inverse Laplace transform to convert the expression from s domain to time domain
				%\STATE Express the decomposed solution as an infinite series
				%\STATE Use the recursive relations to compute the first few terms of the series solution
				\STATE Sum of the  first few terms of the series  to obtain an approximate solution
				\STATE Compare the obtained approximate solution with other existing techniques to test the accuracy
			\end{algorithmic}\label{al2}	
		\end{algorithm}

\section{Experimental Results }\label{sec3}
In this section, we will conduct experiments to evaluate the effectiveness of He's and D-J polynomials in solving various fractional-order PDEs. It is important to note that the problems we are investigating have been previously explored by different  researchers. The contributions of these authors to the field are significant and should not be underestimated (refer to references \cite{costnex1,costGDTM,costADM,cost51a,cost52a,cost52b,costHPM,costAhmad,costRasool,costf1,cost2w1,cost2w2,cost2w3,cost2w4,cost2w5,cost2w6,cost2w7}, and the references inside there). Our analysis will involve a comparison of our results with those already documented in the existing literature. We are not only interested in exploring the results of these problems but also in comparing their computational costs.   Importantly, we achieve more accurate results compared to those previously documented in the literature and cited references up to date. We tested the following problems,
\begin{prb}\label{problem1}
Consider the non-linear fractional order advection equation of the form,
\begin{equation}\label{11}
{D^{\alpha}_{{{{t}}}}{{{{u}}}}}+{{{u}}}{{{{u}}}_{{{{x}}}}}={{{{{x}}}}+{{{{x}}}}{{{t}}}^2}, \;\;\;\;\;\;\;\;\;\; {{{{t}}}>0}, \;\;\;\; {0<\alpha\leq1},
\end{equation}
with initial condition
$$
{{{{u}}}({{{{x}}}},0)=0}.
$$
The exact solution for $\alpha=1$ is ${{{{{u}}}}}({{{x}}},{{{t}}})={{{x}}}{{{t}}}$.
\end{prb}

\begin{table}[H]
				\centering
				\caption{\textbf{Absolute error comparison of He's and D-J polynomials across different iterations ${{{k}}}$,  spaces ${0<{x}<1}$ and time level $0<{{{t}}}<1$  of problem \ref{problem1}.} We observed  that the D-J polynomials work more accurately as compared to He's polynomials. Of course, one can improve the accuracy of D-J polynomials by including additional terms in the series solution, but the statement needs to be validated for He's polynomials.   }\label{table1}
				\resizebox{\textwidth}{!}{
					\begin{tabular}{cccccccccccc}
						\toprule
						\multicolumn{3}{c}{\textbf{ {  AE  \underline{at ${{{k}}}=1$}} }}& \multicolumn{2}{c}{\textbf{ \underline{Absolute error at ${{{k}}}=2$} }}& \multicolumn{2}{c}{\textbf{ \underline{Absolute error at ${{{k}}}=3$} }}& \multicolumn{2}{c}{\textbf{ \underline{Absolute error at ${{{k}}}=4$} }}& \multicolumn{2}{c}{\textbf{\underline{Absolute error at ${{{k}}}=5$}}} \\
						%\multicolumn{2}{c}{\textbf{ {Polynomials} }}& $\left| u_2-u_{exact} \right|$ &  $\left| u_2-u_{exact} \right|$&  $\left| u_3-u_{exact} \right|$&   $\left| u_3-u_{exact} \right|$&  $\left| u_4-u_{exact} \right|$&  $\left| u_4-u_{exact} \right|$&  He's Polynomials&   D-J Polynomials\\
						{${{{t}}}$}& ${{x}}$& {\thead{ D-J \& He's\\ Polynomials}} &{\thead{ He's\\ Polynomials}}  & {\thead{ D-J\\ Polynomials}}&  {\thead{ He's\\ Polynomials}} & {\thead{ D-J\\ Polynomials}}& {\thead{ He's\\ Polynomials}} &  {\thead{ D-J\\ Polynomials}}& {\thead{ He's\\ Polynomials}}   &  {\thead{ D-J\\ Polynomials}}\\ %[0.1ex]
						\midrule
0.1&\num{0.1}&\num{1.33E-07}&\num{5.41E-10}&\num{3.81E-10}&\num{2.20E-12}&\num{8.47E-13}&\num{4.89E-13}&\num{1.54E-15}&\num{4.97E-13}&\num{2.37E-18}\\
0.1&\num{0.3}&\num{4.00E-07}&\num{1.62E-09}&\num{1.14E-09}&\num{6.59E-12}&\num{2.54E-12}&\num{1.47E-12}&\num{4.62E-15}&\num{1.49E-12}&\num{7.11E-18}\\
0.1&\num{0.5}&\num{6.67E-07}&\num{2.71E-09}&\num{1.91E-09}&\num{1.10E-11}&\num{4.24E-12}&\num{2.44E-12}&\num{7.70E-15}&\num{2.48E-12}&\num{1.18E-17}\\
0.1&\num{0.7}&\num{9.34E-07}&\num{3.79E-09}&\num{2.67E-09}&\num{1.54E-11}&\num{5.93E-12}&\num{3.42E-12}&\num{1.08E-14}&\num{3.48E-12}&\num{1.66E-17}\\
0.1&\num{0.9}&\num{1.20E-06}&\num{4.87E-09}&\num{3.43E-09}&\num{1.98E-11}&\num{7.62E-12}&\num{4.40E-12}&\num{1.39E-14}&\num{4.47E-12}&\num{2.13E-17}\\
0.3&\num{0.1}&\num{3.27E-05}&\num{1.21E-06}&\num{8.40E-07}&\num{4.45E-08}&\num{1.68E-08}&\num{8.76E-09}&\num{2.74E-10}&\num{1.03E-08}&\num{3.79E-12}\\
0.3&\num{0.3}&\num{9.82E-05}&\num{3.62E-06}&\num{2.52E-06}&\num{1.34E-07}&\num{5.03E-08}&\num{2.63E-08}&\num{8.23E-10}&\num{3.09E-08}&\num{1.14E-11}\\
0.3&\num{0.5}&\num{1.64E-04}&\num{6.03E-06}&\num{4.20E-06}&\num{2.23E-07}&\num{8.39E-08}&\num{4.38E-08}&\num{1.37E-09}&\num{5.14E-08}&\num{1.90E-11}\\
0.3&\num{0.7}&\num{2.29E-04}&\num{8.45E-06}&\num{5.88E-06}&\num{3.12E-07}&\num{1.17E-07}&\num{6.13E-08}&\num{1.92E-09}&\num{7.20E-08}&\num{2.66E-11}\\
0.3&\num{0.9}&\num{2.95E-04}&\num{1.09E-05}&\num{7.56E-06}&\num{4.01E-07}&\num{1.51E-07}&\num{7.88E-08}&\num{2.47E-09}&\num{9.26E-08}&\num{3.41E-11}\\
0.5&\num{0.1}&\num{4.29E-04}&\num{4.48E-05}&\num{3.04E-05}&\num{4.69E-06}&\num{1.68E-06}&\num{6.72E-07}&\num{7.63E-08}&\num{1.16E-06}&\num{2.93E-09}\\
0.5&\num{0.3}&\num{1.29E-03}&\num{1.34E-04}&\num{9.11E-05}&\num{1.41E-05}&\num{5.05E-06}&\num{2.01E-06}&\num{2.29E-07}&\num{3.48E-06}&\num{8.78E-09}\\
0.5&\num{0.5}&\num{2.15E-03}&\num{2.24E-04}&\num{1.52E-04}&\num{2.35E-05}&\num{8.41E-06}&\num{3.36E-06}&\num{3.81E-07}&\num{5.79E-06}&\num{1.46E-08}\\
0.5&\num{0.7}&\num{3.00E-03}&\num{3.14E-04}&\num{2.13E-04}&\num{3.28E-05}&\num{1.18E-05}&\num{4.70E-06}&\num{5.34E-07}&\num{8.11E-06}&\num{2.05E-08}\\
0.5&\num{0.9}&\num{3.86E-03}&\num{4.03E-04}&\num{2.73E-04}&\num{4.22E-05}&\num{1.51E-05}&\num{6.04E-06}&\num{6.86E-07}&\num{1.04E-05}&\num{2.64E-08}\\
0.7&\num{0.1}&\num{2.37E-03}&\num{5.00E-04}&\num{3.24E-04}&\num{1.06E-04}&\num{3.52E-05}&\num{6.23E-06}&\num{3.12E-06}&\num{3.06E-05}&\num{2.35E-07}\\
0.7&\num{0.3}&\num{7.11E-03}&\num{1.50E-03}&\num{9.73E-04}&\num{3.17E-04}&\num{1.06E-04}&\num{1.87E-05}&\num{9.37E-06}&\num{9.17E-05}&\num{7.05E-07}\\
0.7&\num{0.5}&\num{1.19E-02}&\num{2.50E-03}&\num{1.62E-03}&\num{5.29E-04}&\num{1.76E-04}&\num{3.12E-05}&\num{1.56E-05}&\num{1.53E-04}&\num{1.17E-06}\\
0.7&\num{0.7}&\num{1.66E-02}&\num{3.50E-03}&\num{2.27E-03}&\num{7.41E-04}&\num{2.47E-04}&\num{4.36E-05}&\num{2.19E-05}&\num{2.14E-04}&\num{1.64E-06}\\
0.7&\num{0.9}&\num{2.13E-02}&\num{4.50E-03}&\num{2.92E-03}&\num{9.52E-04}&\num{3.17E-04}&\num{5.61E-05}&\num{2.81E-05}&\num{2.75E-04}&\num{2.11E-06}\\
0.9&\num{0.1}&\num{8.63E-03}&\num{3.13E-03}&\num{1.90E-03}&\num{1.14E-03}&\num{3.43E-04}&\num{6.95E-05}&\num{5.02E-05}&\num{4.37E-04}&\num{6.23E-06}\\
0.9&\num{0.3}&\num{2.59E-02}&\num{9.39E-03}&\num{5.70E-03}&\num{3.41E-03}&\num{1.03E-03}&\num{2.08E-04}&\num{1.51E-04}&\num{1.31E-03}&\num{1.87E-05}\\
0.9&\num{0.5}&\num{4.32E-02}&\num{1.57E-02}&\num{9.50E-03}&\num{5.69E-03}&\num{1.72E-03}&\num{3.47E-04}&\num{2.51E-04}&\num{2.19E-03}&\num{3.12E-05}\\
0.9&\num{0.7}&\num{6.04E-02}&\num{2.19E-02}&\num{1.33E-02}&\num{7.97E-03}&\num{2.40E-03}&\num{4.86E-04}&\num{3.51E-04}&\num{3.06E-03}&\num{4.36E-05}\\
0.9&\num{0.9}&\num{7.77E-02}&\num{2.82E-02}&\num{1.71E-02}&\num{1.02E-02}&\num{3.09E-03}&\num{6.25E-04}&\num{4.52E-04}&\num{3.94E-03}&\num{5.61E-05}\\
						\bottomrule
				\end{tabular}}
			\end{table}

  \begin{figure}[H]
    \centering
    \includegraphics[width=5.3cm]{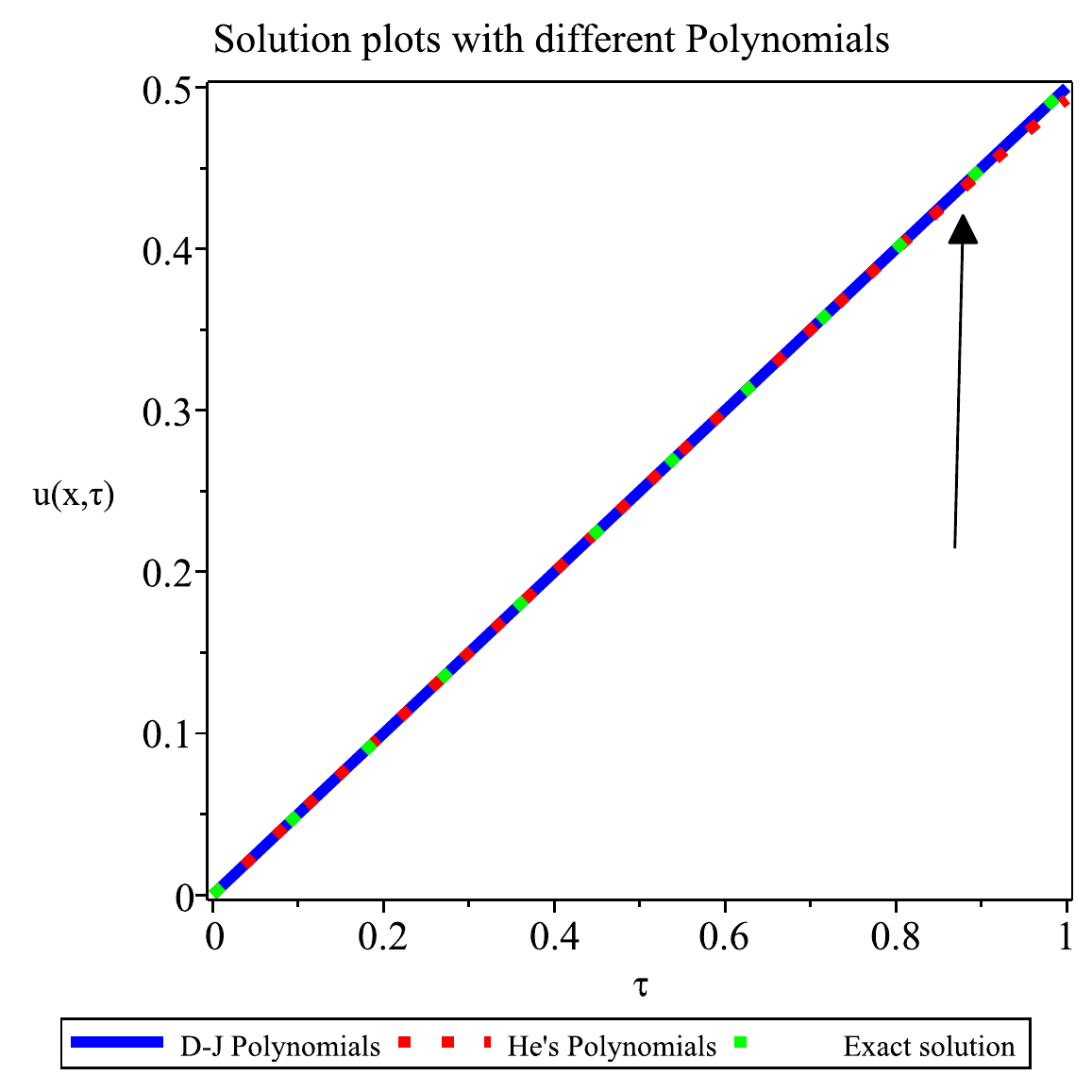}%
    \llap{\raisebox{1.0cm}{%  move next graphics to top right corner
      \includegraphics[width=2.0cm]{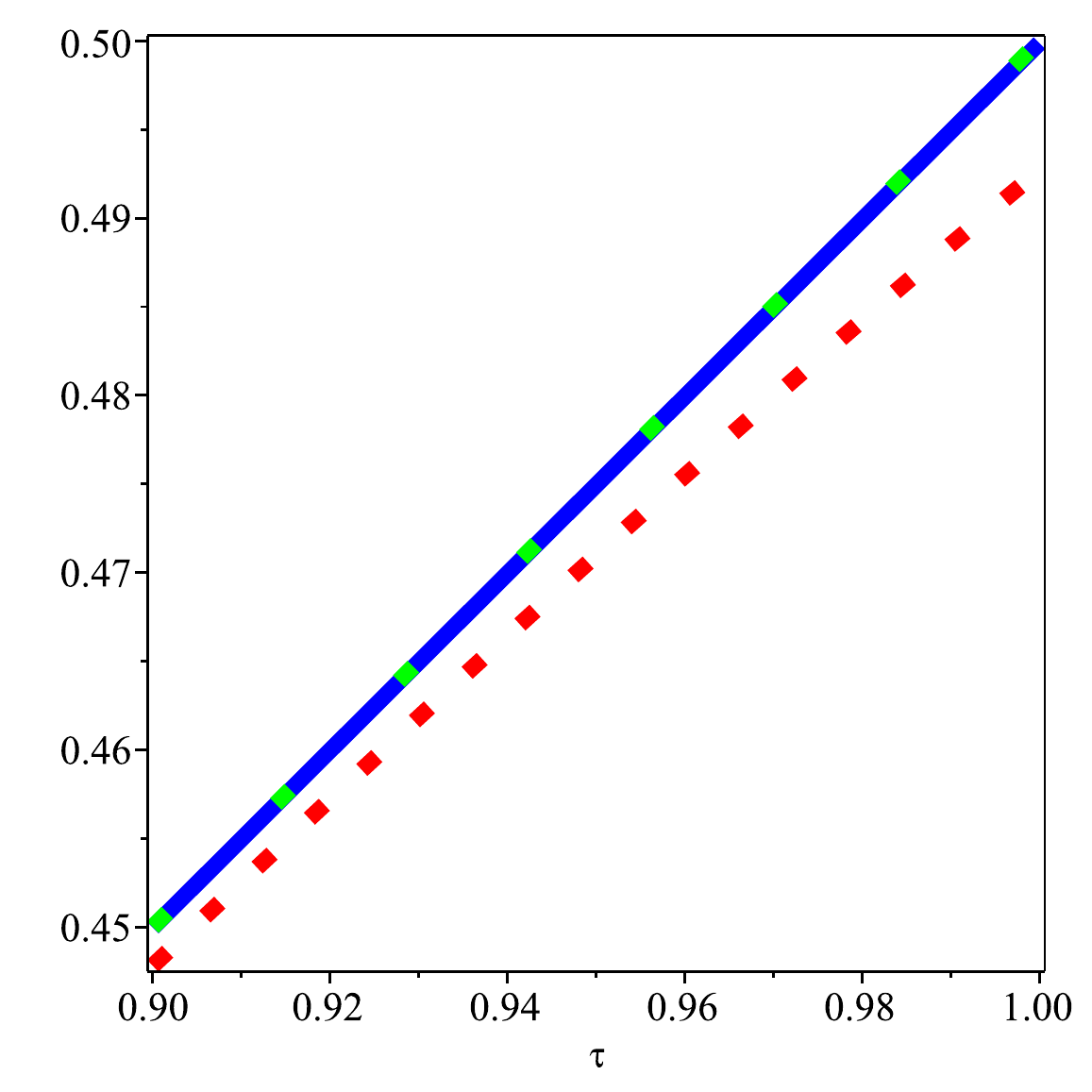}\hspace{0.3em}%
    }}\hspace{0.3em}
     \includegraphics[width=5.3cm]{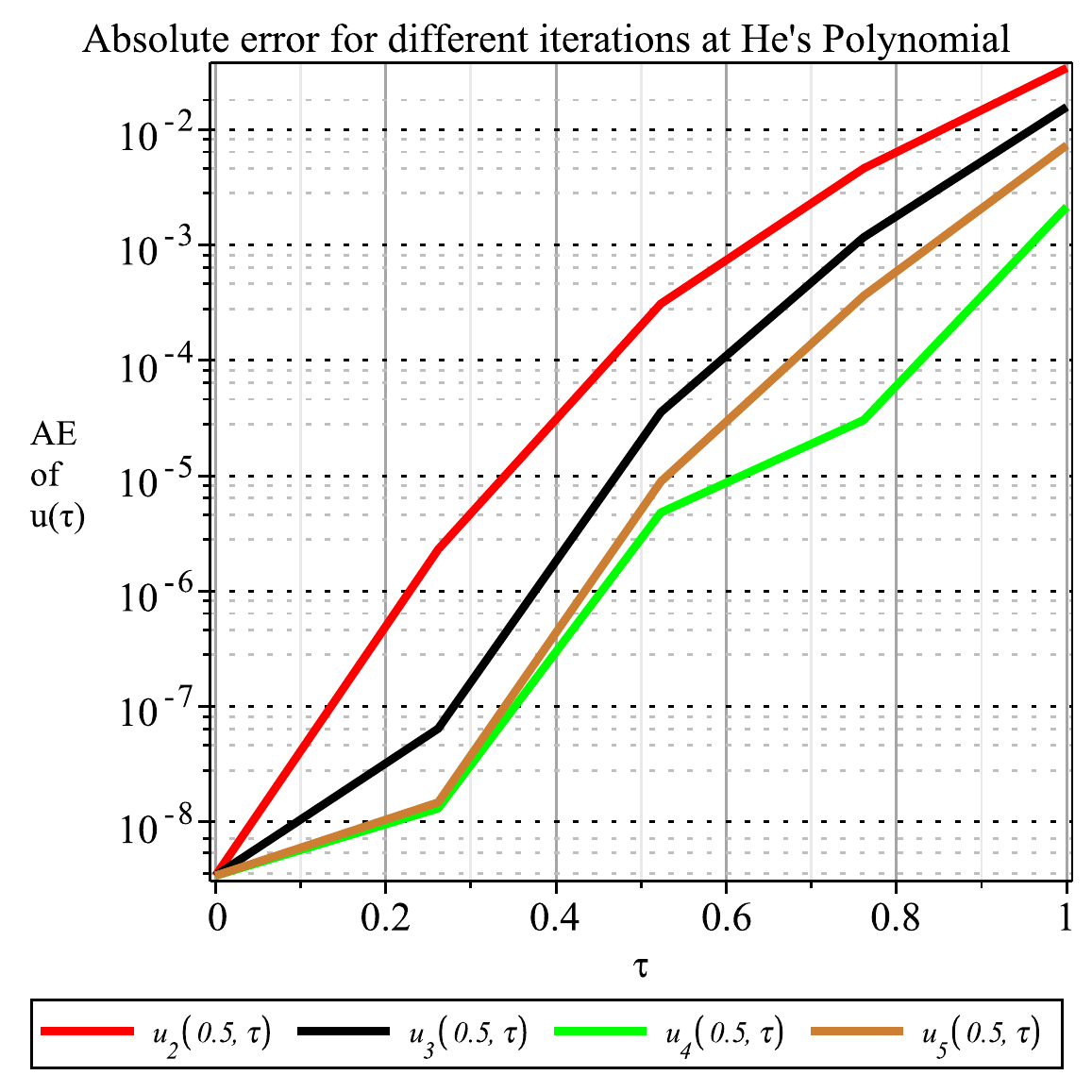}\hspace{0.3em}%
     \includegraphics[width=5.3cm]{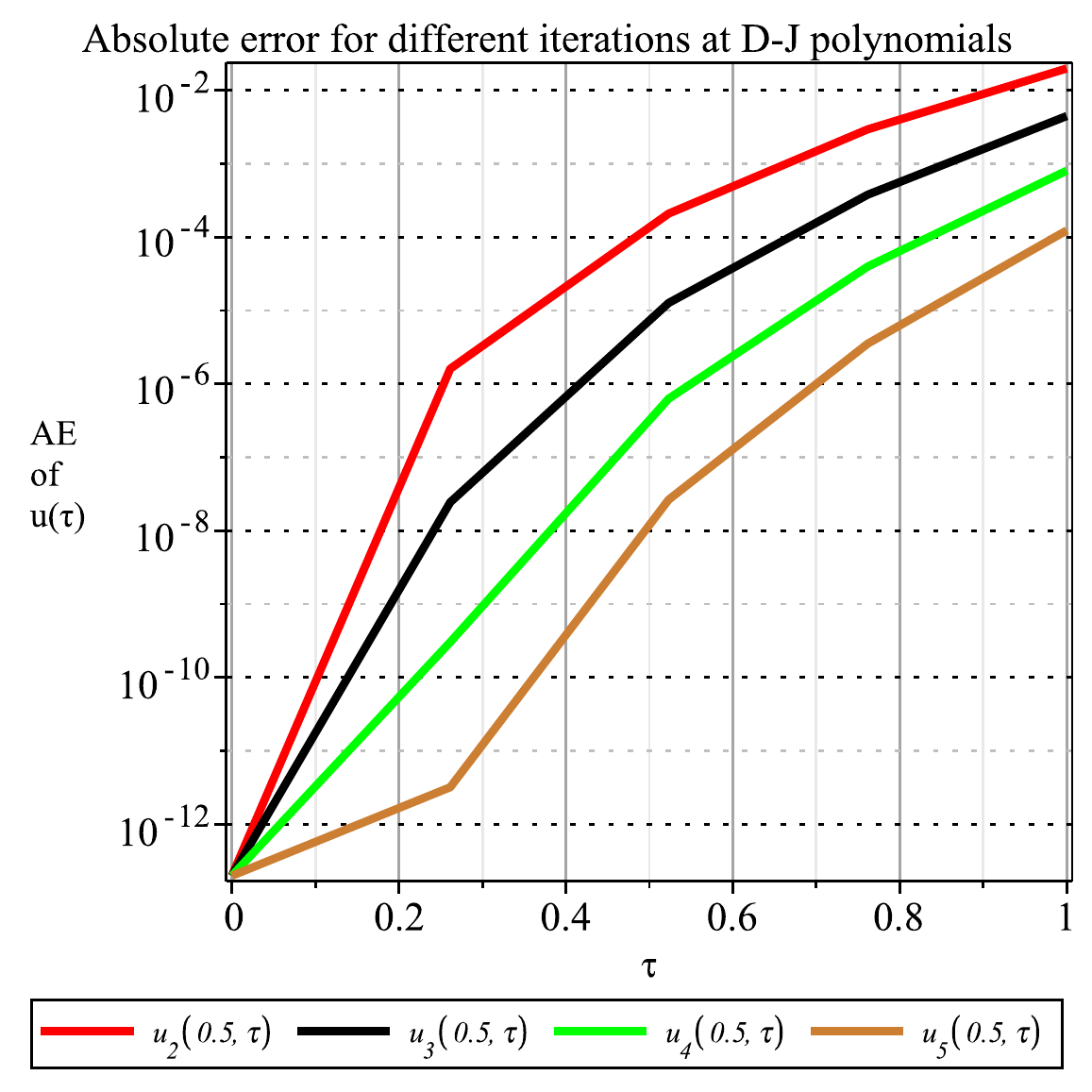}\hspace{0.3em}
    \caption{ \textbf{ Comparison of He's and D-J polynomials across different iterations, along with their associated absolute errors for problem \ref{problem1}.}}\label{fig1}
  \end{figure}

\begin{prb}\label{problem2}
Consider the nonlinear advection equation also known as nonlinear homogenous gas dynamics equation, of the form
\begin{equation}\label{fd}
{D^{\alpha}_{{{{t}}}}{{{{u}}}}}={{{u}}}\frac{\partial {{{u}}}}{\partial {{{{x}}}}}+{{{u}}}-{{{u}}}^2, \;\;\;\;\;\;\;\;\;\;\;\; 0<\alpha\leq1,
\end{equation}
With initial condition,
\begin{equation*}
{{{u}}}({{{{x}}}},0)=e^{-{{{{x}}}}}.
\end{equation*}
The exact solution at $\alpha=1$ is
$$
{{{{u}}}({{{{x}}}},{{{t}}})}=e^{{{{t}}}-{{{{x}}}}}.
$$
\end{prb}

  \begin{table}[H]
				\centering
				\caption{\textbf{Absolute error comparison of He's and D-J polynomials across different iterations ${{{k}}}$,  spaces ${0<{x}<1}$ and time level $0<{{{t}}}<1$  of problem \ref{problem2}.} In this case, we observed no difference in accuracy between the D-J polynomials and He's polynomials. Of course, one can improve the accuracy by including additional terms in the series solution. }\label{table2}
				\resizebox{\textwidth}{!}{
					\begin{tabular}{cccccccccccc}
						\toprule
						%\multicolumn{3}{c}{\textbf{ {  AE  \underline{at ${{{k}}}=0$}} }}& \multicolumn{2}{c}{\textbf{ \underline{Absolute error at ${{{k}}}=1$} }}& \multicolumn{2}{c}{\textbf{ \underline{Absolute error at ${{{k}}}=2$} }}& \multicolumn{2}{c}{\textbf{ \underline{Absolute error at ${{{k}}}=3$} }}& \multicolumn{2}{c}{\textbf{\underline{Absolute error at ${{{k}}}=4$}}} \\
						%\multicolumn{2}{c}{\textbf{ {Polynomials} }}& $\left| u_2-u_{exact} \right|$ &  $\left| u_2-u_{exact} \right|$&  $\left| u_3-u_{exact} \right|$&   $\left| u_3-u_{exact} \right|$&  $\left| u_4-u_{exact} \right|$&  $\left| u_4-u_{exact} \right|$&  He's Polynomials&   D-J Polynomials\\
						{${{{t}}}$}& ${{x}}$& {\thead{  Absolute error\\ at ${{{k}}}=1$ of\\ He's and D-J\\ Polynomials}}& {\thead{  Absolute error\\ at ${{{k}}}=2$ of\\ He's and D-J\\ Polynomials}} &{\thead{  Absolute error\\ at ${{{k}}}=3$ of\\ He's and D-J\\ Polynomials}}  & {\thead{  Absolute error\\ at ${{{k}}}=4$ of\\ He's and D-J\\ Polynomials}}&  {\thead{  Absolute error\\ at ${{{k}}}=5$ of\\ He's and D-J\\ Polynomials}} & {\thead{  Absolute error\\ at ${{{k}}}=6$ of\\ He's and D-J\\ Polynomials}}\\ %[0.1ex]
						\midrule
0.1&\num{0.1}&\num{4.68E-03}&\num{1.55E-04}&\num{3.85E-06}&\num{7.67E-08}&\num{1.27E-09}&\num{1.82E-11}\\
0.1&\num{0.3}&\num{3.83E-03}&\num{1.27E-04}&\num{3.15E-06}&\num{6.28E-08}&\num{1.04E-09}&\num{1.49E-11}\\
0.1&\num{0.5}&\num{3.14E-03}&\num{1.04E-04}&\num{2.58E-06}&\num{5.14E-08}&\num{8.55E-10}&\num{1.22E-11}\\
0.1&\num{0.7}&\num{2.57E-03}&\num{8.49E-05}&\num{2.11E-06}&\num{4.21E-08}&\num{7.00E-10}&\num{9.98E-12}\\
0.1&\num{0.9}&\num{2.10E-03}&\num{6.95E-05}&\num{1.73E-06}&\num{3.45E-08}&\num{5.73E-10}&\num{8.17E-12}\\
0.3&\num{0.1}&\num{4.51E-02}&\num{4.40E-03}&\num{3.25E-04}&\num{1.93E-05}&\num{9.57E-07}&\num{4.08E-08}\\
0.3&\num{0.3}&\num{3.69E-02}&\num{3.60E-03}&\num{2.66E-04}&\num{1.58E-05}&\num{7.83E-07}&\num{3.34E-08}\\
0.3&\num{0.5}&\num{3.02E-02}&\num{2.95E-03}&\num{2.18E-04}&\num{1.29E-05}&\num{6.41E-07}&\num{2.73E-08}\\
0.3&\num{0.7}&\num{2.48E-02}&\num{2.41E-03}&\num{1.78E-04}&\num{1.06E-05}&\num{5.25E-07}&\num{2.24E-08}\\
0.3&\num{0.9}&\num{2.03E-02}&\num{1.98E-03}&\num{1.46E-04}&\num{8.66E-06}&\num{4.30E-07}&\num{1.83E-08}\\
0.5&\num{0.1}&\num{1.35E-01}&\num{2.15E-02}&\num{2.61E-03}&\num{2.57E-04}&\num{2.11E-05}&\num{1.50E-06}\\
0.5&\num{0.3}&\num{1.10E-01}&\num{1.76E-02}&\num{2.14E-03}&\num{2.10E-04}&\num{1.73E-05}&\num{1.22E-06}\\
0.5&\num{0.5}&\num{9.02E-02}&\num{1.44E-02}&\num{1.75E-03}&\num{1.72E-04}&\num{1.42E-05}&\num{1.00E-06}\\
0.5&\num{0.7}&\num{7.39E-02}&\num{1.18E-02}&\num{1.43E-03}&\num{1.41E-04}&\num{1.16E-05}&\num{8.21E-07}\\
0.5&\num{0.9}&\num{6.05E-02}&\num{9.64E-03}&\num{1.17E-03}&\num{1.15E-04}&\num{9.50E-06}&\num{6.72E-07}\\
0.7&\num{0.1}&\num{2.84E-01}&\num{6.22E-02}&\num{1.05E-02}&\num{1.43E-03}&\num{1.64E-04}&\num{1.62E-05}\\
0.7&\num{0.3}&\num{2.32E-01}&\num{5.09E-02}&\num{8.58E-03}&\num{1.17E-03}&\num{1.34E-04}&\num{1.33E-05}\\
0.7&\num{0.5}&\num{1.90E-01}&\num{4.17E-02}&\num{7.03E-03}&\num{9.59E-04}&\num{1.10E-04}&\num{1.09E-05}\\
0.7&\num{0.7}&\num{1.56E-01}&\num{3.41E-02}&\num{5.75E-03}&\num{7.86E-04}&\num{9.00E-05}&\num{8.88E-06}\\
0.7&\num{0.9}&\num{1.28E-01}&\num{2.80E-02}&\num{4.71E-03}&\num{6.43E-04}&\num{7.37E-05}&\num{7.27E-06}\\
0.9&\num{0.1}&\num{5.06E-01}&\num{1.40E-01}&\num{3.00E-02}&\num{5.22E-03}&\num{7.64E-04}&\num{9.66E-05}\\
0.9&\num{0.3}&\num{4.15E-01}&\num{1.15E-01}&\num{2.45E-02}&\num{4.27E-03}&\num{6.26E-04}&\num{7.91E-05}\\
0.9&\num{0.5}&\num{3.39E-01}&\num{9.38E-02}&\num{2.01E-02}&\num{3.50E-03}&\num{5.12E-04}&\num{6.47E-05}\\
0.9&\num{0.7}&\num{2.78E-01}&\num{7.68E-02}&\num{1.64E-02}&\num{2.86E-03}&\num{4.20E-04}&\num{5.30E-05}\\
0.9&\num{0.9}&\num{2.28E-01}&\num{6.29E-02}&\num{1.35E-02}&\num{2.34E-03}&\num{3.43E-04}&\num{4.34E-05}\\
						\bottomrule
				\end{tabular}}
			\end{table}

\begin{figure}[H]
    \centering
    \includegraphics[width=5.3cm]{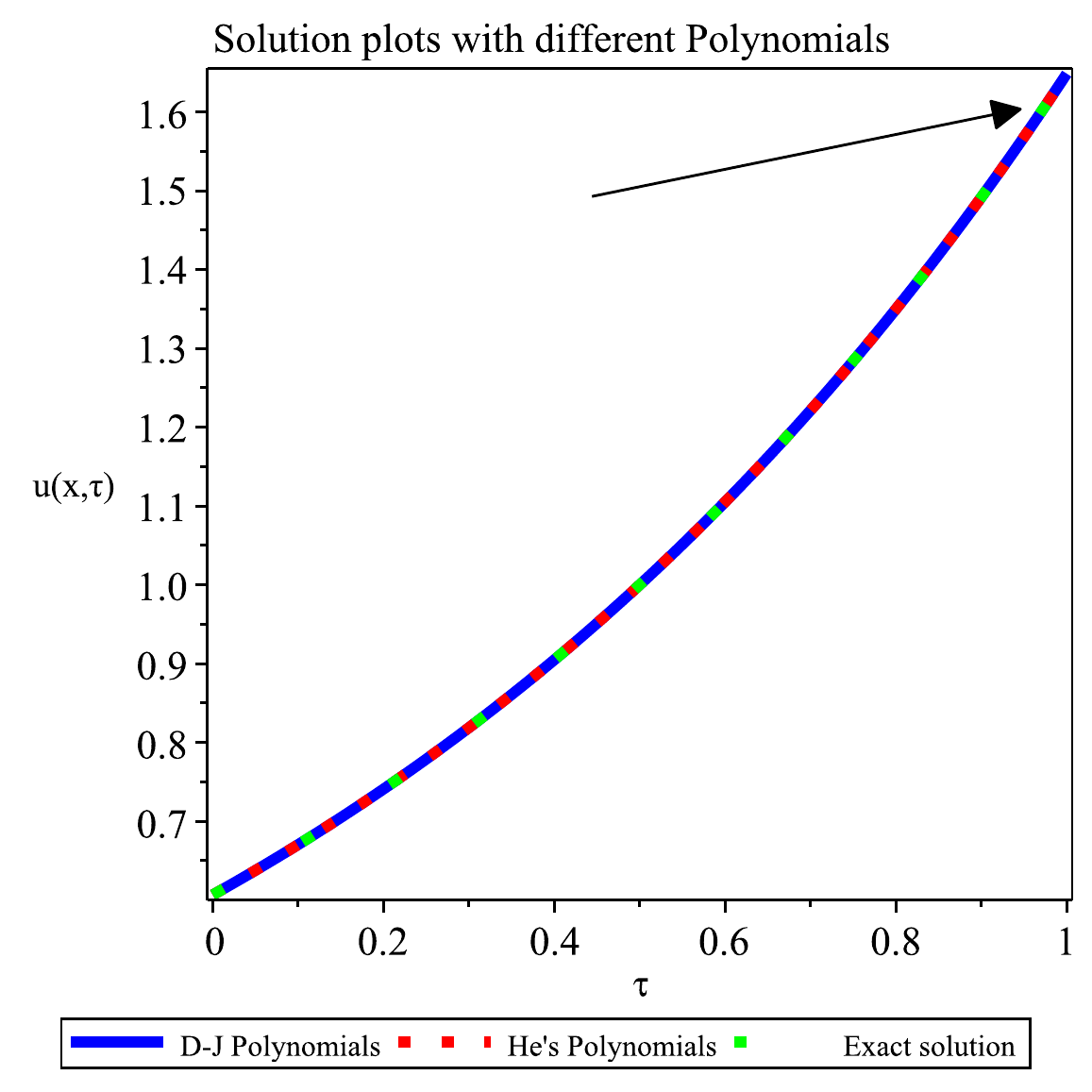}%
    \llap{\raisebox{2.64cm}{%  move next graphics to top right corner
      \includegraphics[width=2.30cm]{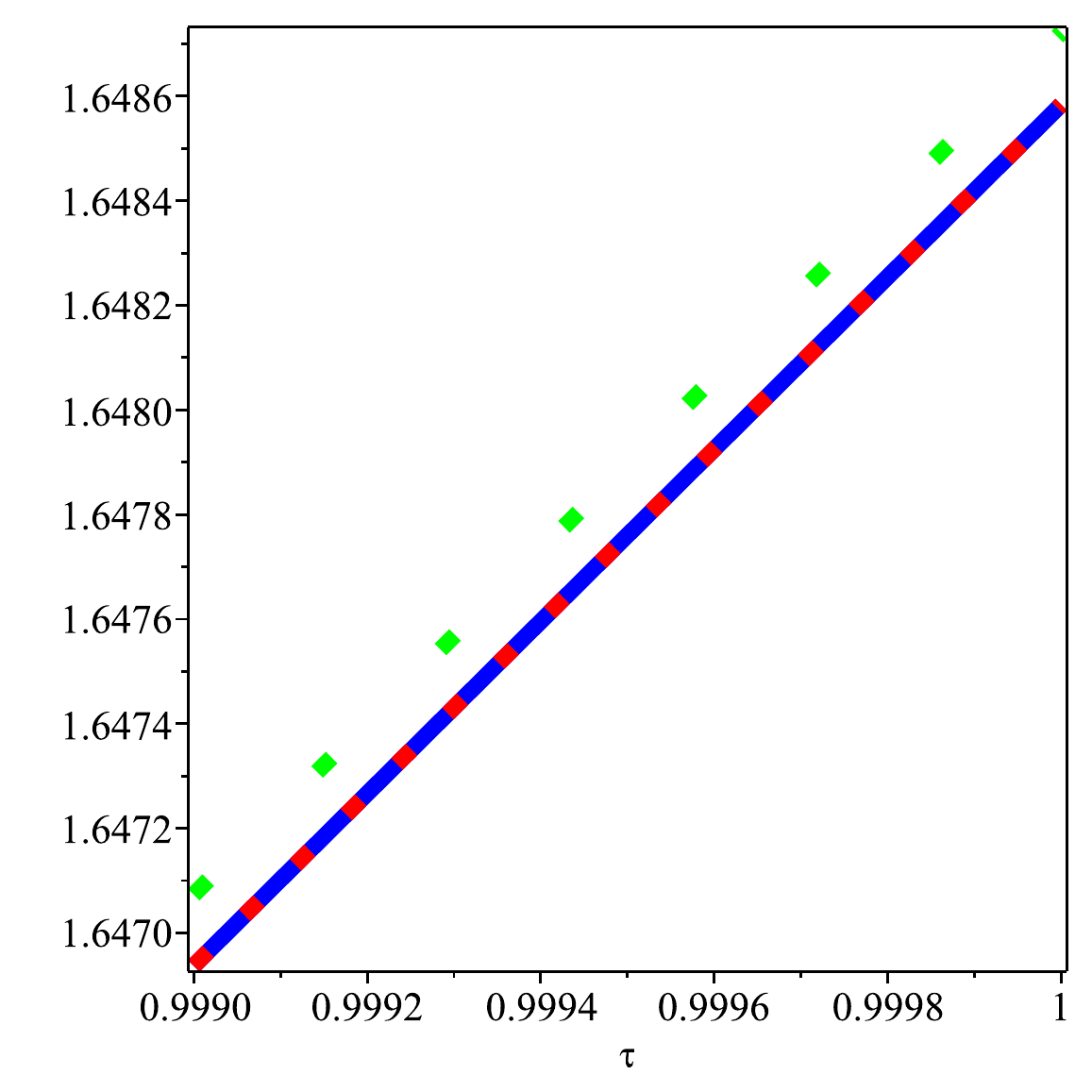}\hspace{4.80em}%
    }}\hspace{0.3em}
     \includegraphics[width=5.3cm]{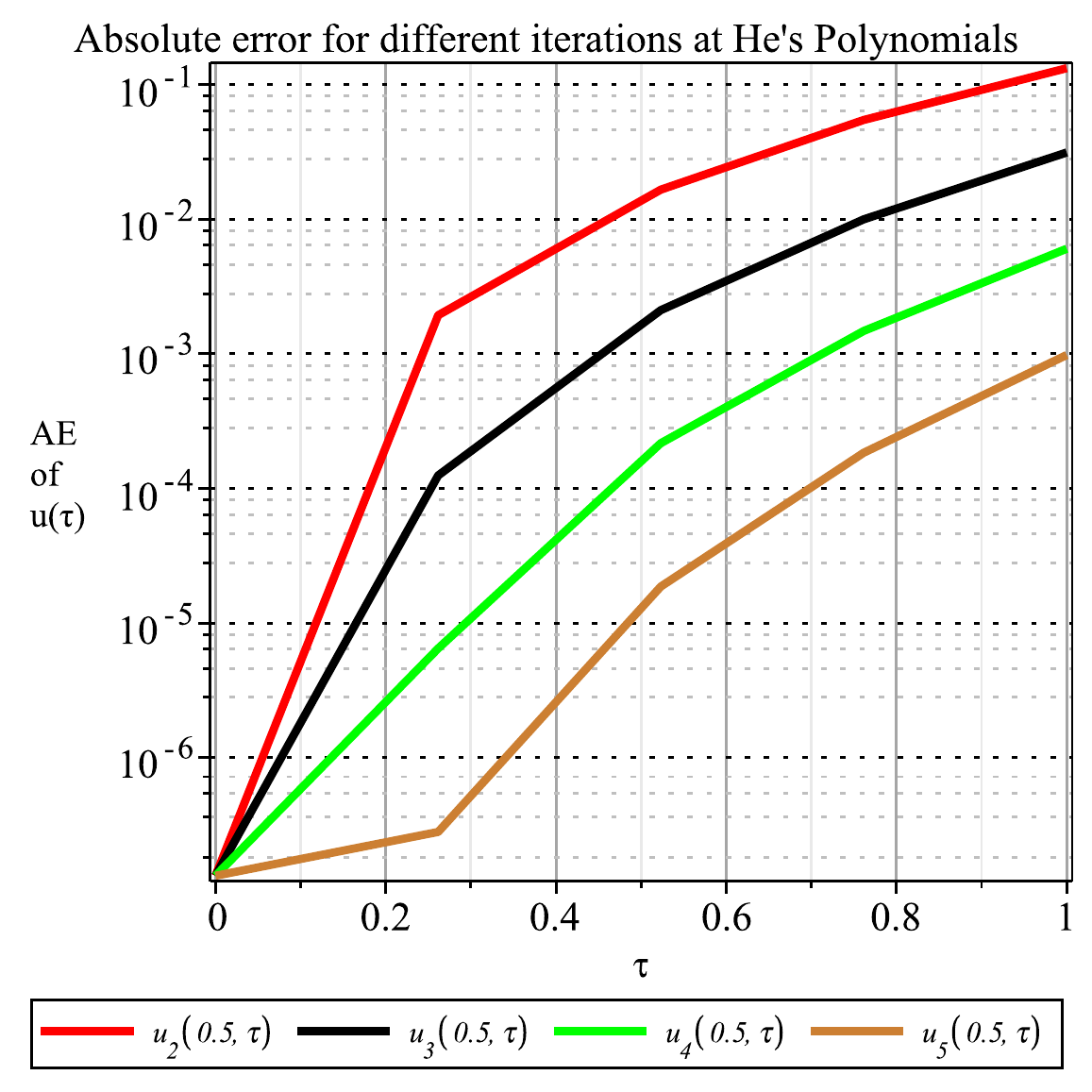}\hspace{0.3em}%
     \includegraphics[width=5.3cm]{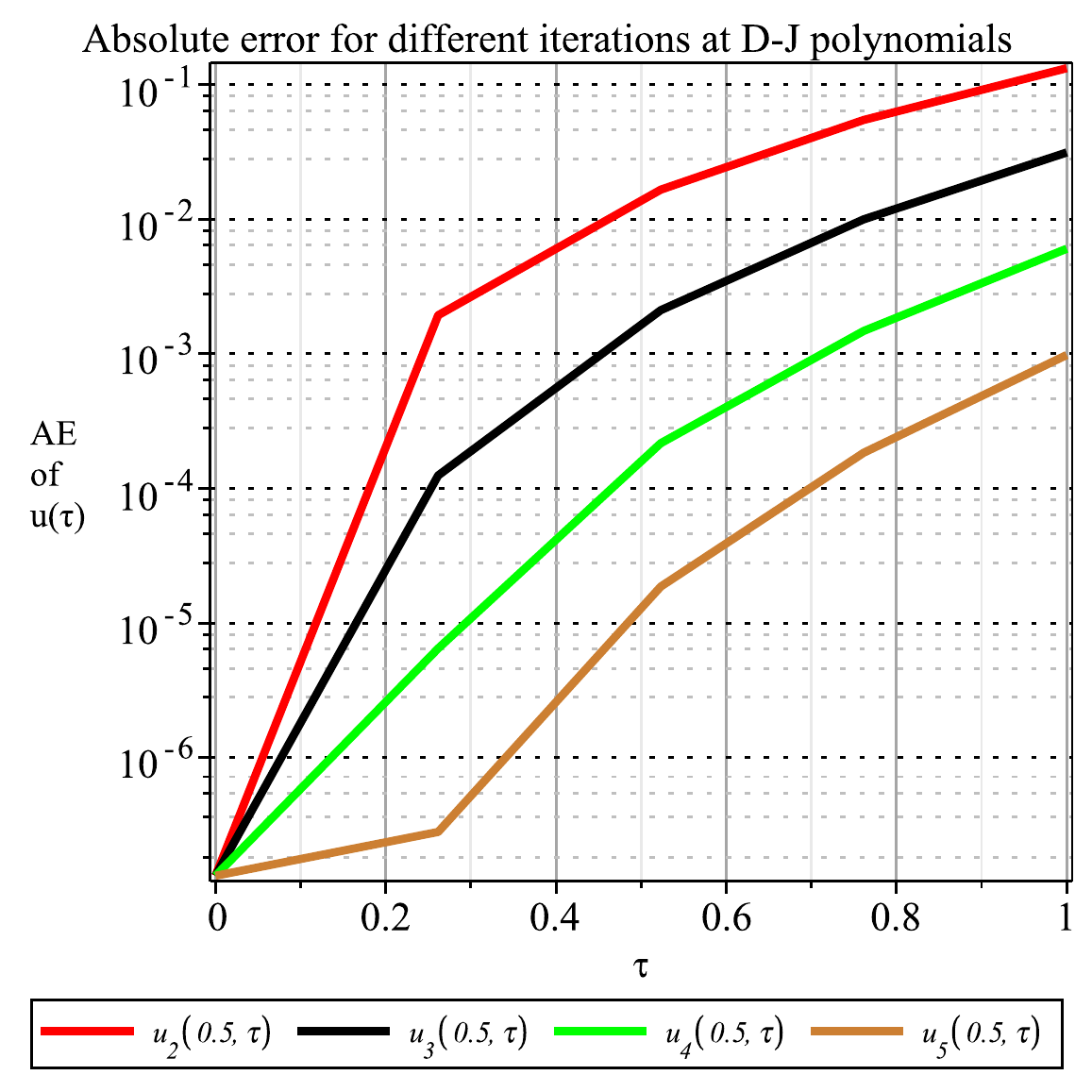}\hspace{0.3em}
    \caption{ \textbf{ Comparison of He's and D-J polynomials across different iterations, along with their associated absolute errors for problem \ref{problem2}.}}\label{fig2}
  \end{figure}

\begin{prb}\label{problem3}
Consider the nonlinear advection equation also known as nonlinear homogenous gas dynamics equation, of the form
\begin{equation}\label{fd}
{D^{\alpha}_{{{{t}}}}{{{{u}}}}}=\frac{\partial^2 {{{u}}}}{\partial {{{{x}}}^2}}+6{{{u}}}-6{{{u}}}^2, \;\;\;\;\;\;\;\;\;\;\;\; 0<\alpha\leq1,
\end{equation}
With initial condition,
\begin{equation*}
{{{u}}}({{{{x}}}},0)=\frac{1}{(1+ e^{{{{x}}}{{{{}}}}})^2}.
\end{equation*}
The exact solution at $\alpha=1$ is
$$
{{{{u}}}({{{{x}}}},{{{t}}})}=\frac{1}{(1+ e^{{{{x}}}-5{{{{t}}}}})^2}.
$$
\end{prb}
\begin{table}[H]
				\centering
				\caption{\textbf{Absolute error comparison of He's and D-J polynomials across different iterations ${{{k}}}$, spaces ${0<{x}<1}$ and time level $0<{{{t}}}<0.1$ of problem \ref{problem3}.} We observed  that the D-J polynomials work more accurately as compared to He's polynomials. Of course, one can improve the accuracy of D-J polynomials by including additional terms in the series solution, but the statement needs to be validated for He's polynomials.}\label{table3}
				\resizebox{\textwidth}{!}{
					\begin{tabular}{cccccccccccc}
						\toprule
						\multicolumn{3}{c}{\textbf{ {  AE  \underline{at ${{{k}}}=1$}} }}& \multicolumn{2}{c}{\textbf{ \underline{Absolute error at ${{{k}}}=2$} }}& \multicolumn{2}{c}{\textbf{ \underline{Absolute error at ${{{k}}}=3$} }}& \multicolumn{2}{c}{\textbf{ \underline{Absolute error at ${{{k}}}=4$} }}& \multicolumn{2}{c}{\textbf{\underline{Absolute error at ${{{k}}}=5$}}} \\
						%\multicolumn{2}{c}{\textbf{ {Polynomials} }}& $\left| u_2-u_{exact} \right|$ &  $\left| u_2-u_{exact} \right|$&  $\left| u_3-u_{exact} \right|$&   $\left| u_3-u_{exact} \right|$&  $\left| u_4-u_{exact} \right|$&  $\left| u_4-u_{exact} \right|$&  He's Polynomials&   D-J Polynomials\\
						{${{{t}}}$}& ${{x}}$& {\thead{ D-J \& He's\\ Polynomials}} &{\thead{ He's\\ Polynomials}}  & {\thead{ D-J\\ Polynomials}}&  {\thead{ He's\\ Polynomials}} & {\thead{ D-J\\ Polynomials}}& {\thead{ He's\\ Polynomials}} &  {\thead{ D-J\\ Polynomials}}& {\thead{ He's\\ Polynomials}}   &  {\thead{ D-J\\ Polynomials}}\\ %[0.1ex]
						\midrule
 0.01&{ 0.1}&\num{1.68E-04}&\num{2.13E-06}&\num{6.76E-07}&\num{6.98E-08}&\num{1.21E-08}&\num{3.88E-10}&\num{1.05E-10}&\num{ 2.45E-11}&\num{3.97E-13
}\\ 0.01&{ 0.4}&\num{1.91E-04}&\num{4.18E-07}&\num{1.44E-06}&\num{6.87E-08}&\num{1.08E-08}&\num{4.40E-10}&\num{1.85E-10}&\num{ 1.93E-11}&\num{3.29E-13
}\\ 0.01&{ 0.7}&\num{1.86E-04}&\num{1.01E-06}&\num{2.09E-06}&\num{4.74E-08}&\num{1.35E-08}&\num{8.99E-10}&\num{1.90E-10}&\num{ 5.61E-12}&\num{1.00E-12
}\\ 0.03&{ 0.1}&\num{1.47E-03}&\num{6.12E-05}&\num{1.45E-05}&\num{5.58E-06}&\num{9.78E-07}&\num{1.06E-07}&\num{2.18E-08}&\num{ 1.77E-08}&\num{3.34E-10
}\\ 0.03&{ 0.4}&\num{1.71E-03}&\num{1.51E-05}&\num{3.51E-05}&\num{5.63E-06}&\num{8.11E-07}&\num{9.72E-08}&\num{4.37E-08}&\num{ 1.44E-08}&\num{4.12E-10
}\\ 0.03&{ 0.7}&\num{1.69E-03}&\num{2.45E-05}&\num{5.37E-05}&\num{3.98E-06}&\num{9.87E-07}&\num{2.15E-07}&\num{4.53E-08}&\num{ 4.59E-09}&\num{5.96E-10
}\\ 0.05&{ 0.1}&\num{3.96E-03}&\num{3.00E-04}&\num{5.08E-05}&\num{4.24E-05}&\num{7.58E-06}&\num{1.51E-06}&\num{2.27E-07}&\num{ 3.77E-07}&\num{7.63E-09
}\\ 0.05&{ 0.4}&\num{4.71E-03}&\num{8.75E-05}&\num{1.45E-04}&\num{4.38E-05}&\num{5.83E-06}&\num{1.12E-06}&\num{5.44E-07}&\num{ 3.17E-07}&\num{1.25E-08
}\\ 0.05&{ 0.7}&\num{4.72E-03}&\num{1.00E-04}&\num{2.35E-04}&\num{3.18E-05}&\num{6.79E-06}&\num{2.72E-06}&\num{5.76E-07}&\num{ 1.09E-07}&\num{9.67E-09
}\\ 0.07&{ 0.1}&\num{7.48E-03}&\num{8.67E-04}&\num{9.61E-05}&\num{1.60E-04}&\num{2.96E-05}&\num{8.92E-06}&\num{9.16E-07}&\num{ 2.80E-06}&\num{5.69E-08
}\\ 0.07&{ 0.4}&\num{9.11E-03}&\num{2.89E-04}&\num{3.49E-04}&\num{1.69E-04}&\num{2.11E-05}&\num{5.29E-06}&\num{2.82E-06}&\num{ 2.43E-06}&\num{1.22E-07
}\\ 0.07&{ 0.7}&\num{9.29E-03}&\num{2.36E-04}&\num{6.07E-04}&\num{1.26E-04}&\num{2.31E-05}&\num{1.43E-05}&\num{3.07E-06}&\num{ 9.01E-07}&\num{4.75E-08
}\\ 0.09&{ 0.1}&\num{1.19E-02}&\num{1.93E-03}&\num{1.17E-04}&\num{4.27E-04}&\num{8.27E-05}&\num{3.40E-05}&\num{2.07E-06}&\num{ 1.25E-05}&\num{2.34E-07
}\\ 0.09&{ 0.4}&\num{1.48E-02}&\num{7.19E-04}&\num{6.37E-04}&\num{4.64E-04}&\num{5.48E-05}&\num{1.59E-05}&\num{9.50E-06}&\num{ 1.12E-05}&\num{6.73E-07
}\\ 0.09&{ 0.7}&\num{1.54E-02}&\num{4.14E-04}&\num{1.20E-03}&\num{3.54E-04}&\num{5.52E-05}&\num{4.90E-05}&\num{1.07E-05}&\num{ 4.43E-06}&\num{9.12E-08}\\
						\bottomrule
				\end{tabular}}
			\end{table}

\begin{figure}[H]
    \centering
    \includegraphics[width=5.3cm]{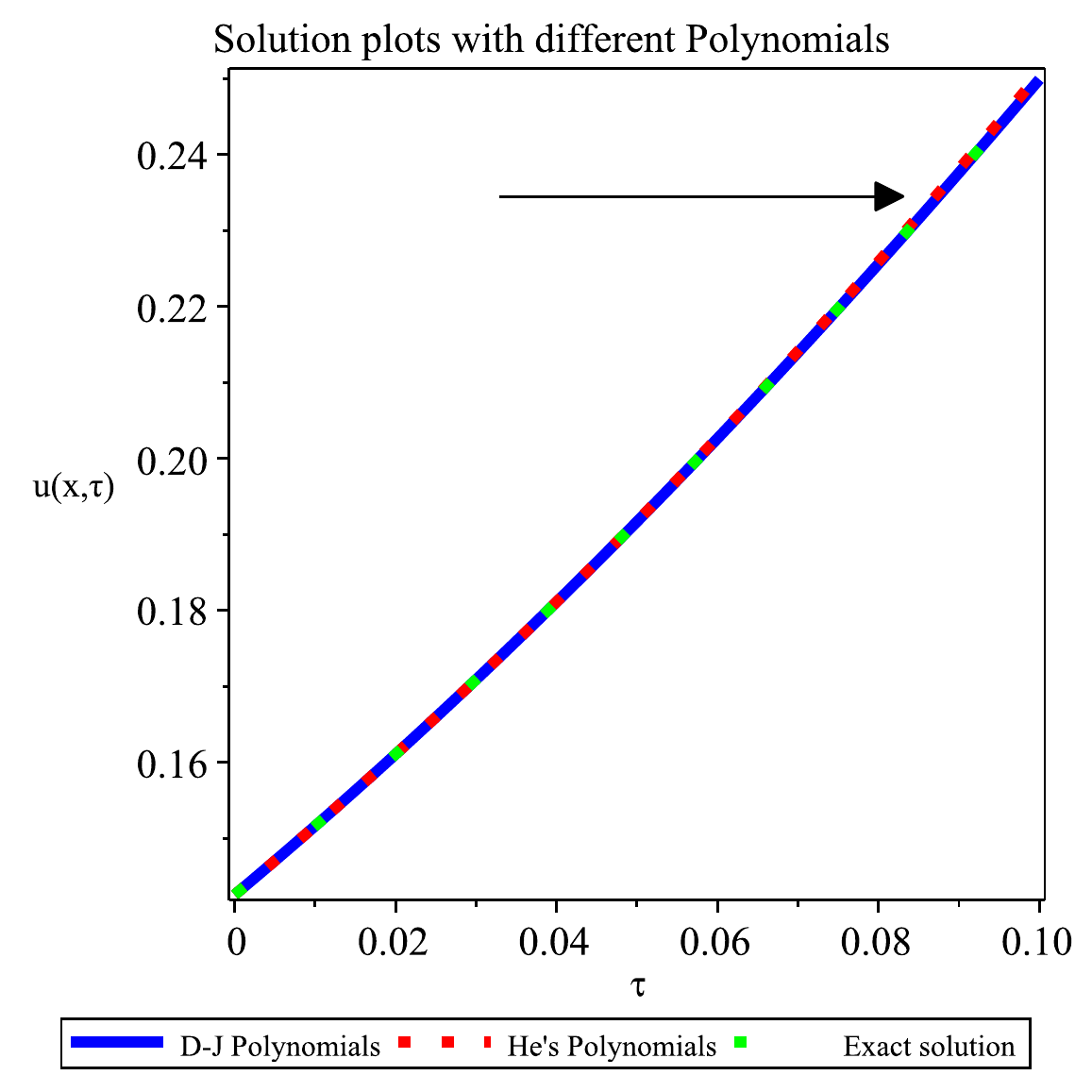}%
    \llap{\raisebox{2.85cm}{%  move next graphics to top right corner
      \includegraphics[width=2.10cm]{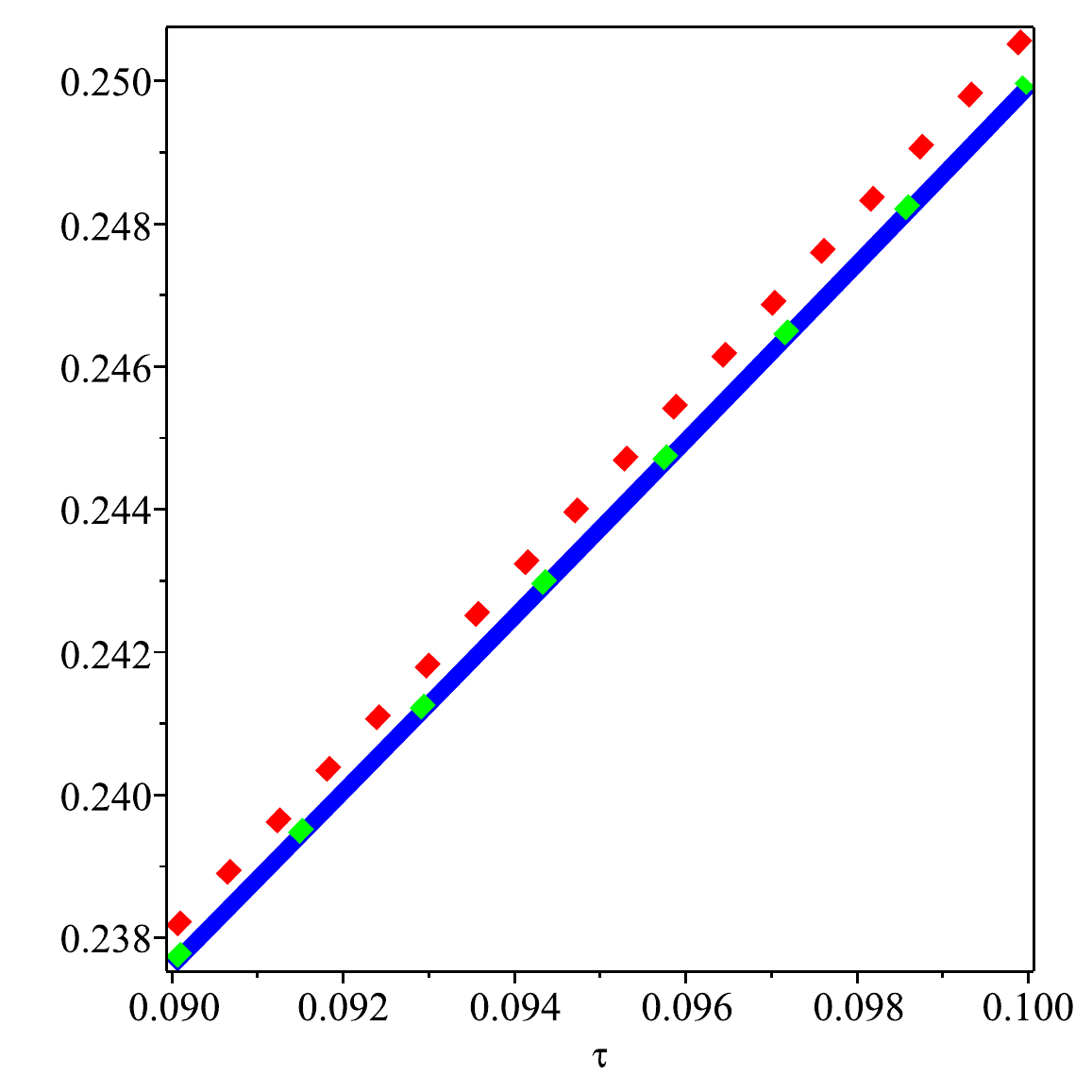}\hspace{5.05em}%
    }}\hspace{0.3em}
     \includegraphics[width=5.3cm]{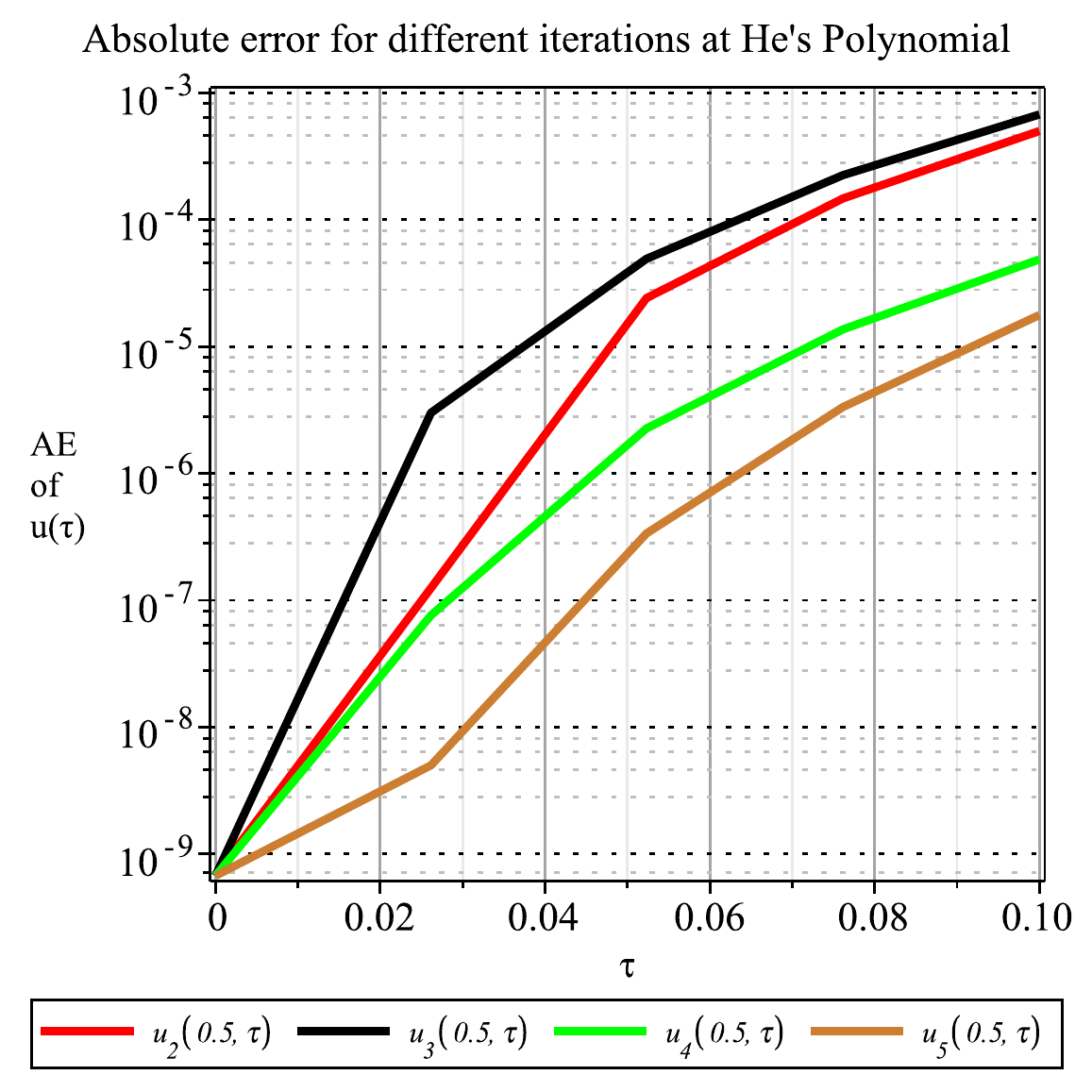}\hspace{0.3em}%
     \includegraphics[width=5.3cm]{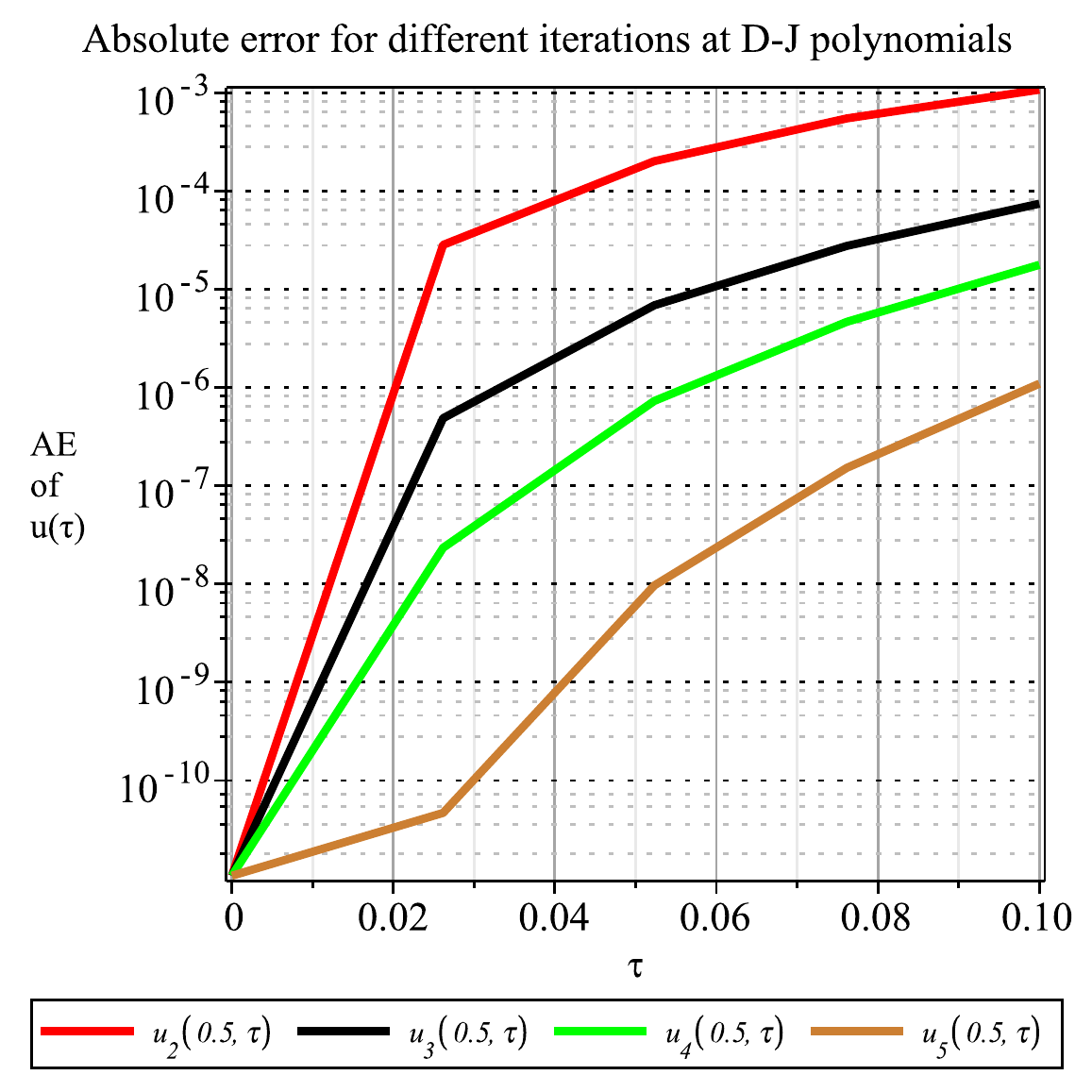}\hspace{0.3em}
    \caption{ \textbf{ Comparison of He's and D-J polynomials across different iterations, along with their associated absolute errors for problem \ref{problem3}.}}\label{fig3}
  \end{figure}
\begin{table}[H]
				\centering
				\caption{\textbf{Solution and error comparison of He's and D-J polynomials with NIM and OAFM \cite{costAhmad} at ${{{t}}}:=0.001$ of problem \ref{problem3}.} We provide more accurate solutions for the same parameters and the same number of iterations compared to the results published by Ahmad et al . \cite{costAhmad}.   We also observed that the D-J polynomials with transformation work more accurately than He's polynomials, NIM, and OAFM. }\label{table5}
				\resizebox{\textwidth}{!}{
					\begin{tabular}{cccccccccccc}
						\toprule
						\multicolumn{3}{c}{\textbf{ {    \underline{Ref to Ahmad et al .\cite{costAhmad} }} }}& \multicolumn{2}{c}{\textbf{ \underline{Proposed solutions} }}&&   \multicolumn{5}{c}{ { \underline{Error comparison of proposed methods with (ref. \cite{costAhmad})}}}  \\
						%\multicolumn{2}{c}{\textbf{ {Polynomials} }}& $\left| u_2-u_{exact} \right|$ &  $\left| u_2-u_{exact} \right|$&  $\left| u_3-u_{exact} \right|$&   $\left| u_3-u_{exact} \right|$&  $\left| u_4-u_{exact} \right|$&  $\left| u_4-u_{exact} \right|$&  He's Polynomials&   D-J Polynomials\\
					  ${{x}}$& {\thead{NIM
 \\ solution}} &{\thead{ OAFM\\ solution}}  & {\thead{He's\\ Polynomials}}&  {\thead{  D-J\\ Polynomials}} & {\thead{ Exact\\ solution}}& {\thead{ AE\\ NIM}} &  {\thead{ AE\\ OAFM}}& {\thead{ He's\\ Polynomials}}   &  {\thead{ D-J\\ Polynomials}}\\ %[0.1ex]
						\midrule
{$0.5$}&{$0.142537$}&{$0.142537$}&{$0.143426115278$}&{$0.143426115270$}&{ 0.143426115271}&\num{1.46E-03}&\num{4.80E-04}&\num{6.27E-12}&\num{1.17E-12
}\\{$0.6$}&{$0.125559$}&{$0.125559$}&{$0.126372035454$}&{$0.126372035447$}&{ 0.126372035448}&\num{1.31E-03}&\num{3.89E-04}&\num{5.53E-12}&\num{1.27E-12
}\\{$0.7$}&{$0.110099$}&{$0.110099$}&{$0.110836873586$}&{$0.110836873580$}&{ 0.110836873582}&\num{1.17E-03}&\num{7.38E-04}&\num{4.66E-12}&\num{1.42E-12
}\\{$0.8$}&{$0.096116$}&{$0.096116$}&{$0.096780772254$}&{$0.096780772249$}&{ 0.096780772250}&\num{1.04E-03}&\num{2.68E-04}&\num{3.74E-12}&\num{1.59E-12
}\\{$0.9$}&{$0.083550$}&{$0.083550$}&{$0.084145873623$}&{$0.084145873618$}&{ 0.084145873620}&\num{9.20E-04}&\num{2.32E-04}&\num{2.80E-12}&\num{1.77E-12
}\\{$1.0$}&{$0.072329$}&{$0.072329$}&{$0.072859838185$}&{$0.072859838181$}&{ 0.072859838183}&\num{8.06E-04}&\num{2.08E-04}&\num{1.90E-12}&\num{1.96E-12}\\
						\bottomrule
				\end{tabular}}
			\end{table}

\begin{table}[H]
				\centering
				\caption{\textbf{Solution comparison of He's and D-J polynomials with other existing techniques at  spaces ${0<{x}\leq1}$ and time level $0<{{{t}}}\leq 0.3$ of problem \ref{problem3}. } We provide more accurate solutions for the same parameters and the same number of iterations compared to the results published in \cite{costADM}. }\label{table4}
				\resizebox{\textwidth}{!}{
					\begin{tabular}{cccccccccccc}
						\toprule
&\multicolumn{3}{c}{\textbf{ \underline{Proposed solutions} }}&    \multicolumn{2}{c}{ \textbf{ \underline{Comparison with \cite{costADM}}}}  \\
						%\multicolumn{3}{c}{\textbf{ {  AE  \underline{at ${{{k}}}=0$}} }}& \multicolumn{2}{c}{\textbf{ \underline{Absolute error at ${{{k}}}=1$} }}& \multicolumn{2}{c}{\textbf{ \underline{Absolute error at ${{{k}}}=2$} }}& \multicolumn{2}{c}{\textbf{ \underline{Absolute error at ${{{k}}}=3$} }}& \multicolumn{2}{c}{\textbf{\underline{Absolute error at ${{{k}}}=4$}}} \\
						%\multicolumn{2}{c}{\textbf{ {Polynomials} }}& $\left| u_2-u_{exact} \right|$ &  $\left| u_2-u_{exact} \right|$&  $\left| u_3-u_{exact} \right|$&   $\left| u_3-u_{exact} \right|$&  $\left| u_4-u_{exact} \right|$&  $\left| u_4-u_{exact} \right|$&  He's Polynomials&   D-J Polynomials\\
						{${{{t}}}$}& ${{x}}$& \textbf{\thead{ D-J \\ Polynomials}} &{\thead{\textbf{ He's}  \\ \textbf{ Polynomials}}}  & \textbf{\thead{ ADM\\ Solution}}& \textbf{\thead{ VIM\\ Solution}} &  \textbf{\thead{ Exact\\ Solution}} \\ %[0.1ex]
						\midrule
0.1&$0.25$&$0.31603224$&$0.31602565$&$0.317948$&$0.315940$&$0.31604242$\\
0.1&$0.50$&$0.24998226$&$0.25004861$&$0.250500$&$0.249926$&$0.25000000$\\
0.1&$0.75$&$0.19167177$&$0.19177604$&$0.190979$&$0.191606$&$0.19168942$\\
0.1&$1.00$&$0.14252275$&$0.14262790$&$0.140979$&$0.142411$&$0.14253696$\\
0.2&$0.25$&$0.46124955$&$0.46002885$&$0.481199$&$0.459320$&$0.46128371$\\
0.2&$0.50$&$0.38698507$&$0.38832484$&$0.396941$&$0.386450$&$0.38745562$\\
0.2&$0.75$&$0.31546712$&$0.31840716$&$0.315266$&$0.315478$&$0.31604242$\\
0.2&$1.00$&$0.24956305$&$0.25285285$&$0.241175$&$0.249092$&$0.25000000$\\
0.3&$0.25$&$0.60693584$&$0.59045763$&$0.681440$&$0.591179$&$0.60419507$\\
0.3&$0.50$&$0.53236177$&$0.53588028$&$0.527635$&$0.527635$&$0.53444665$\\
0.3&$0.75$&$0.45688928$&$0.47510866$&$0.475833$&$0.459719$&$0.46128371$\\
0.3&$1.00$&$0.38369103$&$0.40723342$&$0.372917$&$0.387025$&$0.38745562$\\

						\bottomrule
				\end{tabular}}
			\end{table}

\begin{prb}\label{problem4}
Consider the following nonlinear space and time-fractional hyperbolic equation of the form,
\begin{equation}\label{11a}
{D^{\alpha}_{{{{t}}}}{{{{u}}}}}=\frac{\partial }{\partial {{{{x}}}}}\left({{{u}}}\frac{\partial^\beta {{{u}}}}{\partial {{{{x}}}^\beta}} \right), \;\;\;\;\;\;  1<\alpha\leq2,\;\;\;\;\;  0<\beta\leq1,
\end{equation}
With initial condition
$$
{{{{u}}}({{{{x}}}},0)={{{x}}}^{2\beta}},  \;\;\;\;\;\;  {{{{u}}}_{{{t}}}({{{{x}}}},0)=-{{{x}}}^{2\beta}}.
$$
The exact solution for $\alpha=2$ and $\beta=1$ is ${{{{{u}}}}}({{{x}}},{{{t}}})=\left( \frac{{{{x}}}}{{{{x}}}+1} \right)^2  $.

\end{prb}
\begin{figure}[H]
    \centering
    \includegraphics[width=5.3cm]{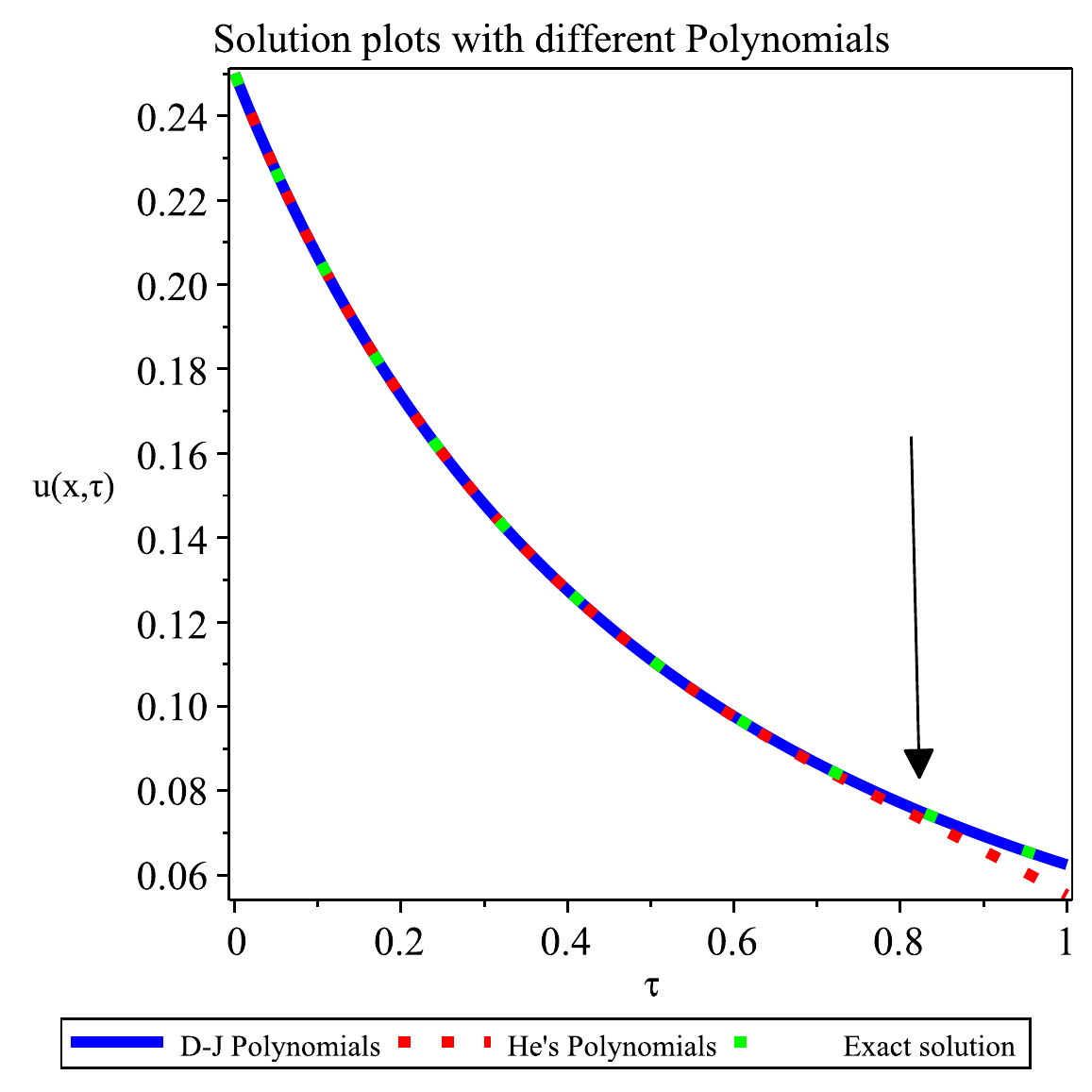}%
    \llap{\raisebox{2.49cm}{%  move next graphics to top right corner
      \includegraphics[width=2.45cm]{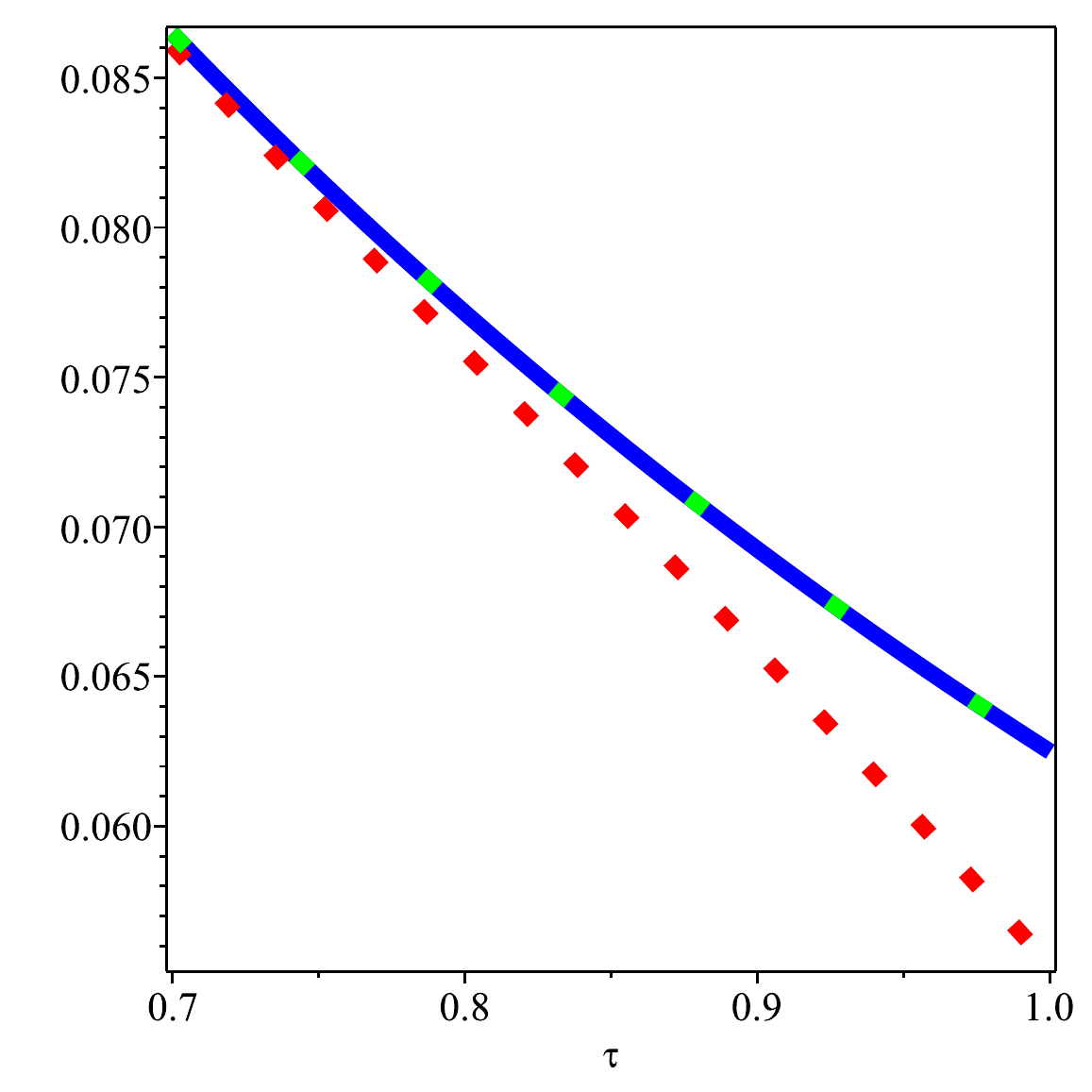}\hspace{0.3em}%
    }}\hspace{0.3em}
     \includegraphics[width=5.3cm]{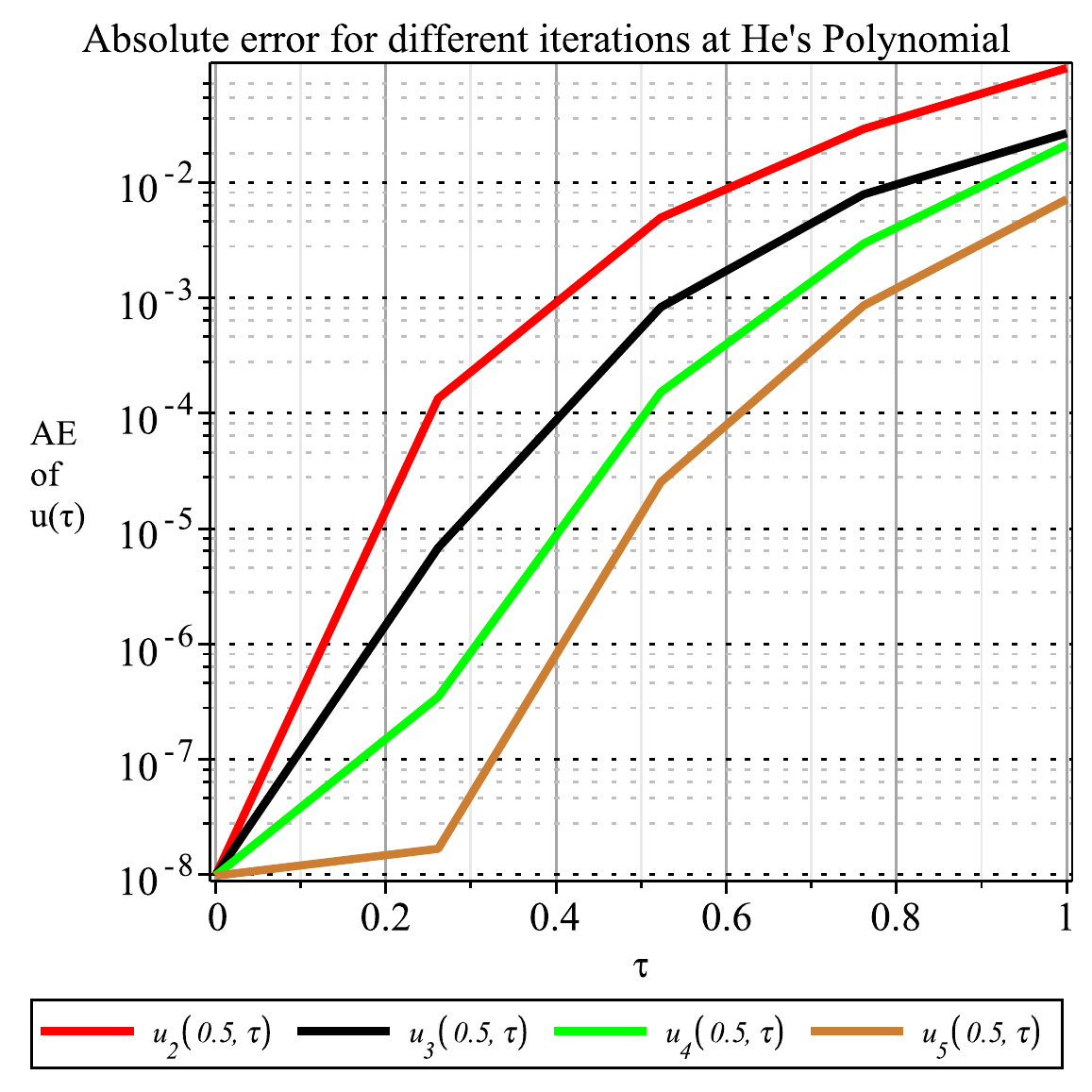}\hspace{0.3em}%
     \includegraphics[width=5.3cm]{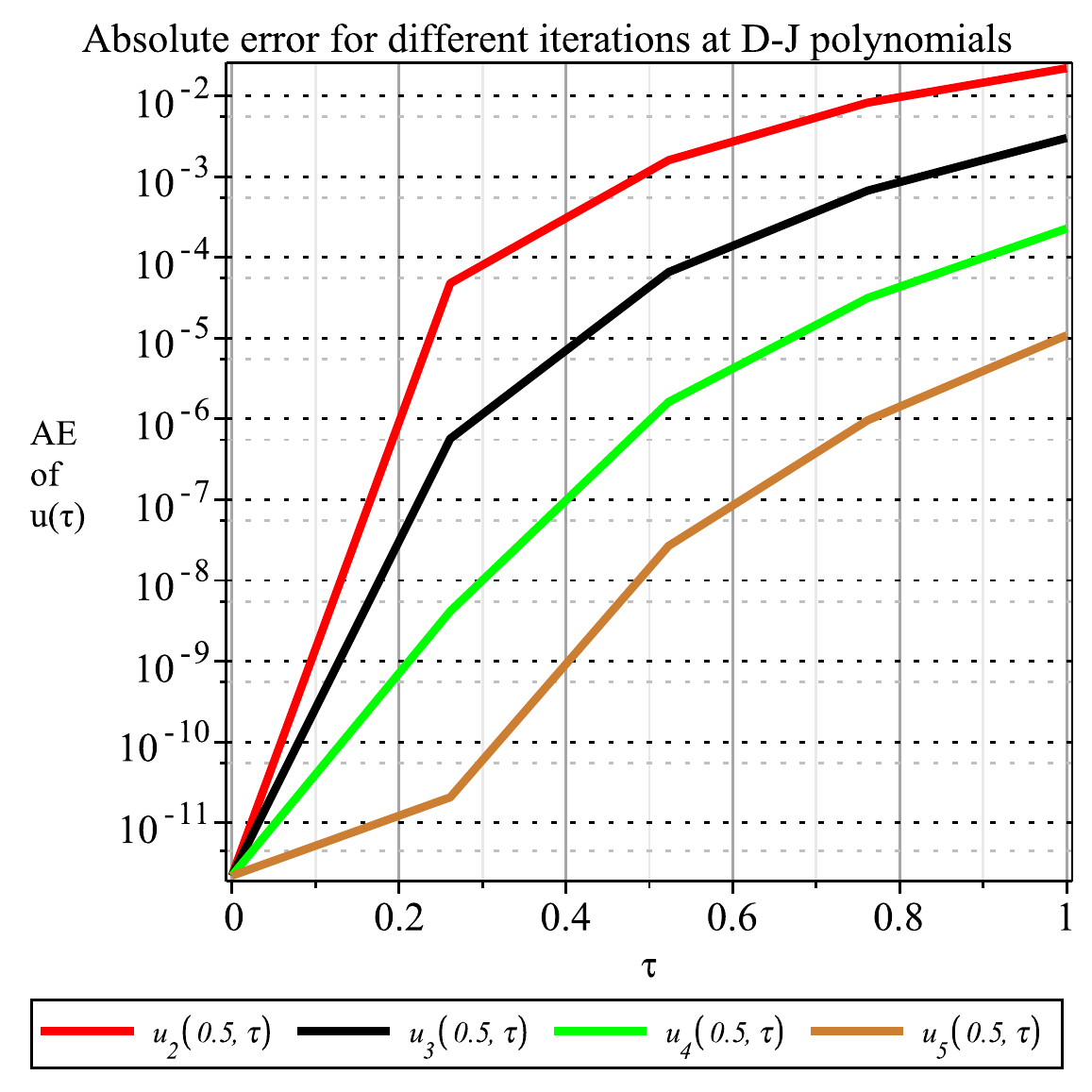}\hspace{0.3em}
    \caption{ \textbf{ Comparison of He's and D-J polynomials across different iterations, along with their associated absolute errors for problem \ref{problem4}.}}\label{fig4}
  \end{figure}

\begin{table}[H]
				\centering
				\caption{\textbf{Absolute error comparison of He's and D-J polynomials across different iterations ${{{k}}}$,  spaces ${0<{x}<1}$ and time level $0<{{{t}}}<1$  of problem \ref{problem4}.} We observed  that the D-J polynomials work more accurately as compared to He's polynomials. Of course, one can improve the accuracy of D-J polynomials by including additional terms in the series solution, but the statement needs to be validated for He's polynomials.}\label{table6}
				\resizebox{\textwidth}{!}{
					\begin{tabular}{cccccccccccc}
						\toprule
					\multicolumn{3}{c}{\textbf{ {  AE  \underline{at ${{{k}}}=1$}} }}& \multicolumn{2}{c}{\textbf{ \underline{Absolute error at ${{{k}}}=2$} }}& \multicolumn{2}{c}{\textbf{ \underline{Absolute error at ${{{k}}}=3$} }}& \multicolumn{2}{c}{\textbf{ \underline{Absolute error at ${{{k}}}=4$} }}& \multicolumn{2}{c}{\textbf{\underline{Absolute error at ${{{k}}}=5$}}} \\
						%\multicolumn{2}{c}{\textbf{ {Polynomials} }}& $\left| u_2-u_{exact} \right|$ &  $\left| u_2-u_{exact} \right|$&  $\left| u_3-u_{exact} \right|$&   $\left| u_3-u_{exact} \right|$&  $\left| u_4-u_{exact} \right|$&  $\left| u_4-u_{exact} \right|$&  He's Polynomials&   D-J Polynomials\\
						{${{{t}}}$}& ${{x}}$& {\thead{ D-J \& He's\\ Polynomials}} &{\thead{ He's\\ Polynomials}}  & {\thead{ D-J\\ Polynomials}}&  {\thead{ He's\\ Polynomials}} & {\thead{ D-J\\ Polynomials}}& {\thead{ He's\\ Polynomials}} &  {\thead{ D-J\\ Polynomials}}& {\thead{ He's\\ Polynomials}}   &  {\thead{ D-J\\ Polynomials}}\\ %[0.1ex]
						\midrule
0.1&\num{0.1}&\num{2.46E-06}&\num{2.40E-08}&\num{9.09E-09}&\num{1.95E-10}&\num{1.79E-11}&\num{1.56E-12}&\num{2.17E-14}&\num{ 1.19E-14}&\num{1.80E-17}\\
0.1&\num{0.3}&\num{2.22E-05}&\num{2.16E-07}&\num{8.18E-08}&\num{1.76E-09}&\num{1.61E-10}&\num{1.40E-11}&\num{1.96E-13}&\num{ 1.07E-13}&\num{1.62E-16}\\
0.1&\num{0.5}&\num{6.16E-05}&\num{5.99E-07}&\num{2.27E-07}&\num{4.88E-09}&\num{4.46E-10}&\num{3.89E-11}&\num{5.44E-13}&\num{ 2.97E-13}&\num{4.51E-16}\\
0.1&\num{0.7}&\num{1.21E-04}&\num{1.17E-06}&\num{4.46E-07}&\num{9.56E-09}&\num{8.75E-10}&\num{7.63E-11}&\num{1.07E-12}&\num{ 5.83E-13}&\num{8.83E-16}\\
0.1&\num{0.9}&\num{1.99E-04}&\num{1.94E-06}&\num{7.37E-07}&\num{1.58E-08}&\num{1.45E-09}&\num{1.26E-10}&\num{1.76E-12}&\num{ 9.64E-13}&\num{1.46E-15}\\
0.3&\num{0.1}&\num{1.35E-04}&\num{1.13E-05}&\num{3.95E-06}&\num{7.37E-07}&\num{6.03E-08}&\num{4.86E-08}&\num{5.64E-10}&\num{ 3.02E-09}&\num{3.57E-12}\\
0.3&\num{0.3}&\num{1.22E-03}&\num{1.02E-04}&\num{3.56E-05}&\num{6.63E-06}&\num{5.42E-07}&\num{4.37E-07}&\num{5.07E-09}&\num{ 2.72E-08}&\num{3.21E-11}\\
0.3&\num{0.5}&\num{3.38E-03}&\num{2.82E-04}&\num{9.88E-05}&\num{1.84E-05}&\num{1.51E-06}&\num{1.21E-06}&\num{1.41E-08}&\num{ 7.55E-08}&\num{8.93E-11}\\
0.3&\num{0.7}&\num{6.62E-03}&\num{5.54E-04}&\num{1.94E-04}&\num{3.61E-05}&\num{2.95E-06}&\num{2.38E-06}&\num{2.76E-08}&\num{ 1.48E-07}&\num{1.75E-10}\\
0.3&\num{0.9}&\num{1.09E-02}&\num{9.15E-04}&\num{3.20E-04}&\num{5.97E-05}&\num{4.88E-06}&\num{3.93E-06}&\num{4.57E-08}&\num{ 2.45E-07}&\num{2.89E-10}\\
0.5&\num{0.1}&\num{6.94E-04}&\num{1.59E-04}&\num{5.16E-05}&\num{2.48E-05}&\num{1.96E-06}&\num{4.20E-06}&\num{4.47E-08}&\num{ 6.46E-07}&\num{6.84E-10}\\
0.5&\num{0.3}&\num{6.25E-03}&\num{1.43E-03}&\num{4.64E-04}&\num{2.23E-04}&\num{1.76E-05}&\num{3.78E-05}&\num{4.02E-07}&\num{ 5.81E-06}&\num{6.16E-09}\\
0.5&\num{0.5}&\num{1.74E-02}&\num{3.97E-03}&\num{1.29E-03}&\num{6.20E-04}&\num{4.90E-05}&\num{1.05E-04}&\num{1.12E-06}&\num{ 1.61E-05}&\num{1.71E-08}\\
0.5&\num{0.7}&\num{3.40E-02}&\num{7.78E-03}&\num{2.53E-03}&\num{1.22E-03}&\num{9.60E-05}&\num{2.06E-04}&\num{2.19E-06}&\num{ 3.16E-05}&\num{3.35E-08}\\
0.5&\num{0.9}&\num{5.62E-02}&\num{1.29E-02}&\num{4.18E-03}&\num{2.01E-03}&\num{1.59E-04}&\num{3.40E-04}&\num{3.62E-06}&\num{ 5.23E-05}&\num{5.54E-08}\\
0.7&\num{0.1}&\num{1.68E-03}&\num{7.95E-04}&\num{2.34E-04}&\num{1.98E-04}&\num{1.63E-05}&\num{6.26E-05}&\num{6.59E-07}&\num{ 1.63E-05}&\num{1.76E-08}\\
0.7&\num{0.3}&\num{1.51E-02}&\num{7.15E-03}&\num{2.10E-03}&\num{1.78E-03}&\num{1.46E-04}&\num{5.64E-04}&\num{5.93E-06}&\num{ 1.47E-04}&\num{1.58E-07}\\
0.7&\num{0.5}&\num{4.20E-02}&\num{1.99E-02}&\num{5.85E-03}&\num{4.95E-03}&\num{4.07E-04}&\num{1.57E-03}&\num{1.65E-05}&\num{ 4.07E-04}&\num{4.40E-07}\\
0.7&\num{0.7}&\num{8.22E-02}&\num{3.89E-02}&\num{1.15E-02}&\num{9.70E-03}&\num{7.97E-04}&\num{3.07E-03}&\num{3.23E-05}&\num{ 7.98E-04}&\num{8.62E-07}\\
0.7&\num{0.9}&\num{1.36E-01}&\num{6.44E-02}&\num{1.89E-02}&\num{1.60E-02}&\num{1.32E-03}&\num{5.07E-03}&\num{5.34E-05}&\num{ 1.32E-03}&\num{1.42E-06}\\
0.9&\num{0.1}&\num{2.51E-03}&\num{2.46E-03}&\num{6.18E-04}&\num{7.21E-04}&\num{6.89E-05}&\num{4.02E-04}&\num{4.32E-06}&\num{ 1.36E-04}&\num{1.73E-07}\\
0.9&\num{0.3}&\num{2.26E-02}&\num{2.22E-02}&\num{5.56E-03}&\num{6.49E-03}&\num{6.20E-04}&\num{3.62E-03}&\num{3.88E-05}&\num{ 1.22E-03}&\num{1.56E-06}\\
0.9&\num{0.5}&\num{6.27E-02}&\num{6.16E-02}&\num{1.54E-02}&\num{1.80E-02}&\num{1.72E-03}&\num{1.00E-02}&\num{1.08E-04}&\num{ 3.40E-03}&\num{4.33E-06}\\
0.9&\num{0.7}&\num{1.23E-01}&\num{1.21E-01}&\num{3.03E-02}&\num{3.53E-02}&\num{3.38E-03}&\num{1.97E-02}&\num{2.11E-04}&\num{ 6.65E-03}&\num{8.49E-06}\\
0.9&\num{0.9}&\num{2.03E-01}&\num{2.00E-01}&\num{5.00E-02}&\num{5.84E-02}&\num{5.58E-03}&\num{3.26E-02}&\num{3.50E-04}&\num{ 1.10E-02}&\num{1.40E-05}\\
						\bottomrule
				\end{tabular}}
			\end{table}

\begin{table}[H]
				\centering
				\caption{\textbf{Solution comparison of He's and D-J polynomials with other existing techniques at  spaces ${0<{x}<1}$ and time level $0<{{{t}}}<1$  of problem \ref{problem4}.} We provide more accurate solutions for the same parameters and the same number of iterations compared to the results published in \cite{costGDTM}.   We also observed that the D-J polynomials with transformation work more accurately than He's polynomials, GDTM, ADM and VIM.}\label{table7}
				\resizebox{\textwidth}{!}{
					\begin{tabular}{cccccccccccc}
						\toprule
						%\multicolumn{3}{c}{\textbf{ {  AE  \underline{at ${{{k}}}=0$}} }}& \multicolumn{2}{c}{\textbf{ \underline{Absolute error at ${{{k}}}=1$} }}& \multicolumn{2}{c}{\textbf{ \underline{Absolute error at ${{{k}}}=2$} }}& \multicolumn{2}{c}{\textbf{ \underline{Absolute error at ${{{k}}}=3$} }}& \multicolumn{2}{c}{\textbf{\underline{Absolute error at ${{{k}}}=4$}}} \\
						%\multicolumn{2}{c}{\textbf{ {Polynomials} }}& $\left| u_2-u_{exact} \right|$ &  $\left| u_2-u_{exact} \right|$&  $\left| u_3-u_{exact} \right|$&   $\left| u_3-u_{exact} \right|$&  $\left| u_4-u_{exact} \right|$&  $\left| u_4-u_{exact} \right|$&  He's Polynomials&   D-J Polynomials\\
						{${{{t}}}$}& ${{x}}$& \textbf{\thead{ D-J \\ Polynomials}} &{\thead{\textbf{ He's}  \\ \textbf{ Polynomials}}}  & \textbf{\thead{ GDTM  \cite{costGDTM} \\ Solution}}& \textbf{\thead{ ADM \cite{costADM}\\ Solution}}& \textbf{\thead{ VIM \cite{costADM}\\ Solution}} &  \textbf{\thead{ Exact\\ Solution}} \\ %[0.1ex]
						\midrule
0.2&$0.25$&${0.043402778}$&$0.043402778$&$0.043403$&$0.043339$&$0.043403$&${0.043402778}$\\
0.2&$0.50$&$0.173611111$&$0.173611110$&$0.173611$&$0.173580$&$0.173611$&$0.173611111$\\
0.2&$0.75$&$0.390625000$&$0.390624998$&$0.390625$&$0.390556$&$0.390625$&$0.390625000$\\
0.2&$1.00$&$0.694444444$&$0.694444441$&$0.694444$&$0.694321$&$0.694444$&$0.694444444$\\
0.4&$0.25$&$0.031887755$&$0.031887347$&$0.031888$&$0.031567$&$0.031888$&$0.031887755$\\
0.4&$0.50$&$0.127551019$&$0.127549388$&$0.127551$&$0.126268$&$0.127551$&$0.127551020$\\
0.4&$0.75$&$0.286989792$&$0.286986122$&$0.286990$&$0.284103$&$0.286990$&$0.286989796$\\
0.4&$1.00$&$0.510204074$&$0.510197550$&$0.510204$&$0.505072$&$0.510204$&$0.510204082$\\
0.6&$0.25$&$0.024414037$&$0.024389883$&$0.024433$&$0.022005$&$0.024414$&$0.024414062$\\
0.6&$0.50$&$0.097656148$&$0.097559534$&$0.097730$&$0.088018$&$0.097656$&$0.097656250$\\
0.6&$0.75$&$0.219726333$&$0.219508951$&$0.219893$&$0.198040$&$0.219727$&$0.219726562$\\
0.6&$1.00$&$0.390624593$&$0.390238135$&$0.390921$&$0.352071$&$0.390625$&$0.390625000$\\
%0.8&$0.25$&$0.019067150$&$0.016726198$&$\cdots$&$\cdots$&$\cdots$&$0.019290123$\\
%0.8&$0.50$&$0.076268599$&$0.066904792$&$\cdots$&$\cdots$&$\cdots$&$0.077160494$\\
%0.8&$0.75$&$0.171604349$&$0.150535781$&$\cdots$&$\cdots$&$\cdots$&$0.173611111$\\
%0.8&$1.00$&$0.305074397$&$0.267619167$&$\cdots$&$\cdots$&$\cdots$&$0.308641975$\\
%1.0&$0.25$&$0.014872859$&$0.008928571$&$\cdots$&$\cdots$&$\cdots$&$0.015625000$\\
%1.0&$0.50$&$0.059491435$&$0.035714286$&$\cdots$&$\cdots$&$\cdots$&$0.062500000$\\
%1.0&$0.75$&$0.133855728$&$0.080357143$&$\cdots$&$\cdots$&$\cdots$&$0.140625000$\\
%1.0&$1.00$&$0.237965739$&$0.142857143$&$\cdots$&$\cdots$&$\cdots$&$0.250000000$\\
						\bottomrule
				\end{tabular}}
			\end{table}

\begin{prb}\label{problem5}
We identified a huge computational cost issue in problem three related to the D-J polynomial, due to its complex initial condition and diffusion term, which involve more iterative terms in the polynomials (see results and discussion). Now, we are examining the time-fractional nonlinear KDV-Burgers equation from a recent preprint (see problem 4.1, Khalil  et al. \cite{costrashid}), which involves even more complex initial conditions with third-order derivatives. Our focus is not on finding a solution but on evaluating the computational cost for D-J polynomials. Of course, we achieve more accurate results compared to Khalil  et al. \cite{costrashid} and the references cited therein.
\begin{equation}\label{fd}
{D^{\alpha}_{{{{t}}}}{{{{u}}}}}=6{{{u}}}\frac{\partial {{{u}}}}{\partial {{{{x}}}}}-\frac{\partial^3 {{{u}}}}{\partial {{{{x}}}^3}}, \;\;\;\;\;\;\;\;\;\;\;\; 0<\alpha\leq1,
\end{equation}
we have the initial condition,
\begin{equation*}
{{{u}}}({{{{x}}}},0)=-\frac{2 e^{c{{{x}}}}c^2}{(1+ e^{{{{x}}}{{{{c}}}}})^2}.
\end{equation*}
The exact solution at $\alpha=1$ is
$$
{{{{u}}}({{{{x}}}},{{{t}}})}=-\frac{2 e^{c{{{x}}}-c^3 {{{t}}}}c^2}{(1+ e^{c{{{x}}}-c^3 {{{t}}}})^2}.
$$
\end{prb}

\begin{figure}[H]
    \centering
    \includegraphics[width=5.3cm]{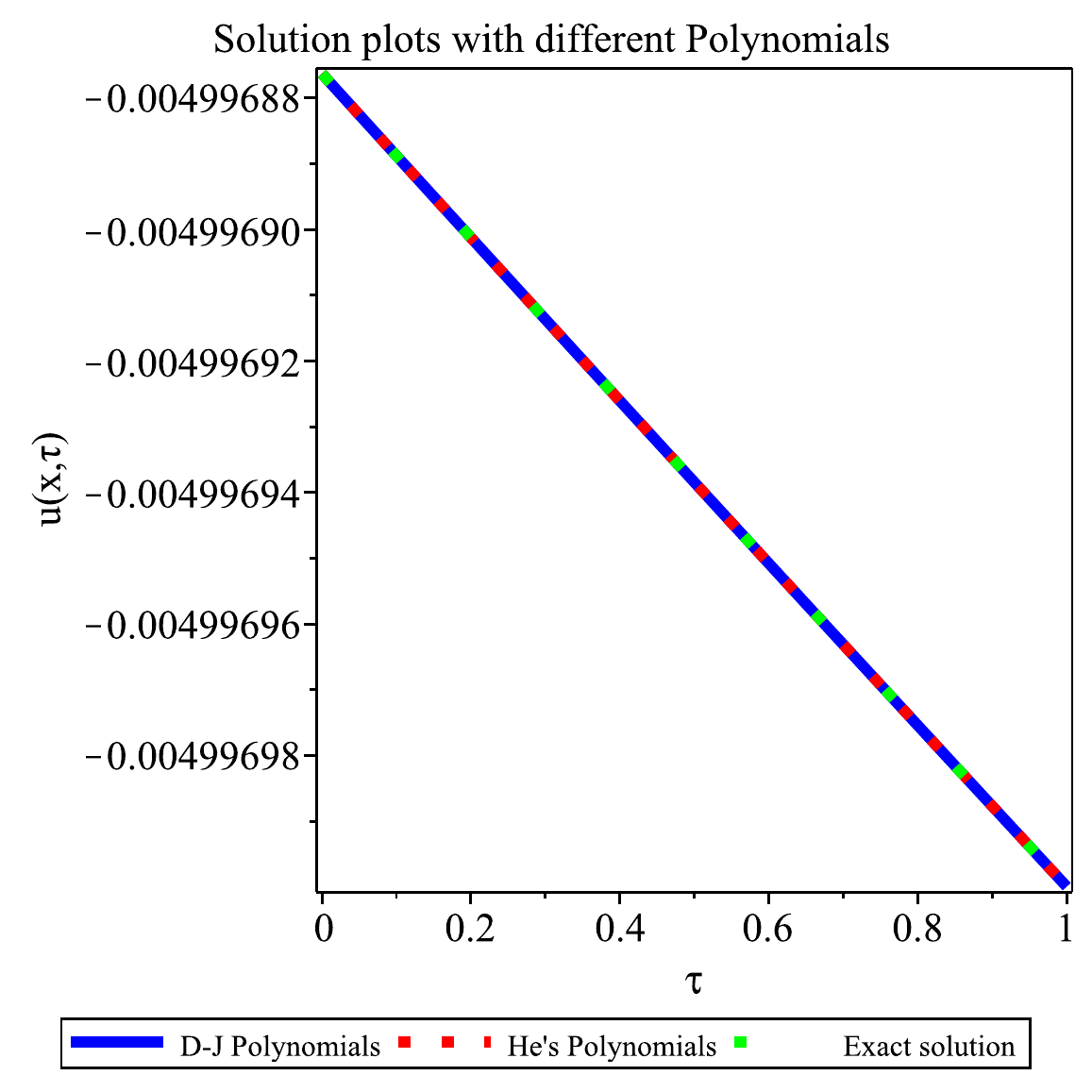}\hspace{0.3em}%
    %\llap{\raisebox{2.49cm}{%  move next graphics to top right corner
%      \includegraphics[width=2.45cm]{na3-eps-converted-to}\hspace{0.3em}%
%    }}\hspace{0.3em}
     \includegraphics[width=5.3cm]{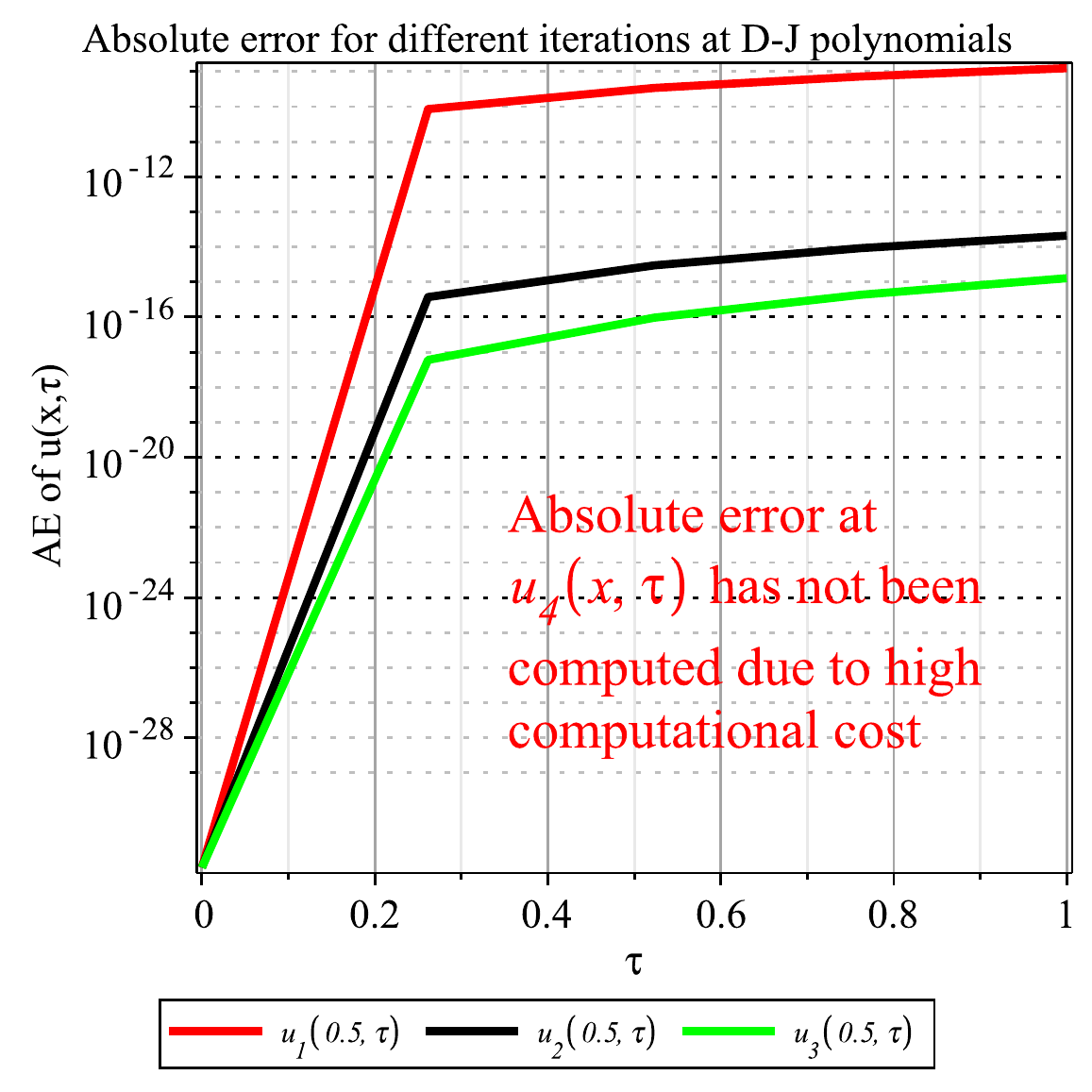}\hspace{0.3em}%
     \includegraphics[width=5.3cm]{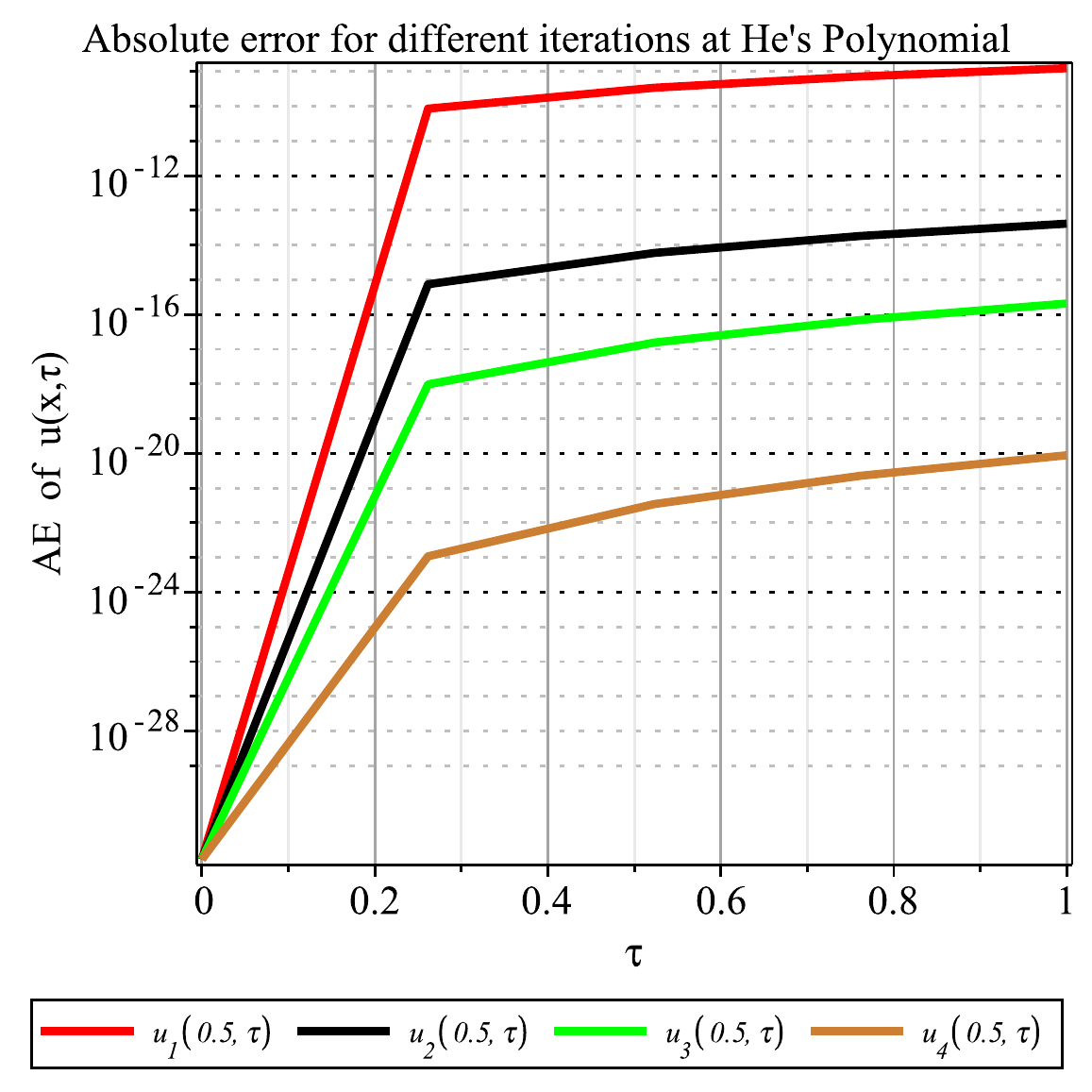}\hspace{0.3em}
    \caption{ \textbf{ Comparison of He's and D-J polynomials across different iterations, along with their associated absolute errors for problem \ref{problem5}.}}\label{fign1}
  \end{figure}

  \begin{table}[H]
				\centering
				\caption{\textbf{Absolute error comparison of He's and D-J polynomials across different iterations ${{{k}}}$,  spaces ${0<{x}<1}$ and time level $0<{{{t}}}<1$  of problem \ref{problem4}.} We observe higher  accuracy for the same parameters   compared to the results published in \cite{costrashid}.}\label{tablen6}
				\resizebox{\textwidth}{!}{
					\begin{tabular}{cccccccccccc}
						\toprule
					\multicolumn{3}{c}{\textbf{ {  AE  \underline{at ${{{k}}}=1$}} }}& \multicolumn{2}{c}{\textbf{ \underline{Absolute error at ${{{k}}}=2$} }}& \multicolumn{2}{c}{\textbf{ \underline{Absolute error at ${{{k}}}=3$} }}& \multicolumn{2}{c}{\textbf{ \underline{Absolute error at ${{{k}}}=4$} }}& \multicolumn{2}{c}{\textbf{\underline{Absolute error from \cite{costrashid}}}} \\
						%\multicolumn{2}{c}{\textbf{ {Polynomials} }}& $\left| u_2-u_{exact} \right|$ &  $\left| u_2-u_{exact} \right|$&  $\left| u_3-u_{exact} \right|$&   $\left| u_3-u_{exact} \right|$&  $\left| u_4-u_{exact} \right|$&  $\left| u_4-u_{exact} \right|$&  He's Polynomials&   D-J Polynomials\\
						{${{{t}}}$}& ${{x}}$& {\thead{ D-J \& He's\\ Polynomials}} &{\thead{ He's\\ Polynomials}}  & {\thead{ D-J\\ Polynomials}}&  {\thead{ He's\\ Polynomials}} & {\thead{ D-J\\ Polynomials}}& {\thead{ He's\\ Polynomials}} &  {\thead{ D-J\\ Polynomials}}& {\thead{ NTIM\\ (see \cite{costrashid})}}   &  {\thead{ q-HAM\\ (see \cite{costrashid}) }}\\ %[0.1ex]
			\midrule			0.1&{0.1}&\num{1.25E-11}&\num{8.31E-18}&\num{4.19E-18}&\num{2.08E-20}&\num{4.27E-19}&\num{1.44E-24}& High & \num{4.99E-9}&\num{ 8.72E-18}\\
0.1&{ 0.2}&\num{5.00E-11}&\num{6.63E-17}&\num{3.37E-17}&\num{3.33E-19}&\num{1.04E-18}&\num{1.15E-23}& computational &\num{9.99E-9}&\num{ 6.65E-17}\\
0.1&{ 0.3}&\num{1.12E-10}&\num{2.23E-16}&\num{1.14E-16}&\num{1.69E-18}&\num{9.08E-18}&\num{3.79E-23}&   cost  &\num{1.12E-8}&\num{ 2.23E-16}\\
0.1&{ 0.4}&\num{2.00E-10}&\num{5.28E-16}&\num{2.72E-16}&\num{5.33E-18}&\num{3.15E-17}&\num{8.54E-23}&   are  required  &\num{1.49E-8}&\num{ 5.28E-16}\\
0.3&{ 0.1}&\num{1.25E-11}&\num{2.50E-17}&\num{1.25E-17}&\num{2.08E-20}&\num{1.61E-18}&\num{1.53E-24}&  to compute &\num{1.99E-8}&\num{ 2.65E-17}\\
0.3&{ 0.2}&\num{5.00E-11}&\num{2.00E-16}&\num{1.00E-16}&\num{3.33E-19}&\num{1.34E-18}&\num{1.14E-23}&   this error  &\num{1.49E-8}&\num{ 2.01 E-16}\\
0.3&{ 0.3}&\num{1.12E-10}&\num{6.73E-16}&\num{3.38E-16}&\num{1.68E-18}&\num{5.48E-18}&\num{3.28E-23}&   calculation&\num{2.99E-8}&\num{ 6.74E-16}\\
0.3&{ 0.4}&\num{4.99E-8}&\num{1.59E-15}&\num{8.03E-16}&\num{5.32E-18}&\num{2.66E-17}&\num{5.75E-23}& N/A&\num{5.99E-8}&\num{ 1.59E-15}\\
0.5&{ 0.1}&\num{1.25E-11}&\num{4.16E-17}&\num{2.07E-17}&\num{2.07E-20}&\num{2.78E-18}&\num{1.63E-24}& N/A&\num{2.49E-8}&\num{ 3.99E-17}\\
0.5&{ 0.2}&\num{4.99E-11}&\num{3.32E-16}&\num{1.66E-16}&\num{3.32E-19}&\num{3.70E-18}&\num{1.13E-23}& N/A&\num{4.99E-8}&\num{ 3.31E-16}\\
0.5&{ 0.3}&\num{1.12E-10}&\num{1.12E-15}&\num{5.61E-16}&\num{1.68E-18}&\num{1.86E-18}&\num{2.77E-23}& N/A&\num{7.48E-8}&\num{1.11E-15}\\
0.5&{ 0.4}&\num{2.00E-10}&\num{2.66E-15}&\num{1.33E-15}&\num{5.31E-18}&\num{2.16E-17}&\num{2.95E-23}& N/A&\num{9.98E-8}&\num{2.66E-15}\\
						\bottomrule
				\end{tabular}}
			\end{table}

\section{Results, Discussion and Computational Cost}\label{result}
In this section, we evaluate the accuracy of the proposed method and thoroughly explore the details of the results obtained. As previously mentioned, we utilized an accurate and straightforward procedure for solving non-linear fractional partial differential equations (PDEs). This method yields highly precise results, which we present with comprehensive numerical data, supported by graphs and tables. All visualizations, including graphs and tables, were generated using Maple Software version 2024. The iterative computations were conducted on a Surface Pro 6 equipped with an Intel Core™ i7-8650U processor. This powerful CPU allowed for efficient processing during the calculations.

For the non-linear part of the analysis, we used He's polynomials and included up to 5 terms to achieve the required accuracy. In parallel, we also employed the 5th term of the D-J polynomials to ensure a robust comparison. This approach allowed us to accurately assess the performance of both types of polynomials in handling the non-linear aspects of the problems.
For better understanding, we have divided this section into two parts: the first subsection covers "{Results and Discussion}," and the second subsection covers "Computational Cost."

\subsection{Results and Discussion}

Figures \ref{fig1} through \ref{fign1} present a comprehensive comparison of He's polynomials and D-J polynomials across various iterations for Problems \ref{problem1} through \ref{problem5}. Each figures contains sub-figures that display the solution plots for both polynomial types alongside their respective exact solutions. In the plots, we observe that the solutions generated by D-J polynomials consistently align more closely with the exact solutions than those produced by He's polynomials. This trend is evident across all five problems, indicating that D-J polynomials provide greater accuracy in approximating the true solutions. Additionally, the accompanying absolute error values for each polynomial type are depicted, further illustrating the superior performance of D-J polynomials in minimizing deviations from the exact solutions. This analysis underscores the effectiveness of D-J polynomials in achieving higher precision in computational tasks compared to He's polynomials.

In Table \ref{table1}, we present the absolute error comparisons of He’s and Daftardar-Jafari (D-J) polynomials across various iterations ${{{k}}}$, for the time interval $0<{{{t}}}<1$ and spatial domain ${0<{x}<1}$ in problem \ref{problem1}. The results indicate that D-J polynomials exhibit superior accuracy compared to He’s polynomials. While it is possible to enhance the accuracy of the D-J polynomials by including additional terms in the series solution, further validation is required for He’s polynomials to ascertain whether similar improvements can be achieved. Table \ref{table2} provides a similar absolute error comparison for He’s and D-J polynomials under the same conditions as in problem \ref{problem2}. Interestingly, we observed no significant difference in accuracy between the two polynomial methods. This suggests that both approaches yield comparable results under these specific conditions. Nevertheless, enhancing accuracy through the inclusion of additional terms in the series solution remains applicable to both methods. In Table \ref{table3}, we analyze the absolute error for He’s and D-J polynomials across different iterations ${{{k}}}$, focusing on the time interval $0<{{{t}}}<1$ and the spatial range ${0<{x}<1}$ in problem \ref{problem3}. Consistent with the findings in problem \ref{problem1}, D-J polynomials are again shown to be more accurate than He’s polynomials. As noted previously, the accuracy of D-J polynomials can be improved by including more terms in the series, while the potential for similar enhancements in He’s polynomials remains to be validated. Table \ref{table4} compares the solutions obtained from He’s and D-J polynomials with other existing techniques at time levels $0<{{{t}}}<0.3$ and spatial domain ${0<{x}<1}$ in problem \ref{problem3}. Our results demonstrate that both polynomial methods provide more accurate solutions than those reported in \cite{costGDTM} for the same parameters and number of iterations. This highlights the effectiveness of both approaches in delivering precise solutions within the specified constraints. Table \ref{table5} presents a comparison of solutions and errors between He’s and D-J polynomials, New Iteration Method (NIM), and Optimal Auxiliary Function Method (OAFM) at ${{{t}}}=0.001$ in problem \ref{problem3}. The findings reveal that both He’s and D-J polynomials yield more accurate solutions than those provided by Ahmad et al. \cite{costAhmad}. Notably, D-J polynomials, particularly when combined with transformation techniques, outperform He’s polynomials, NIM, and OAFM, indicating a robust performance in this context. Table \ref{table6} continues the absolute error comparison for He’s and D-J polynomials across various iterations ${{{k}}}$, over the time interval $0<{{{t}}}<1$ and spatial domain ${0<{x}<1}$ in problem \ref{problem4}. Again, D-J polynomials demonstrate improved accuracy compared to He’s polynomials. While additional terms can enhance the performance of D-J polynomials, the extent to which He’s polynomials can be similarly improved requires further investigation. Table \ref{table7} compares the solutions from He’s and D-J polynomials with other established techniques within the time range $0<{{{t}}}<1$ and spatial domain ${0<{{{x}}}<1}$ in problem \ref{problem4}. Finally, Table \ref{tablen6} presents a comparison of solutions and errors between He’s and D-J polynomials, New Iteration Transform Method (NITM), and q-homotopy analysis method (q-HAM) at various values of ${{{t}}}$ and ${{{x}}}$ in problem \ref{problem5}. The findings reveal that both He’s and D-J polynomials yield more accurate solutions than those provided by Khalil et al. \cite{costrashid}.

The results indicate that the proposed approach provides more accurate solutions than those documented in different articles for identical parameters and iterations. Moreover, D-J polynomials, particularly with transformation, outperform He’s polynomials, Generalized Differential Transform Method (GDTM), Adomian Decomposition Method (ADM), New Iteration Transform Method (NITM), q-homotopy analysis method (q-HAM), Optimal Auxiliary Function Method (OAFM) and Variational Iteration Method (VIM), reinforcing their effectiveness in solving the problems addressed in this study. Overall, the results across the various tables consistently show that D-J polynomials tend to outperform He’s polynomials in terms of accuracy, particularly in nonlinear scenarios.

\subsection{Computational Cost}
Our approximation is based on using two sets of polynomials to manage non-linear terms. This approach also relies on first approximations, which are derived from the initial conditions and source terms. In our case, there is no source term, but we do have a hard initial condition given by $\frac{1}{(1+\exp(x))^2}$ for problem \ref{problem3}, and $-\frac{2 e^{c{{{x}}}}c^2}{(1+ e^{{{{x}}}{{{{c}}}}})^2}$  for problem \ref{problem5} which provides the first approximation. According to the scheme, this first approximation is frequently used at every step.

Now, let's examine the polynomials in detail. We define the following polynomials:

\begin{equation*}
\begin{split}
\mathbf{H}_{0} &:= u_{0}^2, \ \ \   \ \ \   \ \ \   \ \ \   \ \ \   \ \ \   \ \ \   \ \ \   \ \ \   \ \ \  \mathbf{J}_{0} := u_{0}^2,\\
\mathbf{H}_{1} &:= 2u_{0}u_{1}, \ \ \  \ \ \   \ \ \   \ \ \   \ \ \   \ \ \   \ \ \   \ \ \   \ \   \mathbf{J}_{1} := u_{1}(2u_{0} + u_{1}),\\
\mathbf{H}_{2} &:= 2u_{0}u_{2} + u_{1}^2,\ \ \   \ \ \   \ \ \   \ \ \   \ \ \   \ \ \   \          \mathbf{J}_{2} := u_{2}(2u_{0} + 2u_{1} + u_{2}),\\
\mathbf{H}_{3} &:= 2u_{0}u_{3} + 2u_{1}u_{2}, \ \ \  \ \ \   \ \ \   \ \ \ \ \ \          \mathbf{J}_{3} := u_{3}(2u_{0} + 2u_{1} + 2u_{2} + u_{3}),\\
\mathbf{H}_{4} &:= 2u_{0}u_{4} + 2u_{1}u_{3} + u_{2}^2, \ \ \  \ \ \   \  \       \mathbf{J}_{4} := u_{4}(2u_{0} + 2u_{1} + 2u_{2} + 2u_{3} + u_{4}),\\
\mathbf{H}_{5} &:= 2u_{0}u_{5} + 2u_{1}u_{4} + 2u_{2}u_{3}, \ \ \  \       \mathbf{J}_{5} := u_{5}(2u_{0} + 2u_{1} + 2u_{2} + 2u_{3} + 2u_{4} + u_{5}).\\
\end{split}
\end{equation*}
Here, $\mathbf{H}$ represents He's polynomials and $\mathbf{J}$ represents the D-J polynomials.

It is obvious that the first two iterations of these polynomials are identical. However, from the third iteration onward, the D-J polynomials include additional terms. For example, while $\mathbf{H}_1 = 2u_0 u_1$, the corresponding D-J polynomial is $\mathbf{J}_1 = u_{1}(2u_{0} + u_{1})$. Similarly, $\mathbf{H}_{5} := 2u_{0}u_{5} + 2u_{1}u_{4} + 2u_{2}u_{3}$, the corresponding D-J polynomial is $\mathbf{J}_{5} := u_{5}(2u_{0} + 2u_{1} + 2u_{2} + 2u_{3} + 2u_{4} + u_{5})$. This difference implies that the calculation for D-J polynomials requires more computational resources in terms of memory and time.

Let's take a closer look at the computational challenges posed by problems \ref{problem3}, and \ref{problem5}.
In problem \ref{problem3}, we have a diffusion term that involves double derivatives in every approximation. This complexity is also seen in the D-J polynomials, which require more terms to be calculated. As a result, the computer takes longer to process these calculations due to the increased computational workload.
Similarly, problem \ref{problem5} presents even more complexity with initial conditions that involve third-order derivatives. While obtaining results using D-J polynomials without these complex conditions is possible, including a third-order derivative at every step significantly complicates the calculations. This complexity, combined with the intricate initial conditions, leads to a substantial computational cost, making it challenging for software to handle efficiently.
The computational cost rises when initial conditions are more complex, as seen in problems \ref{problem3}, and 5, and involve higher-order derivatives. Although all five problems contain non-linear terms, the initial conditions in problems \ref{problem1}, \ref{problem2}, and \ref{problem4} are relatively simpler compared to those in problems \ref{problem3} and \ref{problem5}. The added complexity in problems \ref{problem3}, and \ref{problem5} necessitates higher computational costs, highlighting the increased difficulty and resource demands for these specific problems.
We have observed that using D-J polynomials demands significantly more memory and processing time compared to He's Polynomials. Although D-J polynomials demonstrate high accuracy in our study, their computational expense is considerable. They sometimes become impractical to compute due to the extensive resources required.
\begin{table}[H]
    \centering
    \caption{\textbf{Comparative Analysis of Computational Costs for the First Five Iteration Terms:} This table illustrates the differences in memory and time requirements between He's and D-J polynomials. Our findings indicate that D-J polynomials demand significantly more computational resources compared to He's polynomials.}
    \label{table8}
    \resizebox{\textwidth}{!}{
        \begin{tabular}{cccccccccccc}
            \toprule
            & \multicolumn{2}{c}{\textbf{\underline{Memory Used in (MB)}}} & \multicolumn{2}{c}{\textbf{\underline{Time taken in (s)}}} \\
            \textbf{\rotatebox[origin=c]{90}{Problems}} & \textbf{\thead{He's \\ Polynomials}} & {\thead{\textbf{D-J} \textbf{Polynomials}}} & \textbf{\thead{He's \\ Polynomials}} & \textbf{\thead{D-J \\ Polynomials}} \\
            \midrule
            $1$ & $58.19$ & $70.18$ & $3.03$ & $3.18$ \\
            $2$ & $38.00$ & $38.00$ & $2.09$ & $2.09$ \\
            $3$ & $22.18$ & \textcolor[rgb]{1.00,0.00,0.50}{$82.61$} & $1.03$ & \textcolor[rgb]{1.00,0.00,0.50}{$30.06$} \\
            $4$ & $40.18$ & $40.18$ & $2.28$ & $2.59$ \\
            $5$ & $72.18$ & \textcolor[rgb]{1.00,0.00,0.00}{$597.21^\nearrow$} & $2.64$ & \textcolor[rgb]{1.00,0.00,0.00}{$8735.39^\nearrow$} \\
            \bottomrule
        \end{tabular}
    }
\end{table}

\begin{figure}[H]
    \centering
    \includegraphics[width=8.3cm]{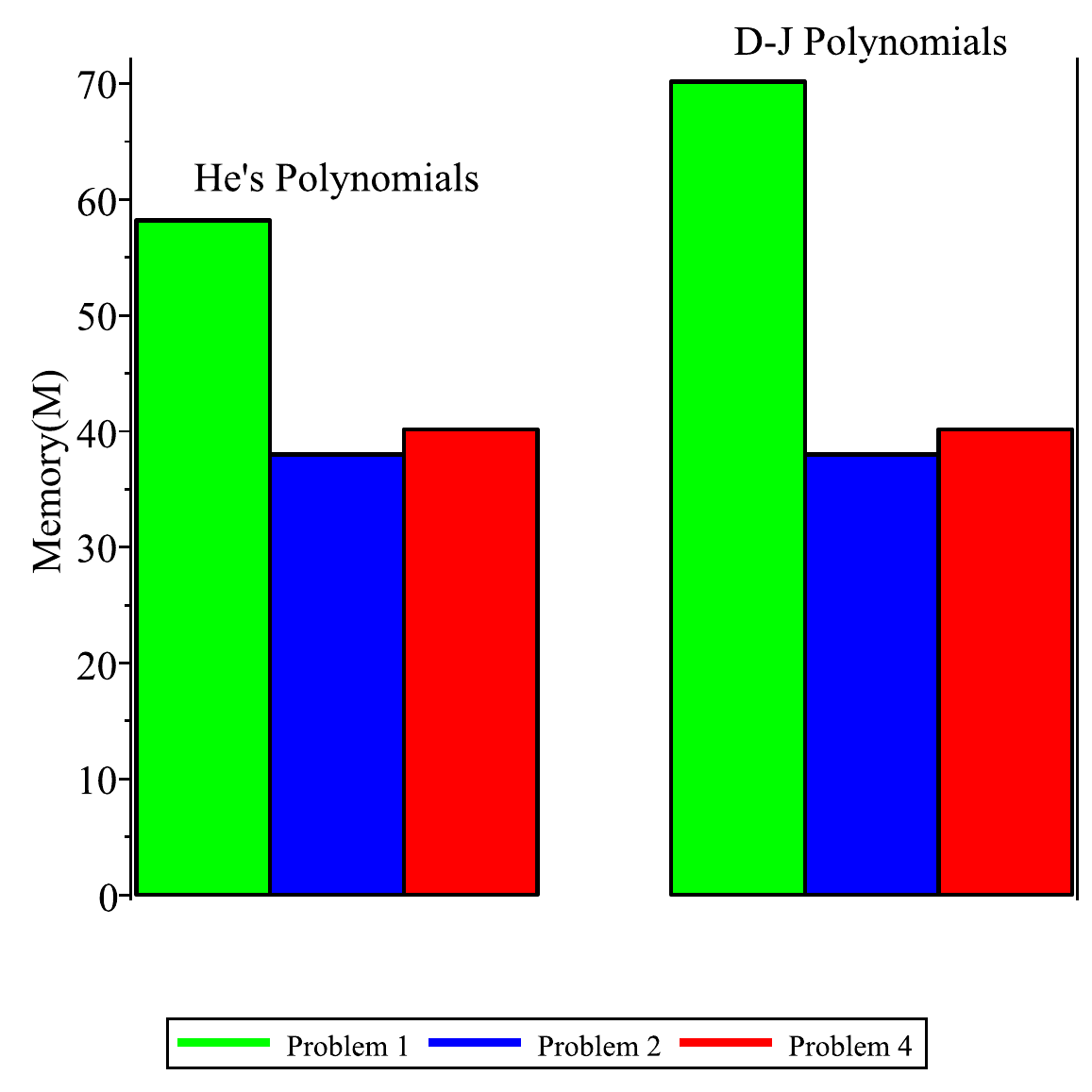}%
    %\llap{\raisebox{0.49cm}{%  move next graphics to top right corner
%      \includegraphics[width=2.45cm]{M_1-eps-converted-to}\hspace{0.1em}
%      \includegraphics[width=2.45cm]{M_2-eps-converted-to}\hspace{0.1em}%
%    }}\hspace{0.3em}
   %  \includegraphics[width=5.3cm]{5a-eps-converted-to}\hspace{0.3em}%
%     \includegraphics[width=5.3cm]{4a-eps-converted-to}\hspace{0.3em}
  \includegraphics[width=8.3cm]{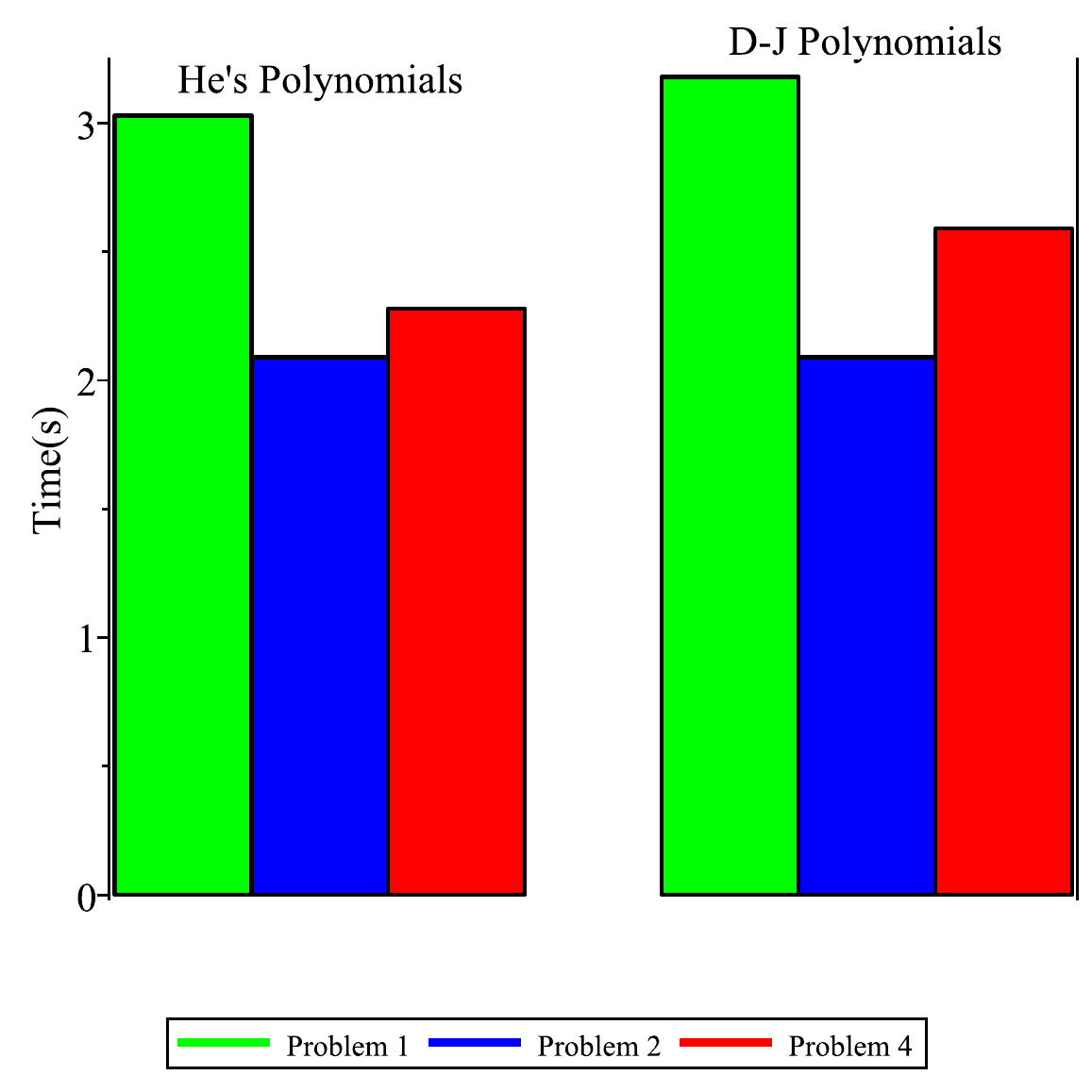}%
    \caption{ \textbf{ The computational cost of He's and D-J polynomials for 1st five terms of iterations.} }\label{fig5}
  \end{figure}

\begin{figure}[H]
    \centering
  \includegraphics[width=4.0cm]{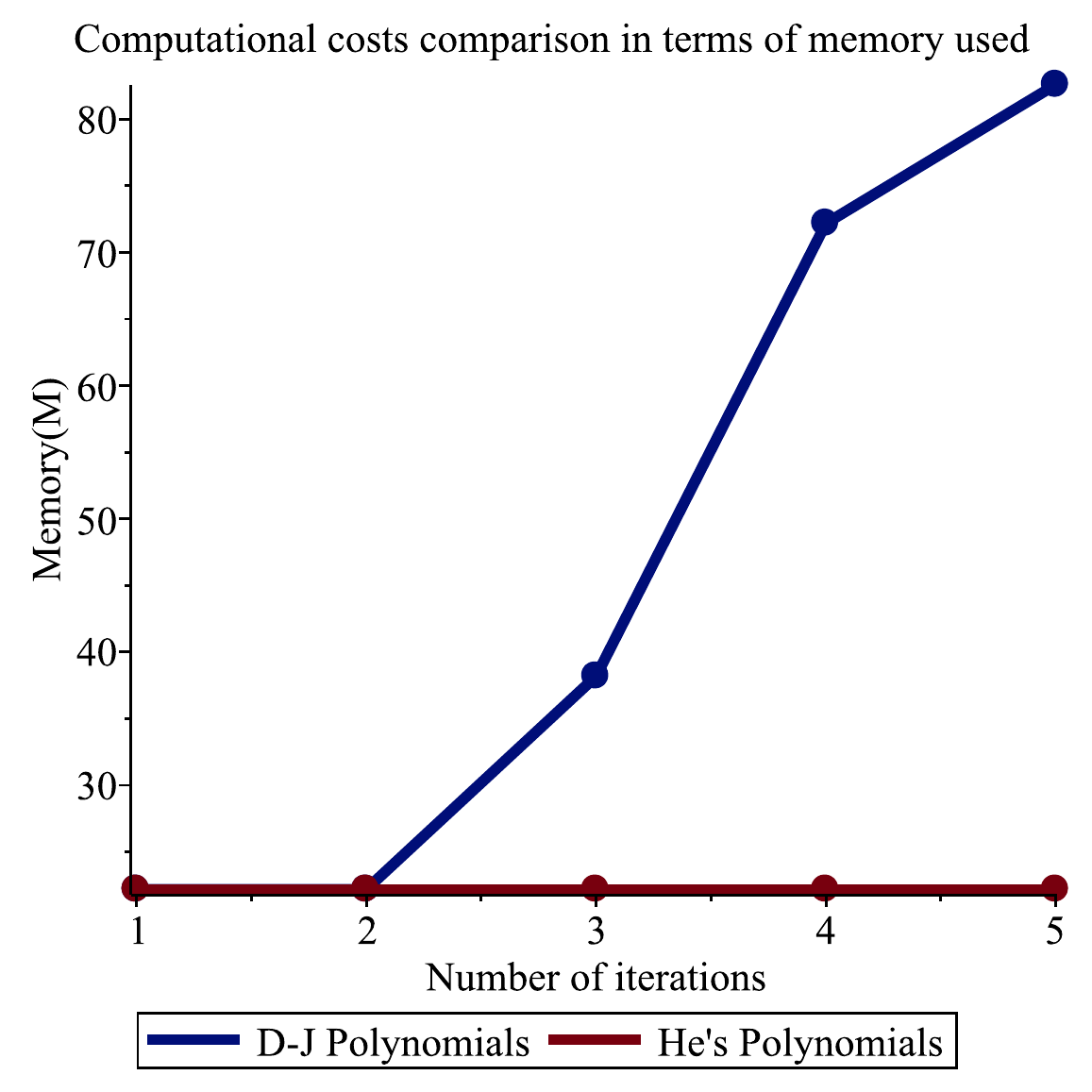}
     \includegraphics[width=4.0cm]{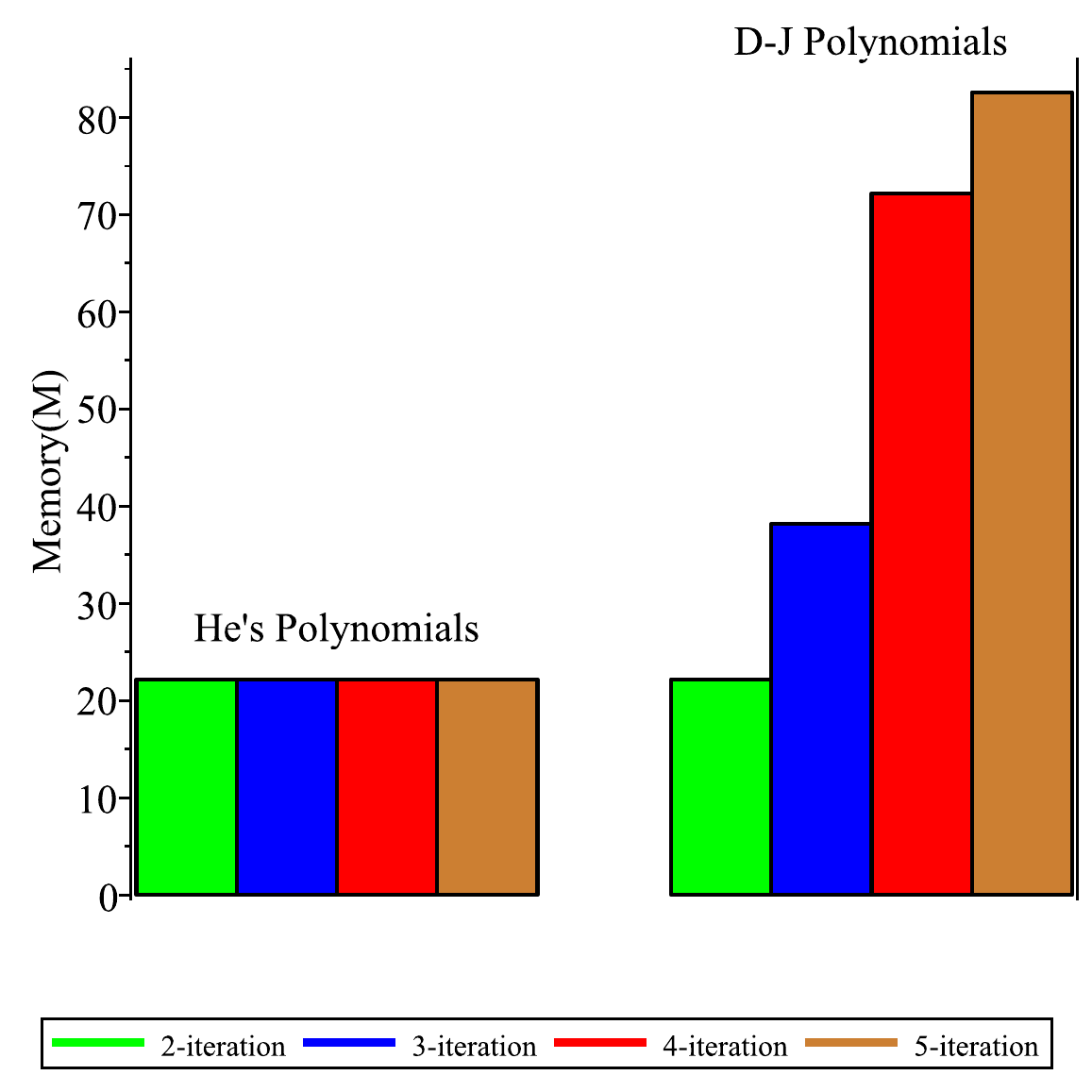}%
   \includegraphics[width=4.0cm]{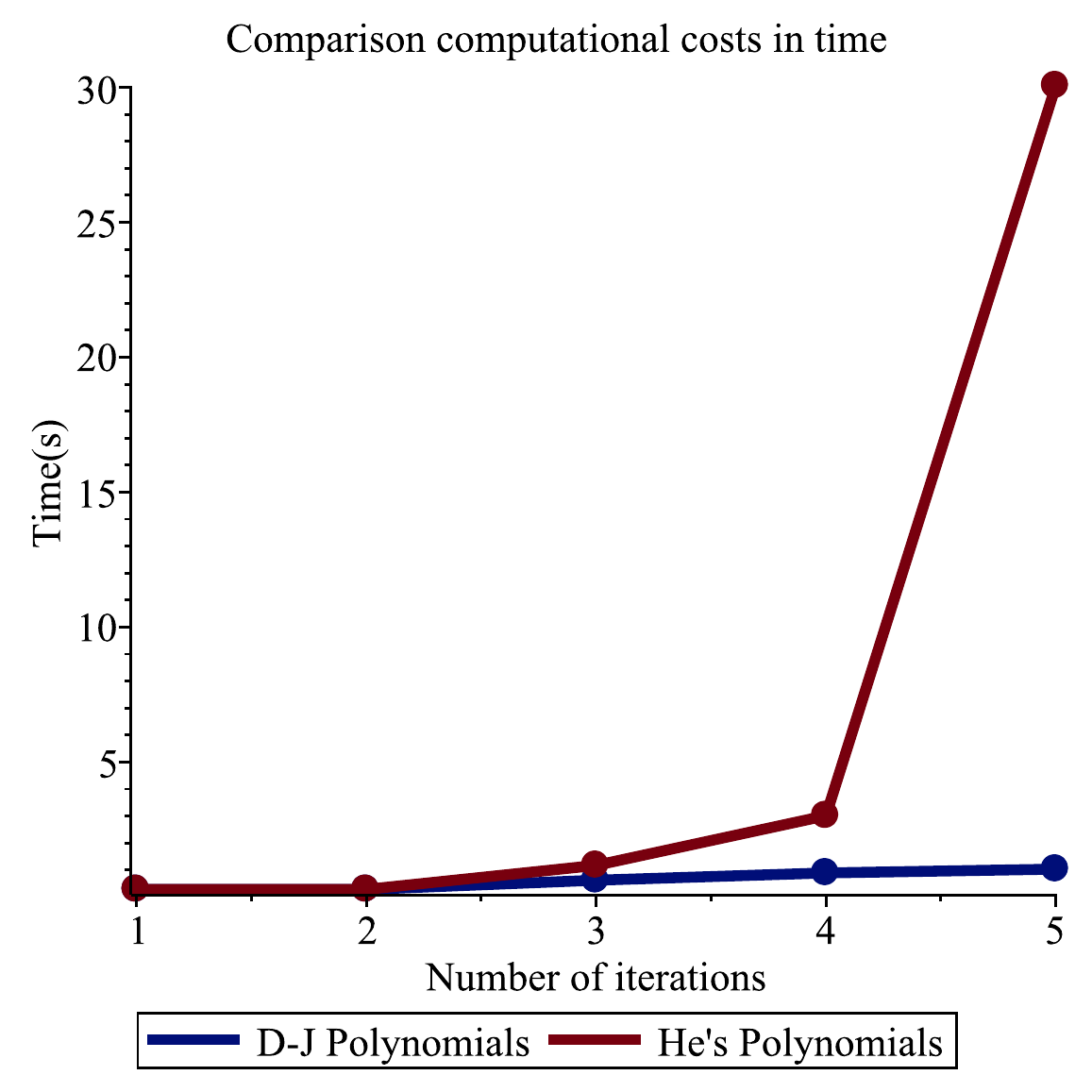}
  \includegraphics[width=4.0cm]{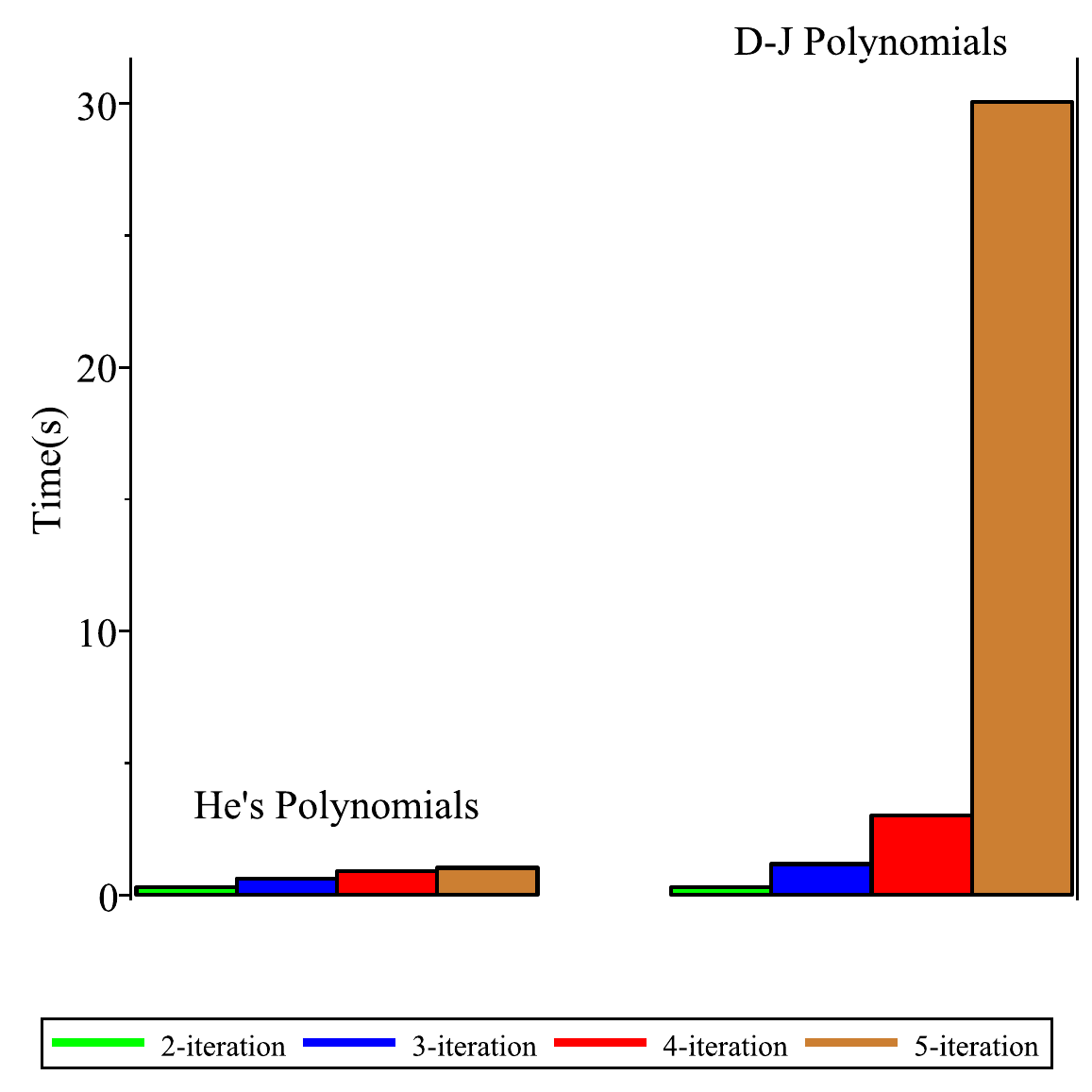}\\%%
  \caption{ \textbf{ The computational cost of He's and D-J polynomials for problem \ref{problem3}.} }\label{fig25}
  \end{figure}

\begin{figure}[H]
    \centering
  \includegraphics[width=4.0cm]{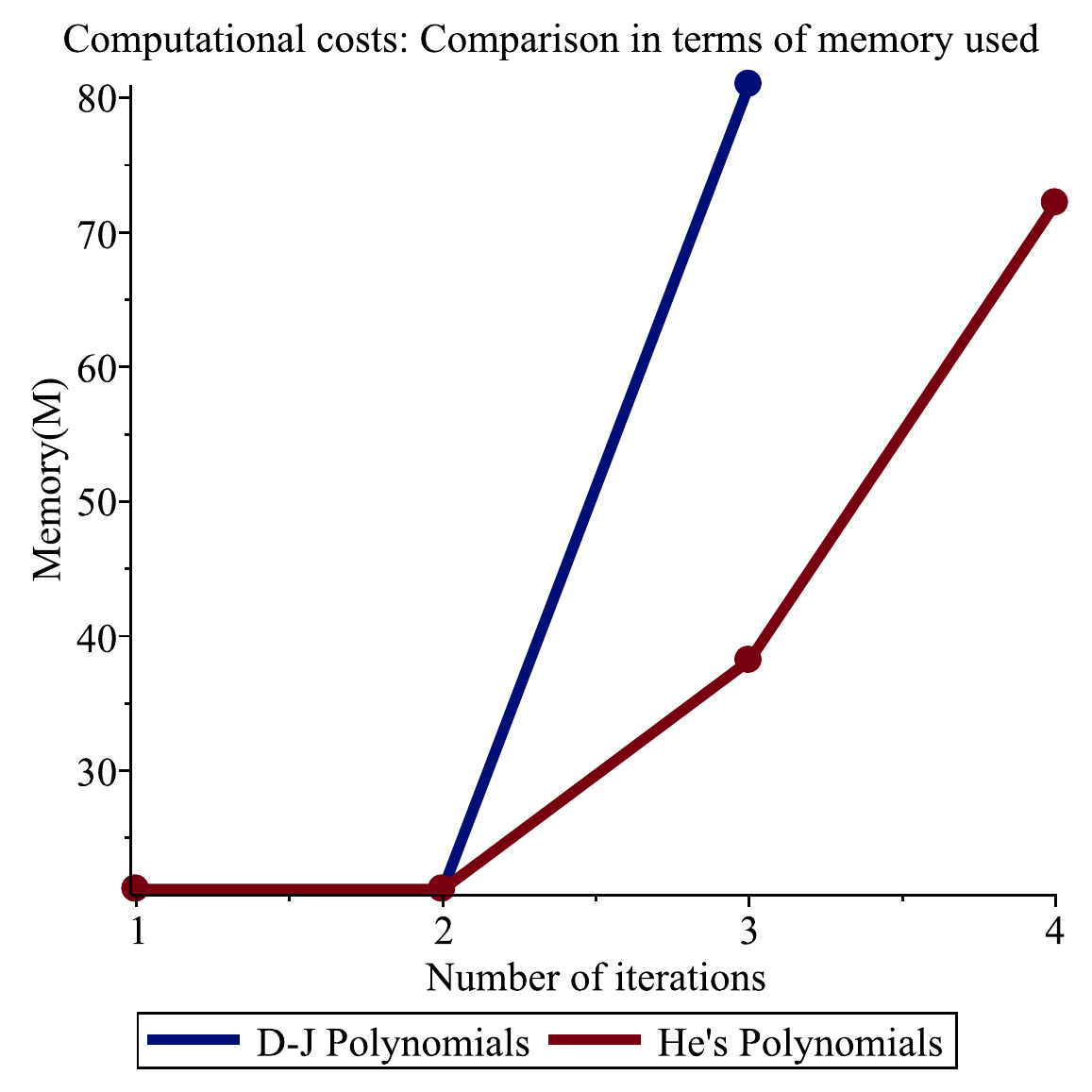}
     \includegraphics[width=4.0cm]{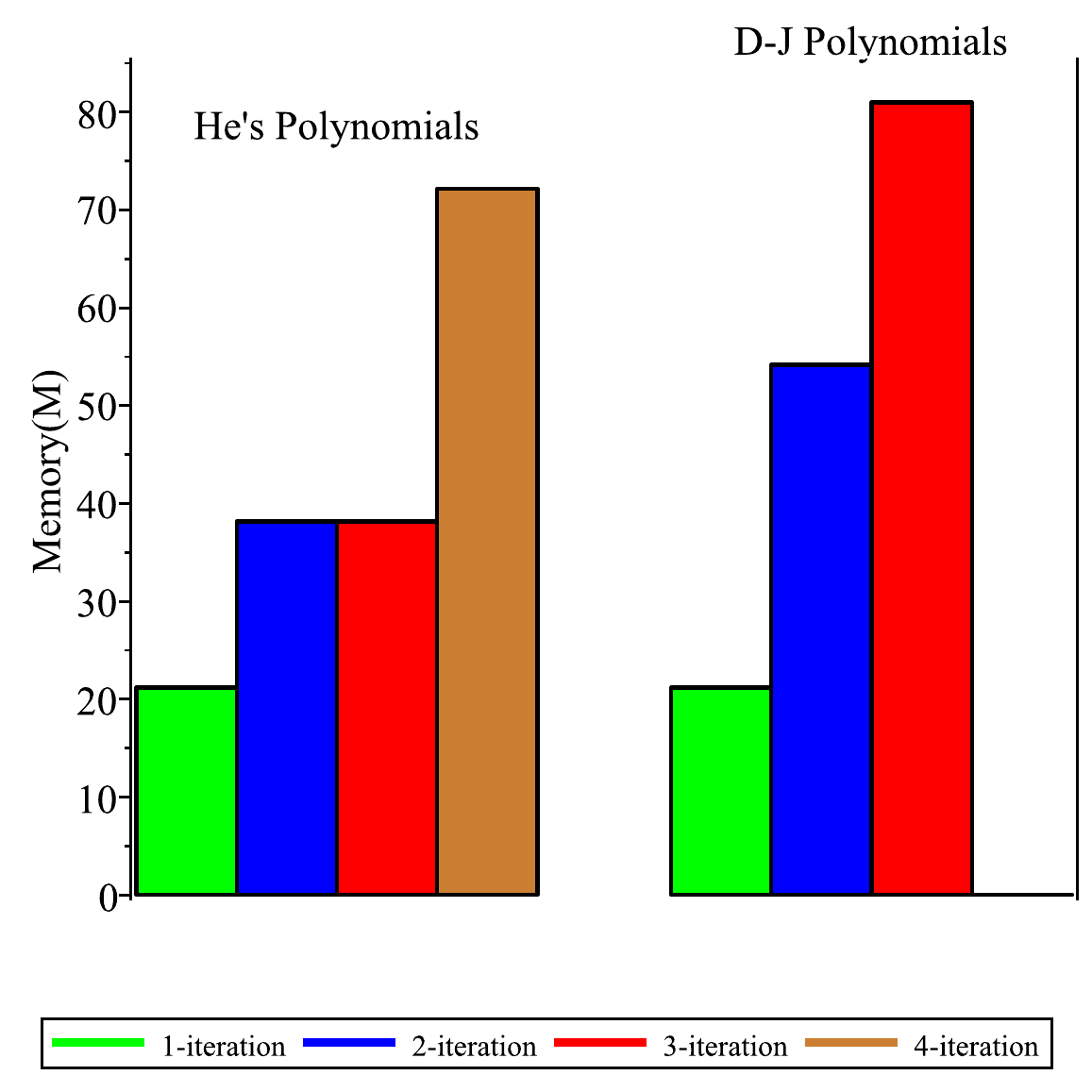}
      \llap{\raisebox{0.89cm}{   %move next graphics to top right corner
       \includegraphics[width=2.45cm,  angle =90 ]{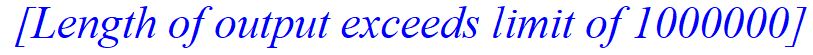}\hspace{0.6em}%
     }}%
   \includegraphics[width=4.0cm]{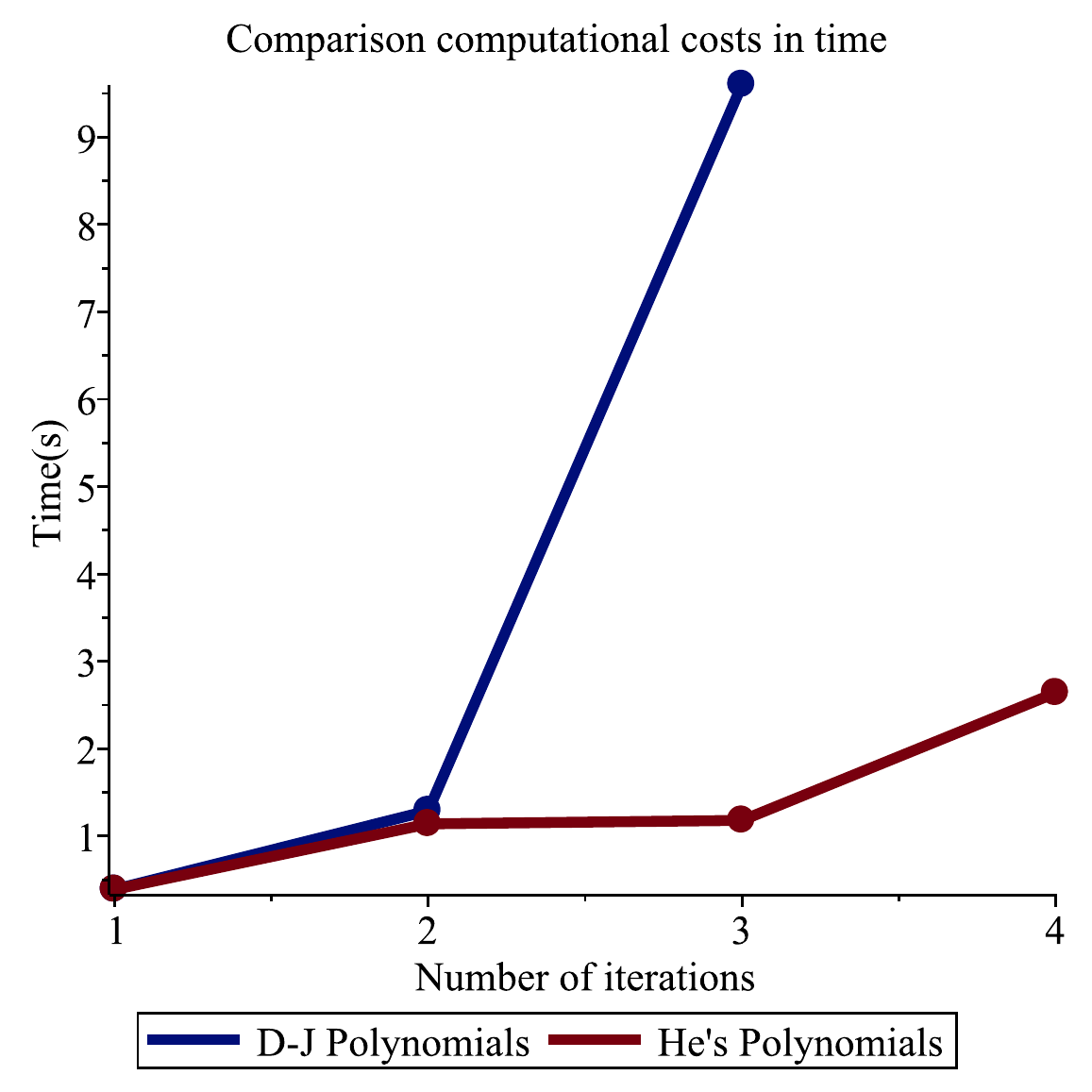}
  \includegraphics[width=4.0cm]{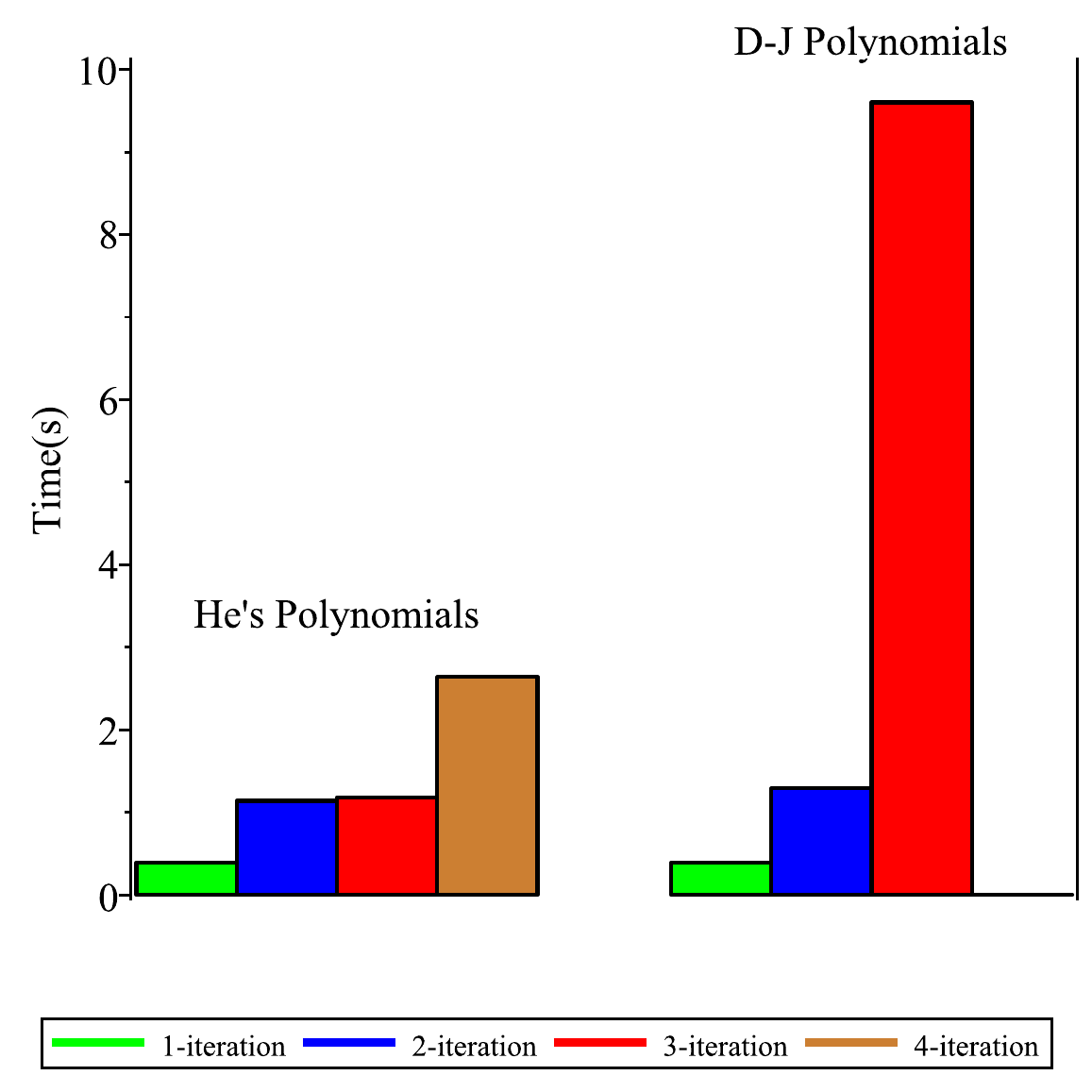}
  \llap{\raisebox{0.89cm}{   %move next graphics to top right corner
       \includegraphics[width=2.45cm,  angle =90 ]{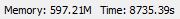}\hspace{0.6em}%
     } }\hspace{0.1em}\\%%
  \caption{\textbf{The computational cost of He's and D-J polynomials for problem \ref{problem5}.}
We used a big amount of memory and spent much more time, but we couldn't execute the fourth iteration for the Jafari polynomial in this particular problem. Additionally, our operating system became overheated and started emitting unusual noises. Maple provided feedback indicating that the length of the output exceeds the limit of 1,000,000, causing the computations to stop. This clearly shows that more complex initial data and higher derivatives overwhelm the system, preventing further iterations.
However, we successfully obtained the fourth and fifth iterations in Problem \ref{problem3}, despite a significant computational jump between both polynomials. In this case, it seems we cannot achieve the same result, which caused a huge computational cost.}\label{fig26}
  \end{figure}

Table \ref{table8} and Figure  \ref{fig5} provide a comparative analysis of the computational costs, including time and memory usage, for D-J polynomials and He's polynomials. In this project, we observed that D-J polynomials require significantly more computational resources than He's polynomials. To illustrate this with specific data, we examined Problems \ref{problem3} and \ref{problem5} in more detail (see figures \ref{fig25} and  \ref{fig26}). For this particular problem, He's polynomials consumed $22.18$ megabytes (MB) of memory and completed the computation in $2.28$ seconds. In contrast, D-J polynomials required $82.61$ MB of memory and took substantially longer, completing the computation in $30.06$ seconds. For problem \ref{problem5}, the 4th iteration term of He's polynomials required 72.18 MB of memory and completed the computation in 2.64 seconds. In contrast, the D-J polynomials failed to produce the 4th iteration term, consuming a substantial 597.21 MB of memory and taking 8735.39 seconds, with the computation still ongoing. The reason for this disparity in computational cost is that D-J polynomials provide more accurate solutions. However, this increased accuracy comes at the expense of needing more advanced and high-memory computing systems and software. Consequently, implementing D-J polynomials necessitates access to advanced hardware and specialized software capable of handling their increased memory and processing requirements.


In summary, we presented a comparative analysis of He's polynomials and D-J polynomials for solving non-linear fractional partial differential equations.
The analysis showed that D-J polynomials give more accurate results but need much more memory and processing time than He's polynomials. This higher computational demand advanced hardware and special software, like GPU systems or machine learning tools. Our findings show that D-J polynomials demand more computational costs for computing but are very accurate for solving non-linear fractional PDEs. Future work could explore further optimization techniques to reduce D-J polynomials' computational limitations, making them more accessible for practical applications in different branches of applied sciences.
\section{Conclusion}\label{conclusion}
This article presents a comprehensive comparison of He's and Daftardar-Jafari (D-J) polynomials in solving fractional partial differential equations (FPDEs) using the Iterative Laplace Transform Method (ILTM). Our analysis reveals that ILTM solutions maintain consistency across different polynomial implementations, providing a reliable framework for comparative evaluation. The findings indicate that D-J polynomials exhibit superior accuracy compared to He's polynomials and other methods, such as the New Iteration Method (NIM) and Optimal Auxiliary Function Method (OAFM). The absolute error values consistently show that D-J polynomials align more closely with exact solutions across various problems, underscoring their effectiveness in minimizing deviations from true solutions. This superior performance makes D-J polynomials a robust choice for achieving higher precision in computational tasks. However, the potential for enhancing He's polynomials through additional terms in the series solution remains a promising avenue for future research. Further validation and investigation are necessary to determine whether He's polynomials can achieve similar improvements in accuracy. The graphical and tabular results validate the high accuracy achieved by ILTM, confirming its effectiveness in solving FPDEs. This study provides valuable insights into the selection of appropriate polynomials for various FPDEs and their systems, thereby broadening the applicability of these methods in future research. While D-J polynomials offer a more accurate approach for solving non-linear fractional PDEs, their higher computational cost should be considered. He's polynomials, on the other hand, offer a more resource-efficient alternative, with significant potential for improvement through further research.

In conclusion, this study highlights the trade-offs between accuracy and computational cost, guiding the selection of polynomial methods based on specific computational constraints and accuracy requirements. The insights gained from this research can assist in making informed decisions when choosing polynomial methods for solving FPDEs, contributing to the advancement of computational techniques in this field.

\subsection*{CRediT authorship contribution statement}
Both authors contributed equally to this work and agreed on the final version of it.
\subsection*{Declaration of competing interest}
				The authors declare that they have no known competing financial interests or personal relationships that could have appeared to influence the work reported in this research paper.
				\subsection*{Data availability}
			This paper does not contain any hidden data. Q. Khan developed and implemented the codes using Maple Version 2024. Data for different numbers of iterations or fractional orders or for  specific values  of ${{{t}}}$ and {{{$x$}}} will be provided upon reasonable request.
				\subsection*{Funding}
This work is partially supported by Hong Kong General Research Fund (GRF) grant project number 18300821, 18300622 and 18300424.
				\subsection*{Acknowledgment }
We are very grateful to our MIT lab mates and the technicians at The Education University of Hong Kong for providing all the necessary facilities for this research work.

\bibliographystyle{unsrt}
\bibliography{ref_cost_paper}
\end{document}